%
%
\font\hd=cmbx10 scaled\magstep1

\ifx\begin\undefined\else\global\advance\srcdepth by
1\expandafter\endinput\fi

\def\begin{}
\newcount\srcdepth
\srcdepth=1

\outer\def\bye{\global\advance\srcdepth by -1
  \ifnum\srcdepth=0
    \def\endcmd{\vfill\supereject\nopagenumbers\par\vfill\supereject\end}
  \else\def\endcmd{}\fi
  \endcmd
}



\def\initialize#1#2#3#4#5#6{
  \ifnum\srcdepth=1
  \magnification=#1
  \hsize = #2
  \vsize = #3
  \hoffset=#4
  \advance\hoffset by -\hsize
  \divide\hoffset by 2
  \advance\hoffset by -1truein
  \voffset=#5
  \advance\voffset by -\vsize
  \divide\voffset by 2
  \advance\voffset by -1truein
  \advance\voffset by #6
  \baselineskip=13pt
  \emergencystretch = 0.05\hsize
  \parskip=3pt plus1pt minus.5pt
  \fi
}

\def\print{\initialize{1095}
  {6.5truein}{8.5truein}{8.5truein}{11truein}{-.0625truein}}

\overfullrule=0pt

\newif\ifblackboardbold

\blackboardboldtrue


\font\sectionfont=cmbx12

\font\scriptit=cmti10 at 7pt
\font\scriptsl=cmsl10 at 7pt
\scriptfont\itfam=\scriptit
\scriptfont\slfam=\scriptsl


\newfam\msamfam  
\font\tenmsam=msam10
\font\sevenmsam=msam7
\font\fivemsam=msam5
\textfont\msamfam=\tenmsam
\scriptfont\msamfam=\sevenmsam
\scriptscriptfont\msamfam=\fivemsam

\newfam\msbmfam
\font\tenmsbm=msam10
\font\sevenmsbm=msam7
\font\fivemsbm=msam5
\textfont\msbmfam=\tenmsbm
\scriptfont\msbmfam=\sevenmsbm
\scriptscriptfont\msbmfam=\fivemsbm

\newfam\eufmfam  
\font\teneufm=eufm10
\font\seveneufm=eufm7
\font\fiveeufm=eufm5
\textfont\eufmfam\teneufm
\scriptfont\eufmfam\seveneufm
\scriptscriptfont\eufmfam\fiveeufm

\newfam\bboldfam
\ifblackboardbold
\font\tenbbold=msbm10
\font\sevenbbold=msbm7
\font\fivebbold=msbm5
\textfont\bboldfam=\tenbbold
\scriptfont\bboldfam=\sevenbbold
\scriptscriptfont\bboldfam=\fivebbold
\def\bbold{\fam\bboldfam\tenbbold}
\else
\def\bbold{\bf}
\fi
\newcount\amsfamcount 
\newcount\classcount   
\newcount\positioncount
\newcount\codecount
\newcount\n             
\def\newsymbol#1#2#3#4#5{               
\n="#2                                  
\ifnum\n=1 \amsfamcount=\msamfam\else   
\ifnum\n=2 \amsfamcount=\msbmfam\else   
\ifnum\n=3 \amsfamcount=\eufmfam
\fi\fi\fi
\multiply\amsfamcount by "100           
\classcount="#3                 
\multiply\classcount by "1000           
\positioncount="#4#5            
\codecount=\classcount                  
\advance\codecount by \amsfamcount      
\advance\codecount by \positioncount
\mathchardef#1=\codecount}              


\font\Arm=cmr9
\font\Ai=cmmi9
\font\Asy=cmsy9
\font\Abf=cmbx9
\font\Brm=cmr7
\font\Bi=cmmi7
\font\Bsy=cmsy7
\font\Bbf=cmbx7
\font\Crm=cmr6
\font\Ci=cmmi6
\font\Csy=cmsy6
\font\Cbf=cmbx6

\ifblackboardbold
\font\Abbold=msbm10 at 9pt
\font\Bbbold=msbm7
\font\Cbbold=msbm5 at 6pt
\fi

\def\smallmath{%
\textfont0=\Arm \scriptfont0=\Brm \scriptscriptfont0=\Crm
\textfont1=\Ai \scriptfont1=\Bi \scriptscriptfont1=\Ci
\textfont2=\Asy \scriptfont2=\Bsy \scriptscriptfont2=\Csy
\textfont\bffam=\Abf \scriptfont\bffam=\Bbf \scriptscriptfont\bffam=\Cbf
\def\rm{\fam0\Arm}\def\mit{\fam1}\def\oldstyle{\fam1\Ai}%
\def\bf{\fam\bffam\Abf}%
\ifblackboardbold
\textfont\bboldfam=\Abbold
\scriptfont\bboldfam=\Bbbold
\scriptscriptfont\bboldfam=\Cbbold
\def\bbold{\fam\bboldfam\Abbold}%
\fi
\rm
}








\newlinechar=`@
\def\forwardmsg#1#2#3{\immediate\write16{@*!*!*!* forward reference should
be: @\noexpand\forward{#1}{#2}{#3}@}}
\def\nodefmsg#1{\immediate\write16{@*!*!*!* #1 is an undefined reference@}}

\def\forwardsub#1#2{\def\newref{{#2}{#1}}}

\def\forward#1#2#3{%
\expandafter\expandafter\expandafter\forwardsub\expandafter{#3}{#2}
\expandafter\ifx\csname#1\endcsname\relax\else%
\expandafter\ifx\csname#1\endcsname\newref\else%
\forwardmsg{#1}{#2}{#3}\fi\fi%
\expandafter\let\csname#1\endcsname\newref}

\def\firstarg#1{\expandafter\argone #1}\def\argone#1#2{#1}
\def\secondarg#1{\expandafter\argtwo #1}\def\argtwo#1#2{#2}

\def\ref#1{\expandafter\ifx\csname#1\endcsname\relax
  {\nodefmsg{#1}\bf`#1'}\else
  \expandafter\firstarg\csname#1\endcsname
  ~\expandafter\secondarg\csname#1\endcsname\fi}

\def\refs#1{\expandafter\ifx\csname#1\endcsname\relax
  {\nodefmsg{#1}\bf`#1'}\else
  \expandafter\firstarg\csname #1\endcsname
  s~\expandafter\secondarg\csname#1\endcsname\fi}

\def\refn#1{\expandafter\ifx\csname#1\endcsname\relax
  {\nodefmsg{#1}\bf`#1'}\else
  \expandafter\secondarg\csname #1\endcsname\fi}



\def\widow#1{\vskip 0pt plus#1\vsize\goodbreak\vskip 0pt plus-#1\vsize}



\def\marginlabel#1{}

\def\showlabelsabove{
\font\labelfont=cmss10 at 6pt
\def\marginlabel##1{\rlap{\smash{\raise 10pt\hbox{\labelfont##1}}}}
}

\newcount\seccount
\newcount\proccount
\seccount=0
\proccount=0

\def\stdskip{\vskip 9pt plus3pt minus 3pt}
\def\stdbreak{\par\removelastskip\penalty-100\stdskip}

\def\proof{\stdbreak\noindent{\sl Proof }}

\def\qed{\vrule height 1.2ex width .9ex depth .1ex}

\def\Box{
  \ifmmode\eqno\qed
  \else\ifvmode\removelastskip\line{\hfil\qed}
  \else\unskip\quad\hskip-\hsize
    \hbox{}\hskip\hsize minus 1em\qed\par
  \fi\stdbreak\fi}

\def\references{
  \removelastskip
  \widow{.05}
  \vskip 24pt plus 6pt minus 6 pt
  \leftline{\sectionfont References}
  \nobreak\stdskip\noindent}

\def\ifempty#1#2\endB{\ifx#1\endA}
\def\makeref#1#2#3{\ifempty#1\endA\endB\else\forward{#1}{#2}{#3}\fi}

\outer\def\section#1 #2\par{
  \removelastskip
  \global\advance\seccount by 1
  \global\proccount=0\relax
                \edef\numtoks{\number\seccount}
  \makeref{#1}{Section}{\numtoks}
  \widow{.05}
  \vskip 24pt plus 6pt minus 6 pt
  \message{#2}
  \leftline{\marginlabel{#1}\sectionfont\numtoks\quad #2}
  \nobreak\stdskip}

\def\proclamation#1#2{
  \outer\expandafter\def\csname#1\endcsname##1 ##2\par{
  \stdbreak
  \global\advance\proccount by 1
  \edef\numtoks{\number\seccount.\number\proccount}
  \makeref{##1}{#2}{\numtoks}
  \noindent{\marginlabel{##1}\bf #2 \numtoks\enspace}
  {\sl##2\par}
  \stdbreak}}

\def\othernumbered#1#2{
  \outer\expandafter\def\csname#1\endcsname##1{
  \stdbreak
  \global\advance\proccount by 1
  \edef\numtoks{\number\seccount.\number\proccount}
  \makeref{##1}{#2}{\numtoks}
  \noindent{\marginlabel{##1}\bf #2 \numtoks\enspace}}}

\proclamation{definition}{Definition}
\proclamation{lemma}{Lemma}
\proclamation{proposition}{Proposition}
\proclamation{theorem}{Theorem}
\proclamation{corollary}{Corollary}
\proclamation{conjecture}{Conjecture}

\othernumbered{example}{Example}
\othernumbered{remark}{Remark}
\othernumbered{construction}{Construction}


\def\figure#1{
 \global\advance\figcount by 1
 \goodbreak
 \midinsert#1\smallskip
 \centerline{Figure~\number\figcount}
 \endinsert}

\def\capfigure#1#2{
 \global\advance\figcount by 1
 \goodbreak
 \midinsert#1\smallskip
 \vbox{\small\noindent {\bf Figure~\number\figcount:} #2}
 \endinsert}

\def\capfigurepair#1#2#3#4{
 \goodbreak
 \midinsert
 #1\smallskip
 \global\advance\figcount by 1
 \vbox{\small\noindent {\bf Figure~\number\figcount:} #2}
 \vskip 12pt
 #3\smallskip
 \global\advance\figcount by 1
 \vbox{\small\noindent {\bf Figure~\number\figcount:} #4}
 \endinsert}


\def\baretable#1#2{
\vbox{\offinterlineskip\halign{
 \strut\kern #1\hfil##\kern #1
 &&\kern #1\hfil##\kern #1\cr
 #2
}}}

\def\gridtablesub#1#2#3{
\vbox{\offinterlineskip\halign{
 \strut\vrule\kern #1\hfil##\hfil\kern #2\vrule
 &&\kern #1\hfil##\kern #2\vrule\cr
 \noalign{\hrule}
 #3
 \noalign{\hrule}
}}}





\input epsf

\newcount\figcount
\figcount=0
\newbox\drawing
\newcount\drawbp
\newdimen\drawx
\newdimen\drawy
\newdimen\ngap
\newdimen\sgap
\newdimen\wgap
\newdimen\egap

\def\drawbox#1#2#3{\vbox{
  \setbox\drawing=\vbox{\offinterlineskip\epsfbox{#2.eps}\kern 0pt}
  \drawbp=\epsfurx
  \advance\drawbp by-\epsfllx\relax
  \multiply\drawbp by #1
  \divide\drawbp by 100
  \drawx=\drawbp truebp
  \ifdim\drawx>\hsize\drawx=\hsize\fi
  \epsfxsize=\drawx
  \setbox\drawing=\vbox{\offinterlineskip\epsfbox{#2.eps}\kern 0pt}
  \drawx=\wd\drawing
  \drawy=\ht\drawing
  \ngap=0pt \sgap=0pt \wgap=0pt \egap=0pt
  \setbox0=\vbox{\offinterlineskip
    \box\drawing \ifgridlines\drawgrid\drawx\drawy\fi #3}
  \kern\ngap\hbox{\kern\wgap\box0\kern\egap}\kern\sgap}}

\def\draw#1#2#3{
  \setbox\drawing=\drawbox{#1}{#2}{#3}
  \global\advance\figcount by 1
  \goodbreak
  \midinsert
  \centerline{\ifgridlines\boxgrid\drawing\fi\box\drawing}
  \smallskip
  \vbox{\offinterlineskip
    \centerline{Figure~\number\figcount}
    \smash{\marginlabel{#2}}}
  \endinsert}

\def\capdraw#1#2#3#4{
  \setbox\drawing=\drawbox{#1}{#2}{#3}
  \global\advance\figcount by 1
  \goodbreak
  \midinsert
  \centerline{\ifgridlines\boxgrid\drawing\fi\box\drawing}
  \smallskip
  \vbox{\offinterlineskip
    \vskip 4pt
    \vbox{\lineskip=3pt\small\noindent {\bf Figure~\number\figcount:} #4}
    \smash{\marginlabel{#2}}}
  \endinsert}

\def\capdrawpair#1#2#3#4#5#6#7#8{
  \goodbreak
  \midinsert
  \setbox\drawing=\drawbox{#1}{#2}{#3}
  \global\advance\figcount by 1
  \centerline{\ifgridlines\boxgrid\drawing\fi\box\drawing}
  \smallskip
  \vbox{\offinterlineskip
    \vskip 4pt
    \vbox{\lineskip=3pt\small\noindent {\bf Figure~\number\figcount:} #4}
    \smash{\marginlabel{#2}}}
  \vskip 12pt
  \setbox\drawing=\drawbox{#5}{#6}{#7}
  \global\advance\figcount by 1
  \centerline{\ifgridlines\boxgrid\drawing\fi\box\drawing}
  \smallskip
  \vbox{\offinterlineskip
    \vskip 4pt
    \vbox{\lineskip=3pt\small\noindent {\bf Figure~\number\figcount:} #8}
    \smash{\marginlabel{#6}}}
  \endinsert}

\def\drawnoname#1#2#3{
  \setbox\drawing=\drawbox{#1}{#2}{#3}
  \global\advance\figcount by 1
  \goodbreak
  \midinsert
  \centerline{\ifgridlines\boxgrid\drawing\fi\box\drawing}
  \smallskip
  \endinsert}

\def\nextfigtoks{%
  \advance\figcount by 1%
  \edef\numtoks{\number\figcount}%
  \advance\figcount by -1}

\newif\ifgridlines
\newbox\figtbox
\newbox\figgbox
\newdimen\figtx
\newdimen\figty

\newdimen\bwd
\bwd=2sp 

\def\hline#1{\vbox{\smash{\hbox to #1{\leaders\hrule height \bwd\hfil}}}}

\def\vline#1{\hbox to 0pt{%
  \hss\vbox to #1{\leaders\vrule width \bwd\vfil}\hss}}

\def\clap#1{\hbox to 0pt{\hss#1\hss}}
\def\vclap#1{\vbox to 0pt{\offinterlineskip\vss#1\vss}}

\def\hstutter#1#2{\hbox{%
  \setbox0=\hbox{#1}%
  \hbox to #2\wd0{\leaders\box0\hfil}}}

\def\vstutter#1#2{\vbox{
  \setbox0=\vbox{\offinterlineskip #1}
  \dp0=0pt
  \vbox to #2\ht0{\leaders\box0\vfil}}}

\def\crosshairs#1#2{
  \dimen1=.002\drawx
  \dimen2=.002\drawy
  \ifdim\dimen1<\dimen2\dimen3\dimen1\else\dimen3\dimen2\fi
  \setbox1=\vclap{\vline{2\dimen3}}
  \setbox2=\clap{\hline{2\dimen3}}
  \setbox3=\hstutter{\kern\dimen1\box1}{4}
  \setbox4=\vstutter{\kern\dimen2\box2}{4}
  \setbox1=\vclap{\vline{4\dimen3}}
  \setbox2=\clap{\hline{4\dimen3}}
  \setbox5=\clap{\copy1\hstutter{\box3\kern\dimen1\box1}{6}}
  \setbox6=\vclap{\copy2\vstutter{\box4\kern\dimen2\box2}{6}}
  \setbox1=\vbox{\offinterlineskip\box5\box6}
  \smash{\vbox to #2{\hbox to #1{\hss\box1}\vss}}}

\def\boxgrid#1{\rlap{\vbox{\offinterlineskip
  \setbox0=\hline{\wd#1}
  \setbox1=\vline{\ht#1}
  \smash{\vbox to \ht#1{\offinterlineskip\copy0\vfil\box0}}
  \smash{\vbox{\hbox to \wd#1{\copy1\hfil\box1}}}}}}

\def\drawgrid#1#2{\vbox{\offinterlineskip
  \dimen0=\drawx
  \dimen1=\drawy
  \divide\dimen0 by 10
  \divide\dimen1 by 10
  \setbox0=\hline\drawx
  \setbox1=\vline\drawy
  \smash{\vbox{\offinterlineskip
    \copy0\vstutter{\kern\dimen1\box0}{10}}}
  \smash{\hbox{\copy1\hstutter{\kern\dimen0\box1}{10}}}}}

\def\figtext#1#2#3#4#5{
  \setbox\figtbox=\vbox{\hbox{#5}\kern 0pt}
  \figtx=-#3\wd\figtbox \figty=-#4\ht\figtbox
  \advance\figtx by #1\drawx \advance\figty by #2\drawy
  \dimen0=\figtx \advance\dimen0 by\wd\figtbox \advance\dimen0 by-\drawx
  \ifdim\dimen0>\egap\global\egap=\dimen0\fi
  \dimen0=\figty \advance\dimen0 by\ht\figtbox \advance\dimen0 by-\drawy
  \ifdim\dimen0>\ngap\global\ngap=\dimen0\fi
  \dimen0=-\figtx
  \ifdim\dimen0>\wgap\global\wgap=\dimen0\fi
  \dimen0=-\figty
  \ifdim\dimen0>\sgap\global\sgap=\dimen0\fi
  \smash{\rlap{\vbox{\offinterlineskip
    \hbox{\hbox to \figtx{}\ifgridlines\boxgrid\figtbox\fi\box\figtbox}
    \vbox to \figty{}
    \ifgridlines\crosshairs{#1\drawx}{#2\drawy}\fi
    \kern 0pt}}}}

\def\swtext#1#2#3{\figtext{#1}{#2}00{#3}}


\def\hpad#1#2#3{\hbox{\kern #1\hbox{#3}\kern #2}}
\def\vpad#1#2#3{\setbox0=\hbox{#3}\vbox{\kern #1\box0\kern #2}}

\def\wpad#1#2{\hpad{#1}{0pt}{#2}}

\def\spad#1#2{\vpad{0pt}{#1}{#2}}

\def\swpad#1#2#3{\spad{#1}{\wpad{#2}{#3}}}




\def\stack#1#2#3{\vbox{\offinterlineskip
  \setbox2=\hbox{#2}
  \setbox3=\hbox{#3}
  \dimen0=\ifdim\wd2>\wd3\wd2\else\wd3\fi
  \hbox to \dimen0{\hss\box2\hss}
  \kern #1
  \hbox to \dimen0{\hss\box3\hss}}}


\def\hexp#1{%
  \setbox0=\hbox{${}^{#1}$}%
  \hbox to .5\wd0{\box0\hss}}

\def\hsub#1{%
  \setbox0=\hbox{${}_{#1}$}%
  \hbox to .5\wd0{\box0\hss}}



\def\bmatrix#1#2{{\left[\vcenter{\halign
  {&\kern#1\hfil$##\mathstrut$\kern#1\cr#2}}\right]}}

\def\rightarrowmat#1#2#3{
  \setbox1=\hbox{\small\kern#2$\bmatrix{#1}{#3}$\kern#2}
  \,\vbox{\offinterlineskip\hbox to\wd1{\hfil\copy1\hfil}
    \kern 3pt\hbox to\wd1{\rightarrowfill}}\,}

\def\leftarrowmat#1#2#3{
  \setbox1=\hbox{\small\kern#2$\bmatrix{#1}{#3}$\kern#2}
  \,\vbox{\offinterlineskip\hbox to\wd1{\hfil\copy1\hfil}
    \kern 3pt\hbox to\wd1{\leftarrowfill}}\,}

\def\rightarrowbox#1#2{
  \setbox1=\hbox{\kern#1\hbox{\small #2}\kern#1}
  \,\vbox{\offinterlineskip\hbox to\wd1{\hfil\copy1\hfil}
    \kern 3pt\hbox to\wd1{\rightarrowfill}}\,}

\def\leftarrowbox#1#2{
  \setbox1=\hbox{\kern#1\hbox{\small #2}\kern#1}
  \,\vbox{\offinterlineskip\hbox to\wd1{\hfil\copy1\hfil}
    \kern 3pt\hbox to\wd1{\leftarrowfill}}\,}








\def\quiremacro#1#2#3#4#5#6#7#8#9{
  \expandafter\def\csname#1\endcsname##1{
  \ifnum\srcdepth=1
  \magnification=#2
  \input quire
  \hsize=#3
  \vsize=#4
  \htotal=#5
  \vtotal=#6
  \shstaplewidth=#7
  \shstaplelength=#8
  \hoffset=\htotal
  \advance\hoffset by -\hsize
  \divide\hoffset by 2
  \ifnum\vsize<\vtotal
    \voffset=\vtotal
    \advance\voffset by -\vsize
    \divide\voffset by 2
  \fi
  \advance\voffset by #9
  \shhtotal=2\htotal
  \baselineskip=13pt
  \emergencystretch = 0.05\hsize
  \horigin=0.0truein
  \vorigin=0.0truein
  \shthickness=0pt
  \shoutline=0pt
  \shcrop=0pt
  \shvoffset=-1.0truein
  \ifnum##1>0\quire{#1}\else\qtwopages\fi
  \fi
}}



\quiremacro{letterbooklet} 
{1000}{4.79452truein}{7truein}{5.5truein}{8.5truein}{0.01pt}{0.66truein}{-.0625t
ruein}

\quiremacro{Afourbooklet}
{1095}{5.25truein}{7truein}{421truept}{595truept}{0.01pt}{0.66truein}{-.0625true
in}

\quiremacro{legalbooklet}
{1095}{5.25truein}{7truein}{7.0truein}{8.5truein}{0.01pt}{0.66truein}{-.0625true
in}

\quiremacro{twoupsub} 
{895}{4.5truein}{7truein}{5.5truein}{8.5truein}{0pt}{0pt}{.0625truein}


\quiremacro{Afourviewsub} 
{1000}{5.0228311in}{7.7625571in}{421truept}{595truept}{0.1pt}{0.5\vtotal}{-.0625
truein}


\quiremacro{viewsub}
{1095}{5.5truein}{8.5truein}{461truept}{666truept}{0.1pt}{0.5\vtotal}{-.125truei
n}


\newcount\countA
\newcount\countB
\newcount\countC

\def\monthname{\begingroup
  \ifcase\number\month
    \or January\or February\or March\or April\or May\or June\or
    July\or August\or September\or October\or November\or December\fi
\endgroup}

\def\dayname{\begingroup
  \countA=\number\day
  \countB=\number\year
  \advance\countA by 0 
  \advance\countA by \ifcase\month\or
    0\or 31\or 59\or 90\or 120\or 151\or
    181\or 212\or 243\or 273\or 304\or 334\fi
  \advance\countB by -1995
  \multiply\countB by 365
  \advance\countA by \countB
  \countB=\countA
  \divide\countB by 7
  \multiply\countB by 7
  \advance\countA by -\countB
  \advance\countA by 1
  \ifcase\countA\or Sunday\or Monday\or Tuesday\or Wednesday\or
    Thursday\or Friday\or Saturday\fi
\endgroup}

\def\timename{\begingroup
   \countA = \time
   \divide\countA by 60
   \countB = \countA
   \countC = \time
   \multiply\countA by 60
   \advance\countC by -\countA
   \ifnum\countC<10\toks1={0}\else\toks1={}\fi
   \ifnum\countB<12 \toks0={\sevenrm AM}
     \else\toks0={\sevenrm PM}\advance\countB by -12\fi
   \relax\ifnum\countB=0\countB=12\fi
   \hbox{\the\countB:\the\toks1 \the\countC \thinspace \the\toks0}
\endgroup}

\def\timestamp{\dayname, \the\day\ \monthname\ \the\year, \timename}


\print


\input diagrams.tex
\def\enma#1{{\ifmmode#1\else$#1$\fi}}

\def\mathbb#1{{\bbold #1}}
\def\mathbf#1{{\bf #1}}


\def\ZZ{\enma{\mathbb{Z}}}



\def\CC{\enma{\mathbf{C}}}

\def\th{{^{\rm th}}}

\def\set#1{\enma{\{#1\}}}
\def\setdef#1#2{\enma{\{\;#1\;\,|\allowbreak
  \;\,#2\;\}}}


\def\im{\mathop{\rm im}\nolimits}

\def\rank{\mathop{\rm rank}\nolimits}

\def\Ext{\mathop{\rm Ext}\nolimits}

\newsymbol\boxtimes1202

\font\tenmsam=msam10
\font\sevenmsam=msam7
\font\fivemsam=msam5
\newfam\msamfam         
\textfont\msamfam\tenmsam
\scriptfont\msamfam\sevenmsam
\scriptscriptfont\msamfam\fivemsam

%

\font\tenmsbm=msbm10
\font\sevenmsbm=msbm7
\font\fivemsbm=msbm5
\newfam\msbmfam         
\textfont\msbmfam\tenmsbm
\scriptfont\msbmfam\sevenmsbm
\scriptscriptfont\msbmfam\fivemsbm

%

\font\teneufm=eufm10
\font\seveneufm=eufm7
\font\fiveeufm=eufm5
\newfam\eufmfam        
\textfont\eufmfam\teneufm
\scriptfont\eufmfam\seveneufm
\scriptscriptfont\eufmfam\fiveeufm

%

%
%
%
%
%
%
%
\newcount\amsfamcount 
\newcount\classcount   
\newcount\positioncount
\newcount\codecount
\newcount\n             
\def\newsymbol#1#2#3#4#5{               
\n="#2                                  
\ifnum\n=1 \amsfamcount=\msamfam\else   
\ifnum\n=2 \amsfamcount=\msbmfam\else   
\ifnum\n=3 \amsfamcount=\eufmfam
\fi\fi\fi
\multiply\amsfamcount by "100           
\classcount="#3                 
\multiply\classcount by "1000           
\positioncount="#4#5            
\codecount=\classcount                  
\advance\codecount by \amsfamcount      
\advance\codecount by \positioncount
\mathchardef#1=\codecount}              
\newcount\famcnt 
\newcount\classcnt   
\newcount\positioncnt
\newcount\codecnt
\def\newmathsymbol#1#2#3#4#5{          
\famcnt=#2                      
\multiply\famcnt by "100        
\classcnt="#3                   
\multiply\classcnt by "1000     
\positioncnt="#4#5              
\codecnt=\classcnt              
\advance\codecnt by \famcnt     
\advance\codecnt by \positioncnt
\mathchardef#1=\codecnt}        

%
\newsymbol\varnothing203F
\def\emptyset{\mathop{\varnothing}}
\newsymbol\SEMI226F
\def\rtimes{\mathop{\SEMI}}
\newsymbol\smallsetminus2272
\def\setminus{\mathop{\smallsetminus}}

\def\Box{
  \ifmmode\eqno\qed
  \else\ifvmode\removelastskip\line{\hfil\qed}
  \else\unskip\quad\hskip-\hsize
    \hbox{}\hskip\hsize minus 1em\qed\par
  \fi\stdbreak\fi}

\def\Sec {\mathop{\rm Sec}\nolimits}
\def\Sym {\mathop{\rm Sym}\nolimits}
\def\Sing{\mathop{\rm Sing}\nolimits}
\def\Join{\mathop{\rm Join}\nolimits}
\def\PGL {\mathop{\rm PGL}\nolimits} 
\def\GL {\mathop{\rm GL}\nolimits} 
\def\SL {\mathop{\rm SL}\nolimits}
\def\PSL {\mathop{\rm PSL}\nolimits}
\def\Pic {\mathop{\rm Pic}\nolimits}
\def\Gr  {\mathop{\rm Gr}\nolimits}
\def\Ext {\mathop{\rm Ext}\nolimits}
\def\Proj {\mathop{\rm Proj}\nolimits}
\def\VSP {\mathop{\rm VSP}\nolimits}
\def\mod {\mathop{\rm mod}}
\def\th {{^{\rm th}}}

\def\diag {\mathop{\rm diag}\nolimits}
\def\mult {\mathop{\rm mult}\nolimits}
\def\vt{\enma{\vartheta}}
\def\boldz{{\bf Z}}
\def\HHH{{\bf H}}
\def\P{{\bf P}}
\def\Pone{{\bf P}^1}
\def\Ptwo{{\bf P}^2}
\def\Pthree{{\bf P}^3}
\def\Pfour{{\bf P}^4}
\def\Pfive{{\bf P}^5}
\def\Psix{{\bf P}^6}
\def\Pseven{{\bf P}^7}
\def\Ptzz{\Pthree_-/\boldz_2\times\boldz_2}
\def\Pzz{\Ptwo_-/\boldz_2\times\boldz_2}
\def\A{{\cal A}}
\def\C{{\cal C}}
\def\M{{\cal M}}

\def\F{{\cal F}}
\def\G{{\cal G}}
\def\H{{\cal H}}
\def\I{{\cal I}}
\def\K{{\cal K}}
\def\L{{\cal L}}
\def\O{{\cal O}}
\def\SS{{\cal S}}
\def\X{{\cal X}}

\def\dual#1{{#1}^{\scriptscriptstyle \vee}}
\def\exact#1#2#3{0\longrightarrow#1\longrightarrow#2
\longrightarrow#3\longrightarrow0}
\def\mapright#1{{\smash{\mathop{\longrightarrow}\limits^{#1}}}}
\def\tenpoint{%
\textfont0=\tenrm \scriptfont0=\sevenrm
\scriptscriptfont0=\fiverm \def\rm{\fam0\tenrm}%
\textfont1=\teni \scriptfont1=\seveni
\scriptscriptfont1=\fivei \def\oldstyle{\fam1\teni}%
\textfont2=\tensy \scriptfont2=\sevensy
\scriptscriptfont2=\fivesy
\textfont\itfam=\tenit \def\it{\fam\itfam\tenit}%
\textfont\slfam=\tensl \def\sl{\fam\slfam\tensl}%
\textfont\ttfam=\tentt \def\tt{\fam\tfam\tentt}%
\textfont\bffam=\tenbf \scriptfont\bffam=\sevenbf
\scriptscriptfont\bffam=\fivebf \def\bf{\fam\bffam\tenbf}%
\abovedisplayskip=12pt plus 3pt minus 9pt
\belowdisplayskip=\abovedisplayskip
\abovedisplayshortskip=0pt plus 3pt
\belowdisplayshortskip=7pt plus 3pt minus 4pt
\smallskipamount=3pt plus 1pt minus 1pt
\medskipamount=6pt plus 2pt minus 2pt
\bigskipamount=12pt plus 4pt minus 4pt
\setbox\strutbox=\hbox{\vrule height8.5pt depth3.5pt width0pt}%
\normalbaselineskip=12pt \normalbaselines \rm}

\def\eightpoint{%
\font\eightrm=cmr8%
\font\sixrm=cmr6%
\font\eighti=cmmi6%
\font\eightit=cmti8%
\font\sixi=cmmi6%
\font\sixit=cmti6%
\font\eightsy=cmsy8%
\font\sixsy=cmsy6%
\font\eightsl=cmsl8%
\font\eighttt=cmtt6%
\font\eightbf=cmbx8%
\font\sixbf=cmbx6%
\textfont0=\eightrm \scriptfont0=\sixrm
\scriptscriptfont0=\fiverm \def\rm{\fam0\eightrm}%
\textfont1=\eighti \scriptfont1=\sixi
\scriptscriptfont1=\fivei \def\oldstyle{\fam1\eighti}%
\textfont2=\eightsy \scriptfont2=\sixsy
\scriptscriptfont2=\fivesy
\textfont\itfam=\eightit \def\it{\fam\itfam\eightit}%
\textfont\slfam=\eightsl \def\sl{\fam\slfam\eightsl}%
\textfont\ttfam=\eighttt \def\tt{\fam\tfam\eighttt}%
\textfont\bffam=\eightbf \scriptfont\bffam=\sixbf
\scriptscriptfont\bffam=\fivebf \def\bf{\fam\bffam\eightbf}%
\abovedisplayskip=9pt plus 3pt minus 9pt
\belowdisplayskip=\abovedisplayskip
\abovedisplayshortskip=0pt plus 2pt
\belowdisplayshortskip=5pt plus 2pt minus 3pt
\smallskipamount=2pt plus 1pt minus 1pt
\medskipamount=4pt plus 2pt minus 2pt
\bigskipamount=9pt plus 4pt minus 4pt
\setbox\strutbox=\hbox{\vrule height 7pt depth 2pt width 0pt}%
\normalbaselineskip=9pt \normalbaselines \rm}


\forward{prelim}{Section}{1}
\forward{1.4}{Section}{2}
\forward{HM}{Section}{3}
\forward{1.6}{Section}{4}
\forward{1.7}{Section}{5}
\forward{CY6}{Theorem}{4.10}
\forward{CY10}{Theorem}{7.4}
\forward{CY14}{Theorem}{2.2}
\forward{CY7}{Proposition}{5.2}
\forward{minors8}{Definition}{6.11}
\forward{singV8y}{Theorem}{6.5}
\forward{defect6}{Remark}{4.11}
\forward{defect8}{Theorem}{6.9}
\forward{kaehler8}{Remark}{6.10}
\forward{rat.map6}{Definition}{4.6}
\forward{HM.wedge}{Theorem}{3.2}
\forward{pencil2}{Proposition}{6.14}
\forward{k10}{Remark}{7.5}

\font\smallrm=cmcsc8

\headline={\ifodd\pageno \ifnum\pageno>1 \smallrm \hfil 
Calabi-Yau Threefolds and Moduli of Abelian Surfaces I
\hfil\folio \else\hfill\fi \else \smallrm \folio \hfill
Mark Gross and Sorin Popescu  
\hfill\fi} \footline={\hss}   

\centerline{\hd Calabi-Yau Threefolds and Moduli of Abelian Surfaces I}
\medskip
\centerline{\it Mark Gross\footnote{*}{Supported by NSF grant DMS-9700761}
and Sorin Popescu\footnote{**}{Supported by
NSF grant DMS-9610205 and MSRI, Berkeley}}
\bigskip

\medskip
{\settabs 3 \columns
\+Mathematics Institute&&
Department of Mathematics\cr
\+University of Warwick&&
Columbia University\cr
\+Coventry, CV4 7AL, UK&&
New York, NY 10027\cr
\+mgross@maths.warwick.ac.uk&&
psorin@math.columbia.edu\cr}

\bigskip
\bigskip


Let $\A_{d}$ denote the moduli space of polarized abelian surfaces of
type $(1,d)$, and let $\A_{d}^{lev}$ be the moduli space of polarized
abelian surfaces with canonical level structure. Both are (possibly
singular) quasi-projective threefolds, and $\A_{d}^{lev}$ is a finite
cover of $\A_{d}$. We will also denote by $\widetilde\A_{d}$ and
$\widetilde\A_{d}^{lev}$ nonsingular models of suitable
compactifications of these moduli spaces.  We will use in the sequel
definitions and notation as in [GP1], [GP2]; see also [Mu1], [LB] and
[HKW] for basic facts concerning abelian varieties and their moduli.
Throughout the paper the base field will be $\CC$.

The main goal of this paper, which is a continuation of [GP1] and
[GP2], is to describe birational models for moduli spaces of these
types for small values of $d$. Since the Kodaira dimension is a
birational invariant, thus independent of the chosen compactification,
we can decide the uniruledness, unirationality or rationality of
nonsingular models of (any) compactifications of these moduli spaces.

Motivation for our project has come from several directions:

Tai, Freitag and Mumford have proved that moduli spaces of principally
polarized abelian varieties are of general type when the dimension $g$
of the abelian varieties is large enough (in fact $\ge 7$).  However
for small dimensions, they are rational or unirational and have nice
projective models: see the work of Katsylo [Kat] for $g=3$, 
van Geemen [vG] and Dolgachev-Ortland [DO] for $g=3$ with level $2$
structure, Clemens [Cle] for $g=4$, and Donagi [Do], 
Mori-Mukai [MM] and Verra [Ver] for $g=5$. 
However, it is an open problem to determine the Kodaira
dimension of the moduli space for $g=6$.

Using a version of the Maass-Kurokawa lifting, Gritsenko [Gri1],
[Gri2] has recently proved the existence of weight 3 cusp forms with
respect to the paramodular group $\Gamma_t$, for almost all values of
$t$. Since one knows the dimension for the space of Jacobi cusp forms
one deduces in this way lower bounds for the dimensions of the spaces
of cusp forms with respect to the paramodular group $\Gamma_t$, and
thus for the plurigenera of the corresponding moduli spaces.  More
precisely he has shown that:
$$\vbox{${\A}_{d}$ is not unirational 
(in fact $p_g(\widetilde{\A}_{d})\ge 1$) if $d\ge 13$ and 
$d\ne 14,\; 15,\; 16,\; 18,\; 20,\; 24,\; 30,$ $36$.}$$
 
In fact it is pointed out in a note by Gritsenko and Hulek [Gri1] that
the same method shows that ${\A}^{lev}_{p}$ is of general type for all
primes $p\ge 37$. See also [Bori], [HS1] and [HS2] for related results.

A few other moduli spaces have known descriptions:

\item{$\bullet$} $\A_1\cong\overline{\M_2}\setminus\Delta_0$ via the
Jacobian map.  $\overline{\M_2}$ is the moduli space of stable curves
of genus 2 and $\Delta_0$ stands for the divisor of the curves with at
least one non-disconnecting node. In particular $\A_1$ is 
rational [I].

\item{$\bullet$} $\A_2$ and $\A_3$ are rational [BL]. It
follows from results of [Ba] that $\A_{(2,4)}^{lev}$, the moduli space
of $(2,4)$-polarized abelian surfaces with canonical level structure,
is birational to $\Pone\times\Pone\times\Pone$.  On the other hand it
doesn't appear to be known whether $\A^{lev}_3$ is unirational.

\item{$\bullet$} $\A_4^{lev}$ is rational. See [BLvS] for a proof of
this and for the geometry of $(1,4)$-polarized abelian surfaces.

\item{$\bullet$} $\overline{\A^{lev}_5}\cong\P(H^0(\F_{HM}(3)))$, 
where bar stands for the toroidal compactification for the Voronoi or
Igusa decomposition, and $\F_{HM}$ is the Horrocks-Mumford bundle on
$\Pfour$. See [HM] and [HKW] for details, and \ref{HM} for a
brief review of the relevant facts.

\item{$\bullet$} $\A^{lev}_7$ is birational to $V_{22}$,
a prime Fano threefold of index 1 and  genus 12 [MS].
In \ref{1.7} we will give  a short proof of this result.
It is likely that in this case $V_{22}$ coincides with
the toroidal compactification.

\item{$\bullet$} $\A_9$ is rational (see [O'G] for unirationality and
[GP2]).  In fact $\A_9^{lev}$ is naturally birational to $\Pthree$
([GP2]).  O'Grady shows also that $\tilde\A_{p^2}$ is a threefold of
general type for all prime numbers $p\ge 11$.

\item{$\bullet$} $\A_{11}^{lev}$ is  birational to the 
Klein cubic threefold 
$$\K=\set{\sum_{i=0}^4x_i^2x_{i+1}=0}\subset\Pfour$$ 
See [GP2] for details.  This is the
unique $\PSL_2(\boldz_{11})$-invariant cubic in $\Pfour$, and
$\PSL_2(\boldz_{11})$ is its full automorphism group [Ad].  The Klein
cubic being smooth is unirational but not rational.  It would be
interesting also in this case to compare $\K$ with the toroidal
compactification.

\medskip
In this paper and its sequel [GP3], we will focus on (most of) the
moduli spaces which Gritsenko has not shown to have non-negative
Kodaira dimension. In particular, we will give details in this paper
as to the structure of $\A_d^{lev}$, $d=6$, $8$ and $10$, and in [GP3]
will discuss $\A^{lev}_{12}$ and $\A_d$ for $d=14,16,18$ and $20$.
In this paper we prove the following results:

\theorem{results}
\item{a)} $\A_6^{lev}$ is birational to a non-singular
quadric hypersurface in $\Pfour$.
\item{b)} $\A_8^{lev}$ is birational to a rational conic bundle over $\Ptwo$.
\item{c)} $\A_{10}^{lev}$ is birational to a quotient 
$\Pthree/\boldz_2\times\boldz_2$,
where the action on $\Pthree$ is given by 
$$(x_1:x_2:x_3:x_4)\mapsto (x_4:x_3:x_2:x_1)\qquad {\rm and}
\qquad (x_1:x_2:x_3:x_4)\mapsto (x_1:-x_2:x_3:-x_4)$$
This quotient is rational and  isomorphic to a (singular)
prime Fano threefold of genus $9$, index $1$ in $\P^{10}$.

Of these, the description of $\A_{10}^{lev}$ is the easiest, following
immediately from [GP1]. The description of $\A_6^{lev}$ follows from
the fact that a $(1,6)$-polarized abelian surface $A\subseteq\Pfive$
is determined (but not completely cut out by) the cubics vanishing on
$A$.  Since the structure of the ideal of a $(1,6)$-polarized abelian
surface was not previously known, we will provide details in
\ref{1.6}.

The structure of $\A_8^{lev}$ is the most difficult to analyze because
it lies at the boundary between those surfaces determined by cubics
and those determined by quadrics. In fact we will show that if
$A\subseteq\Pseven$ is a general $(1,8)$-polarized abelian surface,
then its homogeneous ideal is generated by $4$ quadrics and $16$
cubics. The $4$ quadrics cut out a complete intersection threefold
$X\subset\Pseven$ containing $A$, which in general has as singular
locus 64 ordinary double points. We show that such an $X$ in fact
contains a pencil of $(1,8)$-polarized abelian surfaces, and since the
family of such Heisenberg invariant threefolds turns out to be an open
subset of $\Ptwo$, we obtain a description of $\A_{lev}^8$ as a
$\Pone$-bundle over an open subset of $\Ptwo$.

In the above outline, the threefold $X$ plays a special role. It is a
Calabi-Yau threefold since there is a small resolution $\pi:\widetilde
X\longrightarrow X$ with $\omega_{\widetilde X}\cong \O_{\widetilde
X}$. Furthermore, as $X$ contains a pencil of abelian surfaces, this
small resolution can be chosen so that $\widetilde X$ possesses an abelian
surface fibration. Such fibrations are of independent interest in the
study of Calabi-Yau threefolds. Historically there have been very few
examples; in fact, for a long time the only examples of a Calabi-Yau
threefolds with a fibration whose general fibre is a simple abelian
surface have been the Horrocks-Mumford quintics, whose geometry we
review in \ref{HM}. The Horrocks-Mumford quintic in $\Pfour$ is
defined as the zero locus of the wedge of two independent sections in
$H^0(\F_{HM}(3))$, where $\F_{HM}$ is the Horrocks-Mumford bundle.
Since a (general) section of this bundle vanishes along a
$(1,5)$-polarized abelian surface in $\Pfour$, a Horrocks-Mumford
quintic contains a pencil of such surfaces.

It turns out that in every case considered in this paper, there exist
similar examples of Calabi-Yau threefolds. This is the second focus of
this paper. We find similarities in many of these examples. In
particular, many of them contain second pencils of abelian surfaces of
a different type.  The classic example of this occurs in the case of
Horrocks-Mumford quintics, where a surface obtained via liaison
starting from a $(1,5)$-polarized abelian surface is a non-minimal
$(2,10)$-polarized abelian surface. We will see similar phenomena also
occur in many of the other examples. Another remark is that
all the Calabi-Yau threefolds we discuss in this paper have Euler
characteristic zero.

There is a natural reason for the existence of such Calabi-Yau
threefolds.  In fact, this is motivated by the existence of two series
of degenerate ``Calabi-Yau'' threefolds. The first series, which
occurs in $\P^{n-1}$ for $n\ge 5$, are the secant varieties of
elliptic normal curves in $\P^{n-1}$. These were studied in [GP1], \S
5. Many of the examples of Calabi-Yau threefolds given here are
partial smoothings of these secant varieties. This also works
naturally in case $n=4$, if one thinks of the secant variety as a
double cover of $\Pthree$. As a result, we discuss in \ref{1.4} of
this paper the $(1,4)$ case. We also include a section on
$(1,7)$-polarized abelian surfaces, where the Calabi-Yau threefold
which arises in this fashion (first noted by Alf Aure and Kristian
Ranestad, unpublished) plays an important role in the description of
the moduli.  In particular, we give in \ref{1.7} an alternative
approach to rationality of $\A_7^{lev}$ to that given in [MS].

Another series of degenerate ``Calabi-Yau'' threefolds appear in
$\P^{2n-1}$, $n\ge 3$, as the join of two elliptic normal curves of
degree $n$, lying in disjoint linear subspaces of $\P^{2n-1}$.  These
are threefolds of degree $n^2$ in $\P^{2n-1}$.  We discuss the
geometry of these singular threefolds briefly in \ref{prelim}. This
series gives rise to other classes of Calabi-Yau threefolds, as
partial smoothings, in the $(1,8)$ and $(1,10)$ cases.  Much of this
paper is devoted to describing the geometry of these threefolds.

In brief, the Calabi-Yau threefolds discussed in this paper can be
described as follows. Partial smoothings of secant varieties are
\item{$\bullet$} Double covers of $\Pthree$ branched over 
certain nodal octics, see \ref{CY14}. These double covers contain
pencils of $(1,4)$-polarized abelian surfaces.
\item{$\bullet$} Horrocks-Mumford quintics, containing both a pencil
of $(1,5)$ and a pencil of $(2,10)$-polarized abelian surfaces, see
\ref{HM.wedge}.
\item{$\bullet$} Nodal complete intersections of two cubics in $\Pfive$,
containing both a pencil of $(1,6)$ and a pencil of $(2,6)$-polarized
abelian surfaces, see \ref{CY6}.
\item{$\bullet$} Calabi-Yau threefolds defined by the $6\times 6$
Pfaffians of certain $7\times 7$-skew-symmetric matrix of linear forms
in $\Psix$, see \ref{CY7}. These contain both a pencil of $(1,7)$ and
a pencil of $(1,14)$-polarized abelian surfaces. Some details
concerning the geometry of these threefolds are delayed until [GP3] in
the discussion of $(1,14)$-polarized abelian surfaces.
\item{$\bullet$} Calabi-Yau threefolds defined by the $3\times 3$
minors of a $4\times 4$ matrix of linear forms in $\Pseven$, which
contain a pencil of $(1,8)$-polarized abelian surfaces, see
\ref{minors8}.
\item{$\bullet$} Calabi-Yau threefolds defined by the $3\times 3$
minors of a $5\times 5$ symmetric matrix of linear forms in $\P^9$,
which contain a pencil of $(1,10)$-polarized abelian surfaces, see
\ref{CY10}.

The reader may note that there is a gap corresponding to the partial
smoothing of the secant variety of a degree 9 elliptic normal curve,
which should contain a pencil of $(1,9)$-polarized abelian
surfaces. This does not mean that such a Calabi-Yau threefold does not
exist, but we have not been able to find it!  There is also a nice
analogy between the sequence of Calabi-Yau threefolds described above
and Del Pezzo surfaces. Del Pezzo surfaces of degree 2 through 6 can
be described as a double cover of $\Ptwo$, a cubic hypersurface in
$\Pthree$, a complete intersection of two quadrics in $\Pfour$, a
surface defined by the $4\times 4$-Pfaffians of a $5\times 5$
skew-symmetric matrix of linear forms in $\Pfive$, and a surface
defined by the $2\times 2$ minors of a $3\times 3$ matrix of linear
forms in $\Psix$, respectively. Furthermore, the Del Pezzo surface of
degree $8$ in $\P^8$ isomorphic to $\Pone\times\Pone$ can be described
by the $2\times 2$ minors of a $4\times 4$ symmetric matrix of linear
forms. While we do not explore this analogy further, it is a curious
one! From this perspective the missing Calabi-Yau, which should
contain a pencil of $(1,9)$-polarized abelian surfaces, would be
a ladder determinantal variety of degree $27$ defined by those
$3\times 3$-minors of a $5\times 5$ symmetric matrix of linear forms
in $\P^8$ which do not involve the lower right corner of the matrix.

Partial smoothings of the join of two elliptic curves naturally
occurring in this paper are described as
\item{$\bullet$} Complete intersections of two cubics in $\P^5$, containing
a pencil of $(1,6)$-polarized abelian surfaces.
\item{$\bullet$} Complete intersections of 4 quadrics in $\Pseven$,
containing a pencil of $(1,8)$ and a pencil of $(2,8)$-polarized abelian
surfaces, see \ref{singV8y}, \ref{defect8} and \ref{kaehler8}.
\item{$\bullet$} The proper intersection of ``two copies'' of the
Pl\"ucker embedding of the Grassmannian $\Gr(2,5)\subseteq\P^9$,
containing a pencil of $(1,10)$ and a pencil of $(3,15)$-polarized
abelian surfaces, see \ref{CY10} and \ref{k10}.

Further examples of such Calabi-Yau threefolds will appear in [GP3].

\smallskip 
{\it Acknowledgments:} It is a pleasure to thank Kristian
Ranestad, who joined one of us in discussions leading to some of the
ideas in this paper, and Ciro Ciliberto, Wolfram Decker, Igor
Dolgachev, David Eisenbud, Klaus Hulek, Nicolae Manolache, Christian
Peskine, Greg Sankaran, Frank Schreyer and Alessandro Verra from whose
ideas the exposition has benefited. We are also grateful to Mike
Stillman, Dave Bayer, and Dan Grayson for the programs {\it
Macaulay\/} [BS], and {\it Macaulay2\/} [GS] which have been extremely
useful to us; without them we would perhaps have never been bold
enough to guess the existence of the structures that we describe here.

\section {prelim} {Preliminaries.}

We review our notation and conventions concerning abelian surfaces; more
details can be found in [GP1].

Let $(A,\L)$ be a general abelian surface with a polarization of type
$(1,d)$.  If $d\ge 5$, then $|\L|$ induces an embedding of
$A\subset\P^{d-1}=\P(\dual{H^0(\L)})$ of degree $2d$.  The line bundle
$\L$ induces a natural map from $A$ to its dual $\hat A$,
$\phi_{\L}:A\longrightarrow \hat A$, given by $x\mapsto
t_x^*\L\otimes\L^{-1}$, where $t_x:A\longrightarrow A$ is the morphism
given by translation by $x\in A$.  Its kernel $K(\L)$ is isomorphic to
$\boldz_{d}\times\boldz_{d}$, and is dependent only on the
polarization $c_1(\L)$.

For every $x\in K(\L)$ there is an isomorphism $t_x^*\L
\cong \L$. This induces a projective representation
$K(\L)\longrightarrow {\PGL}(H^0(\L))$, which lifts uniquely to a
linear representation of $K(\L)$ after taking a central extension
of $K(\L)$ 
$$1\longrightarrow{{\CC}^*}\longrightarrow{\G(\L)}
\longrightarrow{K(\L)}\longrightarrow0,
$$ 
whose Schur commutator map is the Weil pairing.  $\G(\L)$ is the {\it
theta group} of $\L$ and is isomorphic to the abstract Heisenberg
group $\H(d)$, while the above linear representation is isomorphic to the
Schr\"odinger representation of $\H(d)$ on $V={\CC}(\boldz_{d})$,
the vector space of complex-valued functions on $\boldz_{d}$.  An
isomorphism between $\G(\L)$ and $\H(d)$, which restricts to the
identity on centers induces a symplectic isomorphism between $K(\L)$
and $\boldz_{d}\times\boldz_{d}$. Such an isomorphism is called a {\it
level structure of canonical type} on $(A,c_1(\L))$. (See [LB],
Chapter 8, \S 3 or [GP1], \S 1.)

A decomposition $K(\L)=K_1(\L)\oplus K_2(\L)$, with $K_1(\L)\cong
K_2(\L) \cong \boldz_{d}$ isotropic subgroups with respect to the Weil
pairing, and a choice of a characteristic $c$ ([LB], Chapter 3, \S 1)
for $\L$, define a unique natural basis $\setdef{\vt^c_x}{ x\in K_1(\L)}$ of
{\it canonical theta functions} for the space $H^0(\L)$ (see [Mu2]
and [LB], Chapter 3, \S 2). This basis allows an identification of
$H^0(\L)$ with $V$ via $\vt_{\gamma}^c\mapsto x_{\gamma}$, where
$x_{\gamma}$ is the function on $\boldz_{d}$ defined by
$x_{\gamma}(\delta)=\cases{1&$\gamma=\delta$\cr
0&$\gamma\not=\delta$\cr}$ for $\gamma,\delta\in\boldz_{d}$. The functions
$x_0,\ldots,x_{d-1}$ can also be identified with coordinates on
$\P(\dual{H^0(\L)})$. Under this identification, the representation
$\G(\L)\longrightarrow \GL(H^0(\L))$ coincides with the Schr\"odinger
representation $\H(d)\longrightarrow \GL(V)$. We will only consider the
action of $\HHH_{d}$, the finite subgroup of $\H(d)\longrightarrow \GL(V)$
generated in the Schr\"odinger representation by $\sigma$ and $\tau$,
where
$$ \sigma(x_i)=x_{i-1}, \qquad \tau(x_i)={\xi}^{-i} x_i,$$ for all
$i\in{\boldz}_{d}$, and where $\xi=e^{2{\pi}i\over d}$ is a primitive
root of unity of order $d$. Notice that $\lbrack\sigma,\tau\rbrack=\xi$, 
thus ${\HHH}_{d}$ is a central extension
$$1 \longrightarrow {\bf\mu_{d}} \longrightarrow {\HHH_{d}} 
\longrightarrow
{\boldz_{d}}\times{\boldz_{d}} \longrightarrow 0.$$
Therefore the choice of a canonical level structure means that
if $A$ is embedded in $\P(\dual{H^0(\L)})$ using as coordinates
$x_{\gamma}=\vt_{\gamma}^c$,  $\gamma\in
\boldz_{d}$, then the image of $A$ will be invariant under the 
action of the Heisenberg group $\HHH_{d}$ via the Schr\"odinger 
representation. (See [LB], Chapter 6, \S 7 for details.)

If moreover the line bundle $\L$ is chosen to be symmetric (and there
are always finitely many choices of such an $\L$ for a given
polarization type), then the embedding via $|\L|$ is also invariant
under the involution $\iota$, where
$$ \iota(x_i)=x_{-i}, \qquad i\in{\boldz}_{d}.$$
This involution restricts to $A$ as the involution $x\mapsto -x$.
We will denote  by $\P_+$ and $\P_-$ the $(+1)$ and
$(-1)$-eigenspaces of the involution $\iota$, respectively.
We will also denote as usual by 
$\HHH_n^e:=\HHH_n\rtimes\langle\iota\rangle$.

We also recall a key result from [GP1]: In that paper, on $\P^{2d-1}\times
\P^{2d-1}$, we have introduced a matrix
$$M_d(x,y)=(x_{i+j}y_{i-j}+x_{i+j+d}y_{i-j+d})_{0\le i,j\le d-1},$$
where the indices of the variables $x$ and $y$ above are all modulo
$2d$. This matrix has the property that if $A\subseteq \P^{2d-1}$
is a Heisenberg invariant $(1,2d)$-polarized abelian surface,
then $M_d$ has rank at most two on $A\times A\subseteq \P^{2d-1}\times
\P^{2d-1}$. Similarly, if $A\subseteq \P^{2d}$
is a Heisenberg invariant $(1,2d+1)$-polarized abelian surface,
then the (Moore) matrix
$$M_{2d+1}'(x,y)=(x_{d(i+j)}y_{d(i-j)})_{i\in\boldz_{2d+1},j\in\boldz_{2d+1}},
$$
on $\P^{2d}\times \P^{2d}$  has rank at most four 
on $A\times A\subseteq \P^{2d}\times\P^{2d}$.
These matrices will prove to be ubiquitous!
\medskip

We include at the end of this section a brief discussion of abelian
surface fibrations on Calabi-Yau threefolds. Throughout this paper
we will use the following terminology:

\definition{CY} A {\it Calabi-Yau threefold} is a non-singular projective
threefold satisfying $\omega_X\cong\O_X$ and $h^1(\O_X)=0$.

\lemma{ABpencil} Let $X$ be a Calabi-Yau threefold, and let $A\subseteq X$
be a minimal abelian surface. Then $A$ is a member of a base-point free
linear system of abelian surfaces which induces a fibration
$\pi:X\longrightarrow \Pone$ with $A$ as a fibre.

\proof:
On $X$, we have
$$\exact{\O_X}{\O_X(A)}{\omega_A},$$
from which we obtain $\dim |A|=1$. It then follows from
[Og] that $|A|$ must be a base-point free linear system
inducing an abelian surface fibration $X\longrightarrow\Pone$. \Box

The Calabi-Yau threefolds we will consider in this paper will be
partial smoothings of two types of singular threefolds. Recall first
[GP1], Proposition 5.1:

\proposition{recallsec} Let $E\subseteq \P^{n-1}$ be an elliptic normal
curve of degree $n$. Then
\item{\rm (1)} $\Sec(E)$ is an irreducible threefold of degree $n(n-3)/2$.
\item{\rm (2)} $\Sec(E)$ is non-singular outside $E$.
\item{\rm (3)} $\omega_{\Sec(E)}\cong\O_{\Sec(E)}$ and $h^1(\O_{\Sec(E)})=0$.

Thus, form a numerical point of view, $\Sec(E)$ is a degenerate
``Calabi-Yau'' threefold, containing a pencil of type II degenerations
of $(1,n)$-polarized abelian surfaces.  Namely, assume that
$E\subset\P^{n-1}$ is Heisenberg invariant under the Schr\"odinger
representation. Then $\Sec(E)$ contains the pencil of Heisenberg
invariant translation scrolls
$$S_{E,\rho}:=\bigcup_{P\in E} \langle P, P+\rho \rangle,$$ 
where $\langle P,P+\rho\rangle$ denotes the line spanned by $P$ and $P+\rho$
and $\rho\in E$ is not a 2-torsion point.  The
fibers of the pencil for $\rho$ a $2$-torsion point on $E$ are
multiplicity two structures on the corresponding smooth 
elliptic scrolls $S_\rho$.

A second series of degenerate ``Calabi-Yau'' threefolds is given as follows:

\proposition{join} Let $L_1,L_2\subseteq\P^{2n-1}$ be two disjoint
linear subspaces of dimension $n-1$, and let $E_1\subseteq L_1$,
$E_2\subseteq L_2$ be two elliptic normal curves of degree $n$. Put
$$\Join(E_1,E_2):=
\bigcup_{e_1\in E_1\atop e_2\in E_2} \langle e_1,e_2\rangle$$
where $\langle e_1,e_2\rangle$ denotes the linear span of $e_1$ and
$e_2$. Then
\item{\rm (1)} $\Join(E_1,E_2)$ is an irreducible threefold of degree
$n^2$ in $\P^{2n-1}$.
\item{\rm (2)} $\Join(E_1,E_2)$ is non-singular outside of $E_1$ and $E_2$.
\item{\rm (3)} $\omega_{\Join(E_1,E_2)}\cong\O_{\Join(E_1,E_2)}$
and $h^1(\O_{\Join(E_1,E_2)})=0$.

\proof: (1) Irreducibility is obvious. The degree of $\Join(E_1,E_2)$ is
the product of the degrees of $E_1$ and $E_2$, hence the join 
is of degree $n^2$.

(2) Let $Y\subseteq E_1\times E_2\times\P^{2n-1}$ be defined by
$$Y=\setdef{(e_1,e_2,y)}{y\in\langle e_1,e_2\rangle}.$$ The projection
$p_{12}:Y\longrightarrow E_1\times E_2$ clearly gives $Y$ the
structure of a $\Pone$-bundle, while the projection
$p_3:Y\longrightarrow\P^{2n-1}$ gives a resolution of singularities of
$\Join(E_1,E_2)$.  Note that $p_3:Y\setminus p_3^{-1}(E_1\cup
E_2)\longrightarrow \P^{2n-1}$ is an embedding. Indeed, $p_3$ is
one-to-one away from $p_3^{-1}(E_1\cup E_2)$; otherwise there exist
distinct points $e_1,e_1'\in E_1$ and $e_2,e_2'\in E_2$ such that
$\langle e_1,e_2\rangle\cap \langle e_1',e_2'\rangle\not=\emptyset$.
But then these two lines span only a $\Ptwo$, so that $\langle
e_1,e_1'\rangle\cap \langle e_2,e_2'\rangle
\not=\emptyset$, contradicting the assumption that 
$L_1\cap L_2=\emptyset$. To see that $p_3$ is
an immersion, a local calculation suffices. Consider an affine patch  
${\CC}^{2n-1}$ of $\P^{2n-1}$, with coordinates
$y_1,\ldots,y_{2n-1}$, in which $L_1$ is defined by
$y_n=\cdots=y_{2n-1}=0$ and $L_2$ is defined by 
$y_1=\cdots=y_{n-1}=0$, and $y_n=1$. If $E_1$ is locally parametrized by
$t\mapsto(\alpha_1(t),\ldots,\alpha_{n-1}(t),0,\ldots,0)$ and $E_2$ is
parametrized by $u\mapsto(0,\ldots,0,1,\beta_1(u),\ldots,\beta_{n-1}(u))$,
then $Y$ has local coordinates $(t,u,s)$ in which $p_3$ is given by
$$(t,u,s)\mapsto (s\alpha_1(t),\ldots,s\alpha_{n-1}(t),(1-s),(1-s)\beta_1(u),
\ldots,(1-s)\beta_{n-1}(u)).$$
Computing the differential of this map one sees easily that  
it is injective for $s\not=0,1$, 
given that $E_1$ and $E_2$ are nonsingular. 
Thus $p_{3}: Y\setminus p_{3}^{-1}(E_1\cup E_2)\longrightarrow \P^{2n-1}$
is an embedding, and so $\Join(E_1,E_2)\setminus(E_1\cup E_2)$ 
is non-singular.

(3) First note that $\Join(E_1,E_2)$ is normal.  Indeed, locally the
singularities of $\Join(E_1,E_2)$ look like
$\hbox{curve}\times\hbox{(cone over elliptic normal curve)}$, hence
are normal. Now the map $p_3:Y\longrightarrow \Join(E_1,E_2)$ contracts
two sections $\sigma_1$ and $\sigma_2$ of the $\Pone$-bundle
$p_{12}:Y\longrightarrow E_1\times E_2$ to $E_1$ and $E_2$,
respectively. Since $p_{12}$ is a $\Pone$-bundle, we can write
$K_Y=-\sigma_1-\sigma_2 +p_{12}^*D$ for some divisor $D$ on $E_1\times
E_2$. But $\sigma_1$ and $\sigma_2$ are both isomorphic to $E_1\times
E_2$ and are disjoint, so by adjunction
$$\eqalign{0=K_{\sigma_1}&=(K_Y+\sigma_1)|_{\sigma_1}\cr
&=(-\sigma_2+p_{12}^*D)|_{\sigma_1}\cr
&=D.\cr}$$
Thus $K_Y=-\sigma_1-\sigma_2$. On the other hand, on $\Join(E_1,E_2)$,
both $E_1$ and $E_2$ are curves of simple elliptic singularities, 
so $\Join(E_1,E_2)$ is Gorenstein and 
$$K_Y=p_{3}^*K_{\Join(E_1,E_2)}-\sigma_1-\sigma_2,$$
from which we conclude that $p_{3}^*K_{\Join(E_1,E_2)}=0$.
Since $p_{3}^*:\Pic(\Join(E_1,E_2))\longrightarrow\Pic(Y)$ is injective,
this shows $K_{\Join(E_1,E_2)}=0$, as desired.

Just as in the proof of [GP1], Proposition 5.1, to show that
$h^1(\O_{\Join(E_1,E_2)})=0$, it is enough to show that 
$\Pic^0(\Join(E_1,E_2))$ is discrete. Note that the image of 
$p_{3}^*:\Pic^0(\Join(E_1,E_2))\longrightarrow\Pic^0(Y)$ is
contained in the subgroup ${\cal P}$ of $\Pic^0(Y)$ given by
$${\cal P}=\setdef{\L\in\Pic^0(Y)}
{\hbox{$\L|_{\sigma_1}=p_1^*\M_1$ and $\L|_{\sigma_2}=p_2^*\M_2$ for
some $\M_i\in\Pic(E_i)$}}.$$
Thus it is enough to show ${\cal P}$ is discrete. Suppose $D\in
{\cal P}$ is a divisor algebraically equivalent to zero.
Then we can write $D\sim n\sigma_1+p_{12}^*C$ for some $n\in\boldz$,
$C\in\Pic(E_1\times E_2)$. Then $n=0$ and $C$
is algebraically equivalent to zero on $E_1\times E_2$. Since $p_{12}^*C\in
{\cal P}$, it follows by restriction to $\sigma_1$ and $\sigma_2$ that
in fact $C=p_1^*F_1$ and $C=p_2^*F_2$ for divisors $F_1$ and $F_2$ on
$E_1$ and $E_2$, respectively. Thus $F_1=F_2=C=0$, so $D=0$ in
$\Pic(Y)$. Thus ${\cal P}$ is a discrete group, as is 
$\Pic^0(\Join(E_1,E_2))$, allowing us to conclude that 
$h^1(\O_{\Join(E_1,E_2)})=0$.
\Box

In fact, in a suitable context, these degenerate ``Calabi-Yau''
threefolds contain a pencil of degenerate Heisenberg invariant
$(1,2n)$-polarized abelian surfaces.  Fix the standard Schr\"odinger
action of $\HHH_{2n}$ on $\P^{2n-1}$, and consider the subgroup
$\HHH''\subseteq \HHH_{2n}$ generated by $\sigma^2$ and $\tau$. Then
$\HHH''/\langle \tau^n\rangle\cong \HHH_n$ acts as the Schr\"odinger
representation of $\HHH_n$ on the subspaces $L_1,L_2\subseteq
\P^{2n-1}$, where $L_1=\set{x_0=x_2=\cdots=x_{2n-2}=0}$, and
$L_2=\set{x_1=x_3=\cdots=x_{2n-1}=0}$.  Let $E\subseteq
L_1\cong\P^{n-1}$ be an $\HHH_n$-invariant elliptic normal curve
(under the Schr\"odinger representation).  Then $\sigma(E)\subseteq
L_2$ is also an $\HHH_n$-invariant elliptic normal curve. Let $0\in E$
denote the origin of $E$. Define, for $\eta\in E$,
$$S_\eta:=\bigcup_{P\in E} \langle P,\sigma(P+\eta)\rangle
\subseteq \Join(E,\sigma(E))\subset\P^{2n-1}.$$

\proposition{newabeldegen} 
\item{\rm (1)} $S_\eta$ is a non-singular elliptic
scroll of degree $2n$ in $\P^{2n-1}$.
\item{\rm (2)} $A_\eta:=S_\eta\cup \sigma(S_\eta)=
S_\eta\cup S_{-\sigma^2(\eta)}$ is
an $\HHH_{2n}$-invariant surface of degree $4n$ and
sectional arithmetic genus $2n+1$. Moreover, 
$A_\eta=A_{\eta'}$ if and only
if $\eta=\eta'$ or $\eta'=-\sigma^2(\eta)$. Thus 
the set $\setdef{A_\eta}{\eta\in E}$ forms a linear 
pencil of surfaces in $\Join(E_1,E_2)$.
\item{\rm (3)} For $\eta\in E$ general, 
there exists a flat family $\A\longrightarrow\Delta$, a point
$0\in\Delta$, along with a $\HHH_{2n}$-invariant embedding
$\A\subseteq \P^{2n-1}_{\Delta}$ such that $\A_0\cong A_\eta$, and
$\A_t$ is a non-singular $(1,2n)$-polarized abelian surface for
$t\in\Delta$, $t\not=0$.

\proof: (1) Non-singularity of $S_\eta$ is straightforward and can be
proved in much the same way as the non-singularity of $\Join(E_1,E_2)$
away from $E_1$ and $E_2$, see \ref{join}. To compute the degree of
$S_\eta$, choose a general hyperplane $H_1\subseteq L_1$, and let
$H=\Join(H_1,L_2)$. Then $H\cap S_\eta =\sigma(E)\cup
l_1\cup\cdots\cup l_n$, where $l_1,\ldots,l_n$ are lines passing
through the points of $E\cap H_1$.  From this we see that $S_\eta$ is
of degree $2n$.

(2) Straightforward.

(3) The elliptic curves $E$ and $\sigma(E)$ are two disjoint
sections of $S_\eta$, thus $S_\eta\cong S=
\P(\O\oplus \L)$, with $\L\in\Pic^0(E)$. This isomorphism can be chosen so that
the section $E$ of $S_{\eta}$ corresponds to the 0-section of
$S$, i.e. the section corresponding to the subbundle $\O$ of
$\O\oplus\L$, while the section $\sigma(E)$ of $S_{\eta}$ corresponds to
the $\infty$-section of $S$, i.e. the section corresponding to
the subbundle $\L$ of $\O\oplus\L$. Then $\O_{S_{\eta}}(E)|_E
=\L$ and $\O_{S_{\eta}}(\sigma(E))|_{\sigma(E)}=\L^{-1}$. Since $E$
maps to an $\HHH_n$-invariant elliptic normal curve in $L_1$, the identification
of $S$ with $S_{\eta}$ is induced by the complete linear system
$|\sigma(E)+n\cdot f_o|$. Now $\sigma(E)+n\cdot f_o\sim E+n\cdot f_p$
for some $p\in E$, and restricting both of these divisors to $E$, we see
that 
$$\L=\O_E(-n([p]-[o])).$$
Now as an abstract surface
$A_\eta$ is a type II degeneration, namely a $2$-cycle of elliptic
ruled surfaces as in the following figure: 
\drawnoname{95}{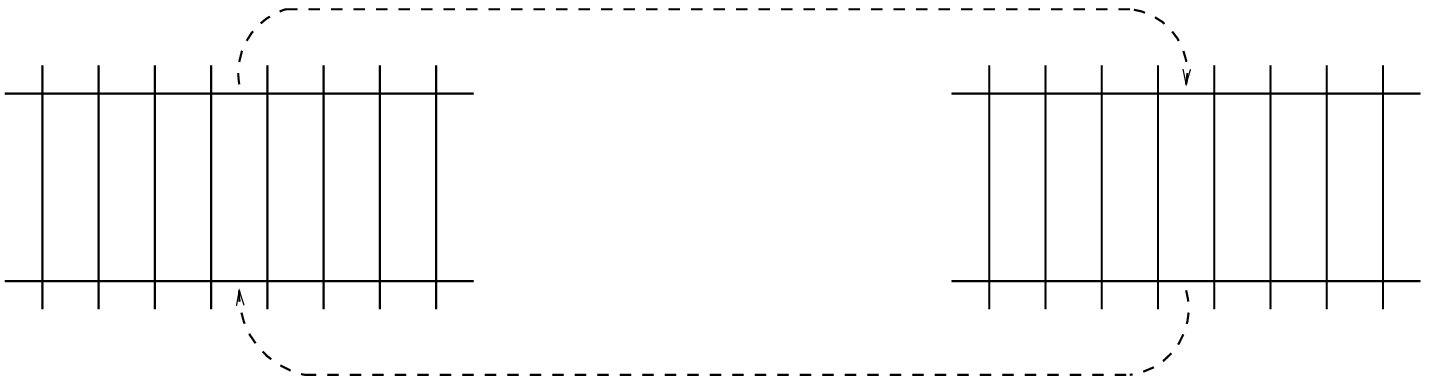}{
\swtext{.375}{1.000}{\swpad{2pt}{2pt}{glue without shift}}
\swtext{.325}{0.730}{\swpad{2pt}{2pt}{$\infty$-section}}
\swtext{.880}{0.730}{\swpad{2pt}{2pt}{$0$-section}}
\swtext{.325}{0.227}{\swpad{2pt}{2pt}{$0$-section}}
\swtext{.880}{0.227}{\swpad{2pt}{2pt}{$\infty$-section}}
\swtext{.116}{0.900}{\swpad{2pt}{2pt}{$S_\eta$}}
\swtext{.766}{0.900}{\swpad{2pt}{2pt}{$\sigma(S_\eta)$}}
\swtext{.375}{0.010}{\swpad{2pt}{2pt}{glue with shift $\sigma^2(2\eta)$}}
}
One starts with the 
surface $S_\eta$ with $E$ and $\sigma(E)$
as its $0$ and $\infty$-sections, respectively, then glues with no shift
the $\infty$-section of $S_\eta$ to the $0$-section $\sigma(E)$ of
the surface $\sigma(S_\eta)\cong\P(\O\oplus \L)$, and then finally
glues with shift $\sigma^2(2\eta)$ the  
$\infty$-section $E$ of $\sigma(S_\eta)$ 
to the $0$-section of the $S_\eta$.  
The restrictions of 
$\O_{S_\eta}(1)$ and $\O_{\sigma(S_\eta)}(1)$ to $E$ coincide
after a shift by $\sigma^2(2\eta)$, which yields
$n([o]-[p])=n([\sigma^2(2\eta)]-[o])$. In other words $p\in E$
is such that $np=-n\sigma^2(2\eta)$ in the group law of $E$, and
$\L=\O_E(-n([p]-[o]))$. As in [HKW], part II, especially Theorem 3.10 and 
Proposition 4.1,
and [DHS], \S3, one may show that such type II degenerations actually occur, 
and in fact all elliptic curves $E$ and all general shifts 
$\eta$ can be realized. 
\Box

\section {1.4} {Moduli of $(1,4)$-polarized abelian surfaces.}

Even though the geometry of the moduli space of $(1,4)$-polarized
abelian surfaces is well understood (see [BLvS]), we wish to partially
review it here as this will lead to another example of a Calabi-Yau
threefold with an abelian surface fibration.  However, we will not
provide many details, leaving further investigations to the interested
reader.

The moduli space $\A_4^{lev}$ was studied in [BLvS]. 
Let $(A,\L)$ be an abelian surface with $\L$ an ample 
line bundle of type $(1,4)$. It is proved in [BLvS] that
if $A$ is sufficiently general, then the morphism $\psi_{|\L|}:A
\longrightarrow \Pthree$ induced by $|\L|$
is birational onto its image, an octic surface. 
Furthermore, [BLvS] give explicitly the defining equation
of such an octic surface. Namely, let $(x_0:\ldots:x_3)$ be the
coordinates on $\Pthree$, on which the Heisenberg group
$\HHH_4$ acts via the Schr\"odinger representation
$$\sigma: x_i\mapsto x_{i-1}\qquad {\rm and}\qquad
\tau:x_i\mapsto \xi^{-i} x_i,$$ 
where $\xi$ is a fixed primitive fourth root of unity. We 
now change coordinates to
$$\eqalign{ z_0= x_0+x_2, &\quad z_2=x_3+x_1,\cr
z_1= x_0-x_2, &\quad z_3=x_3-x_1.\cr}$$
Then the image of $\psi_{|\L|}$ is defined in $\P^3$ 
by the equation  $f=0$,  where
$$f(z_0,\ldots,z_3)=\lambda N {}^t\lambda,$$ with 
$$N={\smallmath
\pmatrix{
z_0^2z_1^2z_2^2z_3^2&0&0&0\cr
0&(z_0^4z_1^4+z_2^4z_3^4)&
(z_0^2z_1^2+z_2^2z_3^2)(-z_0^2z_2^2+z_1^2z_3^2)
&(z_0^2z_1^2-z_2^2z_3^2)(z_0^2z_3^2-z_1^2z_2^2)\cr
0&
(z_0^2z_1^2+z_2^2z_3^2)(-z_0^2z_2^2+z_1^2z_3^2)
&(z_0^4z_2^4+z_1^4z_3^4)
&(z_0^2z_2^2+z_1^2z_3^2)(z_0^2z_3^2+z_1^2z_2^2)\cr
0&
(z_0^2z_1^2-z_2^2z_3^2)(z_0^2z_3^2-z_1^2z_2^2)&
(z_0^2z_2^2+z_1^2z_3^2)(z_0^2z_3^2+z_1^2z_2^2)&
(z_0^4z_3^4+z_1^4z_2^4)\cr}}$$
and for some 
value of the parameter 
$\lambda=(\lambda_0:\ldots:\lambda_3)\in\Pthree$.
 
For a fixed value of $\lambda$ we will denote by 
$A_{\lambda}$ the octic surface defined by 
$\set{f=0}\subset\Pthree$.

Note that 
$(\lambda_0:\ldots:\lambda_3)$
and $(-\lambda_0:\lambda_1:\lambda_2:\lambda_3)$
give the same equation $f=0$, and conversely 
two points $\lambda,\lambda'\in\Pthree$
yield the same equation only if they are related 
in this manner. Thus, if $\boldz_2$ acts on $\Pthree$ via 
$$(\lambda_0:\ldots:\lambda_3)\mapsto
(-\lambda_0:\lambda_1:\lambda_2:\lambda_3),$$ 
then there is a birational map
$$\Theta_4:\A_4^{lev}\rDashto \Pthree/\boldz_2.$$
This latter quotient is isomorphic to the cone over the 
Veronese surface in $\P^5$.

We recall briefly from [BLvS] the structure of the singularities
of $A_{\lambda}\subset\P^3$ for general $\lambda$. The singular locus
of $A_{\lambda}$ is contained in the union of the coordinate
planes $\set{z_0z_1z_2z_3=0}$. In fact, $\Sing(A_{\lambda})\cap
\set{z_i=0}$ is a quartic curve with double points at the three
coordinate vertices on the plane $\set{z_i=0}$. Now
$A_{\lambda}$ has generically ordinary double
points along this curve, with $12$ pinch points in the smooth locus
of the quartic, while $A_{\lambda}$ has quadruple points at the
coordinate vertices (with tangent cone a union of four planes).

\definition{X4l} Let $l\subseteq\Pthree$ be a line. Define
$V_{4,l}^1$ to be the normalization of the hypersurface
$X_l\subseteq\Pthree\times l\subseteq \Pthree\times\Pthree$ given by
the equation
$$\setdef{f(z_0,\ldots,z_3;\lambda_0,\ldots,\lambda_3)=0}
{(\lambda_0,\ldots,\lambda_3)\in l}$$
with $(z_0:\ldots:z_3)$ coordinates on the first $\Pthree$ and
$(\lambda_0:\ldots:\lambda_3)$ coordinates on $l$ on the second $\Pthree$.
Let $X_l\mapright{\psi} Y_l\mapright{}\Pthree$ be the Stein factorization
of the projection $X_l\longrightarrow\Pthree$on the first factor, and define
$V_{4,l}$ to be the normalization of $Y_l$.

Because this case is not of great relevance to the main thrust
of the paper, we state basic properties of these threefolds and only
sketch their proof.

\theorem{CY14} For a general line $l$ in $\Pthree$, 
\item{\rm (1)} $V_{4,l}$ is a double cover of
$\Pthree$ branched over an octic surface $B\subseteq \Pthree$.
\item{\rm (2)} The singular locus of $B$ consists of $148$ ordinary
double points, all of which are contained in $A_{\lambda}$, $\lambda\in l$.
These are of three types: 
\itemitem{\rm (A)}$128$ of them are contained in the
smooth locus of $A_{\lambda}$, 
\itemitem{\rm (B)}$16$ others are contained
in the double point locus of $A_{\lambda}$ for 
$\lambda\in l$ general, 
\itemitem{\rm (C)}The remaining
four occur at the coordinate vertices (at the
quadruple points of $A_{\lambda}$).
\item{\rm (3)} The map $\psi':V^1_{4,l}\longrightarrow V_{4,l}$
induced by the map $\psi:X_l\longrightarrow Y_l$ is a small resolution
of $V_{4,l}$. In particular, $V^1_{4,l}$ is a non-singular Calabi-Yau
threefold which contains a base-point free pencil $\widetilde A_{\lambda}$
of non-singular $(1,4)$-polarized abelian surfaces mapping to the
pencil $A_{\lambda}$, $\lambda\in l$.
\item{\rm (4)} $\chi(V^1_{4,l})=0$ and 
$h^{1,1}(V^1_{4,l})=h^{1,2}(V^1_{4,l})=8$.

\proof: (Sketch) 
(1) Let $l\subset\Pthree$ be parametrized by $\lambda_j=a_{0j}\mu_0
+a_{1j}\mu_1$, $j\in\set{0,1,2,3}$, where $(\mu_0:\mu_1)\in\Pone$, and
$A=(a_{ij})_{0\le i\le 1,0\le j\le 3}$ is a $2\times 4$ matrix. Then
$$f(z_0,\ldots,z_3,a_{00}\mu_0+a_{10}\mu_1,
\ldots,a_{03}\mu_0+a_{13}\mu_1)$$ is a 
quadratic polynomial in $\mu_0,\mu_1$ whose Hessian is the $2\times 2$
matrix $2AN{}^tA$. Thus the discriminant of this quadratic polynomial
is $\det(AN{}^tA)$. By inspection, one finds that every $2\times 2$
minor of the matrix $N$ is divisible by $z_0^2z_1^2z_2^2z_3^2$, and
thus so is $\det(AN{}^tA)$, being a linear combination of minors of
$N$.  Hence one can view $Y_l$ as the double cover of $\Pthree$
branched over $\set{\det(AN{}^tA)=0}$. Since this determinant contains
the factor $z_0^2z_1^2z_2^2z_3^2$, we see $Y_l$ is not normal. Now the
surface
$$B:=\set{(\det AN{}^tA)/z_0^2z_1^2z_2^2z_3^2=0}\subset\Pthree$$ can
be seen to have only $148$ ordinary double points. One may essentially
check this, for a random choice of $l\subset\Pthree$ via a
straightforward computation in {\it Macaulay/Macaulay2\/}.  Thus
$V_{4,l}$, the normalization of $Y_l$, can be obtained by taking a
double cover branched over the octic surface $B$.

The claims in (2) can be verified directly via {\it Macaulay/Macaulay2\/}.

(3) Let $f:V_{4,l}\longrightarrow\Pthree$,
$g:Y_l\longrightarrow\Pthree$ be the double covers.  For
$\lambda\in l$, the surface $A_{\lambda}\times\set{\lambda}
\subseteq X_l$ maps isomorphically to $A_{\lambda}\subseteq \Pthree$
via the first projection, and hence $g^{-1}(A_{\lambda})$ splits as
the union of two surfaces, each being isomorphic to
$A_{\lambda}$. Thus $f^{-1}(A_{\lambda})$ also splits as
$S_{\lambda}\cup S'_{\lambda}$, where $S_{\lambda}$ is the proper
transform of $A_{\lambda}\times\set{\lambda}$.

A local analysis now shows that for general $\lambda$, $S_{\lambda}$ is
non-singular except at the ordinary double points of types (B) and (C).
At nodes of type (B), $S_{\lambda}$ has an improper double point
(i.e. the tangent cone is the union of two planes meeting at a point),
while at the nodes of type (C) the surface $S_{\lambda}$ has tangent cone
a union of four planes meeting at a point. Let $V^2_{4,l}\longrightarrow
V_{4,l}$ be the blow-up of $V_{4,l}$ along $S_{\lambda}$. 
Then $V^2_{4,l}$ is a small resolution of $V_{4,l}$, 
and if all exceptional curves are flopped simultaneously, 
we obtain a small resolution $Z_l\longrightarrow
V_{4,l}$ in which the family of surfaces $S_{\lambda}$ forms 
a base-point free pencil, thus yielding an abelian surface 
fibration $Z_l\longrightarrow\Pone$.
This map along with the natural map $Z_l\longrightarrow\Pthree$ yields a map
$\phi:Z_l\longrightarrow\Pthree\times\Pone$ whose image is clearly 
the hypersurface $X_l$. Hence via the universal property of 
the normalization we get a map $\phi'$ in the diagram
$$
\diagram[tight,width=3em,height=2.5em,midshaft]
Z_l\\
\dTo^{\phi'}&\rdTo^{\phi}\\
V_{4,l}^1&\rTo & X_l\\
\dTo^{\psi'}&&\dTo_{\psi}\\
V_{4,l}&\rTo & Y_l\\
\enddiagram
$$
The morphism $\phi'$ is a birational map between normal varieties, and
it is also clear that it doesn't contract any positive dimensional
components.  Thus it is an isomorphism, and (3) follows, with
$\widetilde A_{\lambda}$ the proper transform of $S_{\lambda}$ in
$V^1_{4,l}$.

(4) Now $\chi(V^1_{4,l})=0$ follows from the fact that the Euler
characteristic of a non-singular double cover of $\Pthree$ branched
along a smooth octic is $-296$, and that $V_{4,l}$ has $148$ ordinary
double points (each ordinary double point increases by two the Euler
characteristic with respect to that of a double solid branched over a
smooth octic).  The calculation of the Hodge numbers may be done in a
general example via {\it Macaulay/Macaulay2\/} using standard techniques
(see [Scho] and [We], or \ref{defect6} for details).
\Box

\remark{misc4} We note here that $X_{4,l}$ can be viewed as a partial
smoothing of the secant variety of an elliptic normal curve in $\Pthree$,
in the sense that through the general point of $\Pthree$ pass precisely
two secants of an elliptic normal curve, so we can think of the
secant variety as a double cover of $\Pthree$. The branch locus is
in fact the union of the four quadric cones containing the elliptic curve.

\section {HM} {Moduli of $(1,5)$-polarized abelian surfaces.}

We review here certain aspects of the geometry of $(1,5)$-polarized
abelian surfaces. We will see in many ways that this case is a
paradigm for many of the higher degree cases.  We first review briefly
the well-known description of $\A_5^{lev}$. See [HKW] and references
therein for proofs and details.

The main result of [HM] is that every $(1,5)$-polarized Heisenberg
invariant abelian surface in $\Pfour$ is the zero locus of a section
of the (twisted) Horrocks-Mumford bundle $\F_{HM}(3)$, and conversely,
when the zero locus is smooth.  Here, $H^0(\F_{HM}(3))$ is a four
dimensional vector space. Thus, if $U\subseteq \A_5^{lev}$ is the open
set consisting of triples $(A,H,\alpha)$, with $\alpha$ the level
structure and $H$ very ample on $A$, we obtain a morphism
$$\Theta_5:U\longrightarrow \P(H^0(\F_{HM}(3)))$$
defined by
$$(A,H,\alpha)\mapsto [s]\in\P(H^0(\F_{HM}(3))),$$ where $s$ is a 
section (uniquely determined up to scalar multiple) 
of $\F_{HM}(3)$ vanishing on $A$. The results of [HKW] 
then show that this morphism extends to a morphism
(denoted exactly as the previous one)
$$\Theta_5:\overline{\A_5^{lev}}\longrightarrow \P(H^0(\F_{HM}(3)))$$ 
where $\overline{A_5^{lev}}$ is the Igusa toroidal 
compactification of $\A_5^{lev}$.
$\Theta_5$ is a birational morphism, and [HKW] gives a complete
description of its structure, as well as a description of
$\A_5^{lev}\subseteq \overline{\A_5^{lev}}$. In theory, for many
of the cases in this paper, one could produce a similar fine structure
theory, and thus provide a biregular description of $\overline{\A^{lev}_n}$,
but we shall not even begin to attempt this. See also
the introduction for comments regarding the $(1,7)$ and
$(1,11)$-polarizations.

\definition{HMquintic} Paraphrasing {\rm [Moo]}, 
(see also {\rm [Au], [ADHPR1]} and {\rm [ADHPR2]})
we define for $y\in\Pfour$ the Horrocks-Mumford quintic 
$X_{5,y}:=
\set{\det(M_5'(x,y))=0}\subseteq\Pfour,$ 
whenever this determinant does not vanish identically,
where $M_5'(x,y)=(x_{3(i+j)}y_{3(i-j)})_{i,j\in\boldz_5}$
as in \ref{prelim}.

This is not the usual definition of the Horrocks-Mumford quintics:
more standard is to choose two independent sections $s,s'\in
H^0(\F_{HM}(3))$ and consider the vanishing locus of $s\wedge s'\in
H^0(\wedge^2\F_{HM}(3))
\cong H^0(\O_{\Pfour}(5))^{\HHH_5}$. However, these two definitions coincide.
An argument is given in Remark 4.1 and preceding discussion of [ADHPR2].
Summarizing that argument, we define a rational map
$$\overline{\Theta}:\Pfour\rDashto\P(H^0(\O_{\Pfour}(5))^{\HHH_5})
\cong\P(\wedge^2 H^0(\F_{HM}(3)))$$ by taking 
$y\in\Pfour$ to the $\HHH_5$-invariant quintic $\det(M'_5(x,y))$. This
map is defined outside of the so-called Horrocks-Mumford lines, the
$\HHH_5$-orbit of $\Pone_-$. The image of this map is the Pl\"ucker
quadric of decomposable tensors in 
$\P(\wedge^2 H^0(\F_{HM}(3)))$. 
This is worked out explicitly for instance in Remark 4.1 of [ADHPR2].

\theorem{HMwedge} For a general $y\in\Pfour$,
\item{\rm (1)} $X_{5,y}$ is a 
quintic hypersurface in $\Pfour$ whose singular locus
consists of $100$ ordinary double points.
\item{\rm (2)} There is a small resolution $X^1_{5,y}\longrightarrow X_{5,y}$
such that $X^1_{5,y}$ is a Calabi-Yau threefold with
$\chi(X^1_{5,y})=0$, $h^{1,1}(X^1_{5,y})=h^{1,2}(X^1_{5,y})=4$.
In addition, $X^1_{5,y}$ is fibred in $(1,5)$-polarized abelian surfaces.
\item{\rm (3)} $X^1_{5,y}$ also contains a pencil of abelian surfaces with
a polarization of type $(2,10)$, blown up in $25$ points.

\proof: These are all well-known. Sketching the ideas here, for (1),
two general sections of the Horrocks Mumford bundle $\F_{HM}(3)$
vanish along two smooth abelian surfaces that meet transversally in
$100$ points, the nodes of the corresponding quintic (see [HM], [Au],
[Hu2] for details). For (2), first note that if $A\subseteq
\Pfour$ is an abelian surface  and $y\in A$, then $M'_5(x,y)$ has
rank at most $4$ on $A$ by [GP1], Corollary 2.8. Thus $A\subseteq
X_{5,y}$. Blowing up $A$ produces a small resolution, and flopping the 100
exceptional curves, we obtain $X^1_{5,y}$, in which $A$ moves in a 
base-point free pencil, by \ref{ABpencil}.
The invariants of $X^1_{5,y}$ are well-known, see [Au] for details.

For (3), it is well-known that if  $A\subseteq X_{5,y}$ as above,
and $X'$ is another general quintic hypersurface containing $A$, then
$$X_{5,y}\cap X'=A\cup A',$$
where $A'\subseteq \Pfour$ is a non-singular surface of degree $15$,
in fact an abelian surface with a $(2,10)$-polarization blown up in 
$25$ points (Ellingsrud-Peskine unpublished, see [Au] or [ADHPR1],
[ADHPR2] for details). It follows from \ref{ABpencil} that 
there is a pencil of such surfaces on $X_{5,y}$.\Box

\remark{misc5} (1) The Horrocks-Mumford quintics $X_{5,y}$ can be viewed
as partial smoothings of the secant varieties of elliptic normal curves in
$\Pfour$, as [GP1] Theorem 5.3 or [Hu1], pg. 109 shows that such
secant varieties are also Horrocks-Mumford quintics.\hfill\break
(2) The K\"ahler and moving cones of $X^1_{5,y}$ are well-studied: see
[Bor1], [Bor2], [Scho], [Scho2] and [Fry]. 
In particular there are an infinite number of
minimal models, and $X^1_{5,y}$ contains an infinite number of
pencils of (birationally) abelian surfaces, of both types.
\medskip

We now collect a number of results about Horrocks-Mumford quintics and
$(2,10)$-polarized abelian surfaces we will need later which do not
appear to be in the literature.  We first describe a certain family of
Horrocks-Mumford (HM-) quintic hypersurfaces $X_{5,y}\subset\Pfour$, 
where the parameter point $y$ lies in $\Ptwo_+\subseteq\Pfour$. 
This result will be needed for understanding the singularity structure of
Calabi-Yau threefolds arising in the $(1,10)$ case (see \ref{CY10}):

\proposition{symmoore} Let $M_5'(x,y)={(x_{3(i+j)}y_{3(i-j)})}_{i\in\boldz_5}$,
and let $X_{5,y}$ be the (symmetric) HM-quintic hypersurface given by 
$\set{\det M_5'(x,y)=0}\subset\Pfour(x)$, for a  parameter point 
$y\in\Ptwo_+$ (here and in the proof of the proposition $\pm$ are again
with respect to the Heisenberg involution $\iota$ acting in $\Pfour$). 
\item{\rm (1)} If $y\in B\subset\Ptwo_+$, where $B$ is the Brings curve 
(the curve swept by the non-trivial $2$-torsion points of $\HHH_5$
invariant elliptic normal curves in $\P^4$, see {\rm [BHM]}), 
then $X_{5,y}$ is the secant variety of an elliptic normal curve in $\Pfour$.
\item{\rm (2)} If $y\in C_+\subset\Ptwo_+$, 
where $C_+=\set{y_0^2+4y_1y_2=0}$ is the 
modular conic (cf. {\rm [BHM]}), then $X_{5,y}$ is the trisecant variety 
of an elliptic quintic scroll in $\Pfour$.
\item{\rm (3)} For a general $y\in\Ptwo_+\setminus(C_+\cup B)$, the singular 
locus of $X_{5,y}$ is the union of two elliptic curves of degree $10$, 
meeting along the $\HHH_5$-orbit of the parameter point $y$. See {\rm (5)}
below for the nature of the singularities. 
\item{\rm (4)} Let $\widetilde X_{5,y}\subseteq \Pfour(x)\times\Pfour(z)$
be defined by
$$\widetilde X_{5,y}=\setdef{(x,z)\in\Pfour\times\Pfour}{M'_5(x,y)z=0}.$$
Then for general $y\in\Ptwo_+\setminus(C_+\cup B)$,
$\widetilde X_{5,y}$ has $50$ ordinary double points, 
and $\Sing(\widetilde X_{5,y})$
maps to the Heisenberg orbit of $y$ under projection to $\Pfour(x)$.
\item{\rm (5)} For general $y\in\Ptwo_+\setminus (C_+\cup B)$,
$X_{5,y}$ has only $cA_1$ singularities away from the Heisenberg orbit
of $y$, while at a point of the Heisenberg orbit of $y$,
$X_{5,y}$ has a $cA_3$ singularity.
\item{\rm (6)} If $p_1,p_2:\Pfour\times\Pfour\rightarrow\Pfour$
are the two projections, then 
$X_{5,y}':=p_2(\widetilde X_{5,y})\subseteq\Pfour$
is a $\HHH_5$-invariant quintic which, for general $y\in\Ptwo_+\setminus
(C_+\cap B)$, has singular locus a union of two elliptic quintic
normal curves in $\Pfour$. Furthermore, $X_{5,y}'\subset\P^4(z)$ is defined by the equation
$\set{\det L(z,y)=0}$, where
$L(z,y)$ is the $5\times 5$ matrix given by
$$L(z,y):={(z_{2i-j}y_{i-j})}_{i,j\in\boldz_5}.$$

\proof: 
First we note that (1) follows immediately from 
[GP1], Theorem 5.3. On the other hand, (2) follows from
[ADHPR1], Proposition 24.

For the rest,
recall that $\Gamma=\wedge^2H^0(\F_{HM}(3))={H^0(\O_{\Pfour}(5))}^{\HHH_5}$
is the $6$-dimensional space of $\HHH_5$-invariant quintics in
$\P^4$. The decomposable vectors in $\Gamma$ correspond to the
Horrocks-Mumford quintics, namely the quintic hypersurfaces in
$\Pfour$ whose equations are the determinants of matrices of type
$M'_5$, see \ref{HMquintic} and the discussion thereafter.  We will
need the following standard facts concerning HM-quintics and
$\HHH_5^e=\HHH_5\rtimes\langle\iota\rangle$-invariant elliptic
quintic scrolls in $\Pfour$, most of which can be found in [BHM]:

\item{\rm (I)} There exists a 1-dimensional family $X_t$, 
$t\in\Pone_-$, of $\HHH^e_5$-invariant elliptic quintic scrolls 
in $\Pfour$; the smooth ones correspond to $t\in\Pone_-\setminus
\set{(0:1),(1:0),((1\pm\sqrt{5})\xi^k:2), k\in\boldz_5}$, the
singular are cycles of planes.
The ruling of a smooth elliptic scroll $X_t$ over the origin of 
the base curve maps to  a line $l\subset\Ptwo_+$, which is tangent to 
the conic $C_+$ at a point which corresponds to the point 
$t\in X(5)\cong C_+$ via the standard identification 
[BHM], or [ADHPR2], Proposition 4.3. 
The scroll $X_t$ is embedded in $\Pfour$ by the linear
system $|C_0+2l|$, where $C_0$ is the (unique) section of 
the scroll with self-intersection $1$ meeting $\Pone_-$.
\item{\rm (II)} For any smooth elliptic scroll $X_t$, 
there exist exactly three pairs $(E_i,\tau_i)$, $i =0,1,2$, of
$\HHH^e_5$-invariant elliptic normal curves in $\P^4$ and $2$-torsion
points, such that $X_t$ is the $\tau_i$-translation scroll of $E_i$,
see the discussion after \ref{recallsec} and [Hu2] for exact definitions.
\item{\rm (III)} The linear system $|2E_i|=|-2K_{X_t}|=|4C_0-2l|$ is a base-point
free elliptic pencil, whose only singular fibres are the double curves
$2E_i$.
\item{\rm (IV)} The linear system $|C_0+2E_i|=|5C_0-2l|$ on $X_t$, which has
as base locus the three origins (on $\Pone_-$) of the elliptic normal 
curves $E_i$, 
defines a $2:1$ rational map from $X_t$ onto $\Ptwo$. The map is branched along
$l$ and the three exceptional lines (over the base points), it contracts
$C_0$ and the elliptic normal curves $E_i$ to $p$ and $p_i$, respectively, 
and maps the other elements of the pencil $|-2K_{X_t}|$ to the pencil 
of lines through $p$ (cf. [BHM], Propositions 5.4 and 5.5).
\item{\rm (V)} By [BHM], proof of Proposition 5.9, the restriction of $\Gamma$,
the space of $\HHH_5$-invariant quintics, to a $\HHH^e_5$-invariant
elliptic quintic scroll $X_t\subset\P^4$ is always
$3$-dimensional. Furthermore, the kernel of this restriction is
$s_t\wedge H^0(\F_{HM}(3))$, where $s_t\in H^0(\F_{HM}(3))$ is the
unique section of the Horrocks-Mumford bundle vanishing on a double
structure on the elliptic scroll $X_t$.  The sections $s_t\in
H^0(\F_{HM}(3))$ vanishing on a double structure on the elliptic
scroll $X_t$ are parameterized by a smooth rational sextic curve
$C_6\subset\Pthree=\P(H^0(\F_{HM}(3)))$ (cf. [BM], 1.2).
\item{\rm (VI)} The linear system $\Gamma$ induces a rational map 
$\Theta:\Pfour\rDashto\Omega\subset\Pfive$, where $\Omega$ is a smooth
quadric, which is $100:1$ and is not defined exactly on the so called
Horrocks-Mumford lines (the $\HHH_5$ orbit of $\Pone_-$).  In this
setting, Horrocks-Mumford quintics correspond to pullbacks via
$\Theta$ of the tangent hyperplanes to the Pl\"ucker quadric $\Omega$.
\item{\rm (VII)} $\Theta$ restricted to $X_t$ factors as 
$\Theta:X_t\rightarrow X_t/\boldz_5\times\boldz_5\cong X_t
\rDashto\Omega$, and the latter map is induced by
the linear system $|C_0+2E_i|$, and thus by (IV) above,
maps the scroll onto a linear subspace $\Ptwo\subset\Omega\subset\Pfive$.

\noindent
Next we identify the quintics in part (3) of the statement of 
\ref{symmoore}:

\lemma{} Symmetric HM-quintics $X_{5,y}\subset\P^4$, for a parameter
$y\in\Ptwo_+\setminus C_+\cup B$, correspond to wedge products
$s_1\wedge s_2$, where $s_i\in H^0(\F_{HM}(3))$, $i=1,2$, are two
sections of the Horrocks-Mumford bundle each vanishing on a
$\HHH^e_5$-invariant elliptic quintic scroll $X_i\subset\P^4$.  Such
quintic hypersurfaces can also be characterized as the unique quintics
in $\Pfour$ containing the union of the two elliptic quintic scrolls
$X_i$: If $l_i$ are the rulings over the origin of $X_i$, $i=1,2$,
then we may take as parameter $y$ of the matrix $M'_5(x,y)$ the point
$\set{y}=l_1\cap l_2\in\P^2_+$. In particular, this allows us to
identify $\Sym^2(C_6)$ with $\Ptwo_+$.

\proof: Let $y\in\Ptwo_+\setminus (C_+\cup B)$, 
and let $l_1$ and $l_2$ be the two
tangent lines to the modular conic $C_+$ that pass through the point
$y$.  By fact (I) above $l_1$ and $l_2$ are each rulings over the
origin for two distinct $\HHH_5^e$-invariant quintic elliptic scrolls
$X_i\subset\P^4$, $i=1,2$.  On the other hand, 
for a fixed point $x\in B$, 
the quintic $\set{\det M'_5(x,y)=0}\subset\Pfour(y)$
is the secant variety of the  elliptic normal curve in $\P^4(y)$ passing
through $x$. Now if $E_1^i, E_2^i, E_3^i\subseteq X_i$ are the three
elliptic normal curves of fact (II) above, then $E^i_j\cap l_i$ consists of two
distinct points, and thus $B\cap l_i=\bigcup_{j=1}^3 (E^i_j\cap l_i)$ consists
of six distinct points (see [BHM], \S 6). Thus $l_i$ is a secant
to each $E^i_j$, and thus
$\det M_5'(x,y)$ vanishes at each $x\in B\cap l_i$.
By B\'ezout's theorem, it
follows that $X_{5,y}$ must vanish along both rulings $l_1$ and
$l_2$. Further, by Heisenberg invariance, $X_{5,y}$ must then vanish on
the $\HHH_5$ orbits of $l_1$ and $l_2$, and thus on both elliptic
scrolls $X_1$ and $X_2$ again by B\'ezout's theorem.  The first claim
in the lemma follows now from fact (V) above.  For the second claim,
observe that by fact (VII) above, the map $\Theta:\P^4\rDashto\Omega$
maps the two scrolls $X_i$ onto two planes in $\Pfive$ meeting at a
point, and thus spanning a unique hyperplane in $\Pfive$. This
concludes the proof of the lemma.\Box

\noindent {\it Proof of \ref{symmoore} continued:} Suppose that
$y\in\Ptwo_+\setminus (C_+\cup B)$, and thus that $X_{5,y}=\set{s_1\wedge
s_2=0}\subset\P^4(x)$, with $s_i\in H^0(\F_{HM}(3))$ vanishing doubly
on elliptic quintic scrolls $X_i\subset\P^4(x).$ The elliptic quintic
scrolls $X_1$ and $X_2$ meet only in $25$ points, namely the $\HHH_5$
orbit of $y\in\Ptwo_+$. Moreover, since $\rank M'_5(y,y)\le 3$, it
follows that all the $\HHH_5$ translates of $y$ are in the singular
locus of the quintic $X_{5,y}$.

By fact (III) above, on each scroll $X_i$, the linear system $|-2K_{X_i}|$ is
a base-point free pencil, with general member a smooth elliptic curve
of degree $10$. We denote by $D_i$ the unique degree $10$ elliptic
curve in the pencil $|-2K_{X_i}|$ which passes through the parameter
point $y$. Each curve in the pencil $|-2K_{X_i}|$ is Heisenberg
invariant, since at least three curves in the pencil 
(the doubled elliptic quintic
curves) are Heisenberg invariant, and only the identity
automorphism on $|-2K_{X_i}|\cong\Pone$ has $\ge 3$ fixed points.
Thus $D_1\cap D_2=X_1\cap X_2$ is the $\HHH_5$ orbit of
$y$.

We show now that $\rank M'_5(x,y)\le 3$  for all $x\in D_1\cup D_2$,
and thus that $X_{5,y}$ is singular along the union of these
two elliptic curves. The space of Heisenberg invariant quintics
containing the elliptic curve $D_i$ (the elliptic scroll $X_i$)
is $4$-dimensional (respectively, $3$-dimensional), so $D_i$ is
linked on $X_i$ to a $\HHH_5$-invariant curve $G_i$ of degree $15$.
The elliptic curve  $G_i$ is a section of the scroll $X_i$,
described in [ADHPR2], Proposition 4.10 (iii), and is the unique 
$\HHH_5\rtimes\langle\iota\rangle$-invariant curve of class
$C_0+12l$ on $X_i$. By [ADHPR2], Proposition 4.12 (ii), 
for a fixed point $z\in G_i$, not
on a Horrocks-Mumford line, the quintic 
$\set{\det M'_5(x,z)=0}\subset\Pfour(x)$
is the trisecant variety of the elliptic scroll $X_i$. It then follows,
as in [ADHPR1], Proposition 24, that $\rank M'_5(x,y)
\le 3$ for all $x$ on the unique curve of numerical equivalence class
$4C_0-2l$ passing through $z$. However, if $z\in D_i$, this curve is
precisely $D_i$.
On the other hand the matrices $M'_5(x,y)$ and
$M'_5(y,x)$ coincide up to a permutation of their columns, and so
do their collection of $4\times 4$ minors. Thus it follows that
$\rank M'_5(x,y)\le 3$  for all $x\in D_i\cap G_i$. Now on the scroll
$X_i$ we have $D_i\cdot G_i=50$, and since both these curves are
$\HHH_5$-invariant, 
$D_i\cap G_i$ must consist of at least $25$ distinct points.
By [ADHPR2], Proposition 4.8, 
$G_1\cap G_2=\emptyset$, and since $\rank M'_5(x,y)\le 3$ 
for all $x\in D_1\cap D_2$, it follows that
each $4\times 4$ minor of $M'_5(x,y)$ vanishes in at 
least $50$ points along $D_i$, and thus by B\'ezout's theorem,
vanishes identically on the curve $D_i$. 
It follows that $\rank M'_5(x,y)\le 3$ for all $x\in D_1\cup D_2$, 
which means that the quintic $X_{5,y}\subset\P^4$ is singular along
the union of the two elliptic curves $D_i$. 

One checks now in {\it Macaulay/Macaulay2\/} that for a general
$y\in\Ptwo_+\setminus (C_+\cup B)$ the quintic $X_{5,y}$ 
is singular only along the union $D_1\cup D_2$.  

To show (4) and (5), note that it is well-known that for general
$y\in\Pfour$, $\widetilde X_{5,y}\subseteq \Pfour\times\Pfour$ has
precisely 50 ordinary double points. (See for example [Au] or
[Bor1].) Thus, if this is also the case for one $y\in \Ptwo_+$,
it will be the case for general such $y$. But this can be checked using
{\it Macaulay/Macaulay2\/}. Now $p_1:\widetilde X_{5,y}\rightarrow X_{5,y}$
is an isomorphism away from $p_1^{-1}(\Sing(X_{5,y}))$, since $\rank M'_5(x,y)
\le 3$ only for points of $\Sing(X_{5,y})$. 
Furthermore, $p_1^{-1}(\Sing(X_{5,y}))
\rightarrow \Sing(X_{5,y})$ is a $\Pone$-bundle, since $\rank M'_5(x,y)=3$
on $\Sing(X_{5,y})$. Since generically $p_1$ resolves the singularities of
$X_{5,y}$, $X_{5,y}$ must generically have $cA_1$ singularities. Because
$\widetilde X_{5,y}$ has only ordinary double points, it is then clear that
$X_{5,y}$ has only $cDV$ singularities, and by deformation theory
of du Val singularities, $y$ must be at least $cA_3$, and at least two
ordinary double points of $\widetilde X_{5,y}$ are sitting over $y$. By
Heisenberg invariance, there are precisely two ordinary double points
over each point in the Heisenberg orbit of $y$, and thus each
of these points is a $cA_3$ point, while away from the Heisenberg orbit
of $y$, points of $\Sing(X_{5,y})$ are $cA_1$ points.

Finally, to show (6), first note that the equation $M_5'(x,y)z=0$ of
$\widetilde X_{5,y}$ is equivalent to the equation ${}^tL(z,y)x=0$, given
that $y\in\Ptwo_+$. Thus
the equation of $X_{5,y}'=p_2(\widetilde X_{5,y})$ 
is given by $\set{\det L(z,y)=0}$.
Now let $N(x,y)$ be a $4\times 5$ matrix of linear forms whose rows
are general linear combinations of the rows of $M_5'(x,y)$, and let
$$\Delta=\setdef{x\in\Pfour}{\rank N(x,y)\le 3}.$$
Then it is well-known (see [Au] or [Bor1]) that the linear system
$|4H-\Delta|$ on $X_{5,y}$ induces a birational map 
$\varphi:X_{5,y}
\rDashto X_{5,y}'$ which lifts to the projection 
$p_2:\widetilde X_{5,y}\rightarrow
X_{5,y}'$. Now $X_{5,y}$ contains the elliptic quintic scrolls $X_1$ and $X_2$.
Each scroll $X_i$ is embedded via $|C_0+2l|$ by fact (I) above,
and hence $\varphi|_{X_i}$ is induced by a subsystem of 
$|4C_0+8l-(\Delta\cap X_i)|$. But $\Delta$ certainly includes the curve
$D_i$, which is of class $-2K_{X_i}=4C_0-2l$. Furthermore, for
a general $(1,5)$-polarized abelian surface $A\subseteq X_{5,y}$, 
$\Delta\cap A$ is
degree $20$. It then follows from the intersection theory on
a resolution of $X_{5,y}$ that
that $\Delta\cap X_i=D_i\cup l^1_i\cup\cdots\cup l^5_i$, where
$l^j_i$ are  
lines of the ruling on $X_i$, and thus $\varphi|_{X_i}$ is
induced by a linear system numerically equivalent to $5l$. Thus 
$\varphi|_{X_i}$ maps $X_i$ to an elliptic normal quintic curve
in $X_{5,y}'\subseteq \Pfour$, and $\Sing(X_{5,y}')$ contains  at least two
elliptic normal curves. However, checking for general
$y\in\Ptwo_+\setminus (C_+\cup B)$, one finds via {\it Macaulay/Macaulay2\/}
precisely this singular locus, so this describes the singular locus
for $X_{5,y}'$ for general $y\in P^2_+$.
\Box

Next we study $(2,10)$-polarized abelian surfaces in $\Pfour$.

\proposition{nonmin15} Let $A$ be a general abelian surface with an
ample line bundle $\L=(\L')^{\otimes 2}$ of type $(2,10)$, and let 
$x_0\in A$. Let $\pi:\tilde A\rightarrow A$ be the blow-up of $A$
at the $25$ points of the set $x_0+K(\L')$, and suppose there exists elements
$f_0,\ldots,f_4\in H^0(\tilde A,\pi^*\L\otimes\omega_{\tilde A}^{-1})$
inducing a morphism 
$$f=(f_0,\ldots,f_4):\tilde A\rightarrow\Pfour.$$ Suppose
furthermore this map is equivariant with respect to the action of
$K(\L')$ on $\tilde A$ and the Schr\"odinger action of $\HHH_5$ on $\Pfour$,
and is also equivariant with respect to negation on $\tilde A$ and
the Heisenberg involution $\iota$ on $\Pfour$. 
Then, possibly after changing the origin on $A$,
we can take $x_0=0$, $\L$ a line bundle of characteristic $0$ with respect
to some decomposition. Furthermore, $f$ is an embedding, and  for
any $y\in f(\tilde A)$ the quintic hypersurface $\set{\det M'_5(x,y)=0}$ 
contains $f(\tilde A)$.

\proof: The fact that negation on $A$ lifts to negation on 
$\tilde A$ tells us that $x_0$ is a two-torsion point. Furthermore,
$(-1)^*\L\cong\L$, so $\L$ is a symmetric line bundle. Since $x_0$
is two-torsion and $A_2\subseteq K(\L)$, we have $t_{x_0}^*\L\cong \L$, and
in particular if we change the origin of $A$ to be $x_0$, $\L$ will
still be a symmetric line bundle. Indeed, if $(-1)'$ denotes negation
on $A$ with origin $x_0$, then $(-1)'^*=t_{-x_0}^*\circ (-1)^*\circ
t_{x_0}^*$, so $(-1)'^*\L\cong (-1)^*\L\cong\L$. In addition, the
action of $(-1)'$ on sections of $\L$ is the same as that of the negation
$(-1)$, so changing the origin to $x_0$ does not affect the hypotheses of
the proposition.

Now let $H\subseteq\Pfour$ be a hyperplane containing $\Ptwo_+$, let
$\tilde D=f^*H$ be the corresponding divisor on $\tilde A$, 
and set $D=\pi_*\tilde D$.
Let
$$A_2^-(D)=\setdef{x\in A_2}{\mult_x(D)\equiv 1\mod 2}.$$
By [LB], IV, (7.6), if $\L$ is of characteristic zero with respect to some
decomposition, then 
$\#A_2^-(D)$ is $0$ or $16$, while if $\L$ is not of characteristic zero,
then $\#A_2^-(D)=8$. We claim the latter does not occur. 
Indeed, since $\tilde D^2=15$, we expect $f^{-1}(\Ptwo_+)$ 
to consist of $15$ points with $D$
vanishing to order $1$ on each point, ruling out the latter case. To
show that $f^{-1}(\Ptwo_+)$ in fact consists of $15$ distinct points,
we proceed as follows. We have $\pi(f^{-1}(\Ptwo_+))\subseteq A_2$,
so $f^{-1}(\Ptwo_+)$ consists of $2$-torsion points and possibly a number 
of copies of $E_1=\pi^{-1}(0)$. Now if $p$ is an isolated point in
$f^{-1}(\Ptwo_+)$, then $\iota$ acts on the tangent plane $T_{A,p}$
to $A$ at $p$ as negation. In particular, $f_*T_{A,p}$ must intersect
$\Ptwo_+$ transversally, so that $p$ occurs with multiplicity one
in $f^{-1}(\Ptwo_+)$. Similarly, if $E_1\subseteq f^{-1}(\Ptwo_+)$, then
$E_1$ must occur with multiplicity one. 
Thus either $E_1\not\subseteq f^{-1}(\Ptwo_+)$, and then 
$f^{-1}(\Ptwo_+)$ consists of $15$ distinct points,
or $E_1\subseteq f^{-1}(\Ptwo_+)$, and we also have a residual
$(\tilde D-E_1)^2=12$ points. In either case, we have $\#A_2^-(D)>8$,
and hence $\#A_2^-(D)=16$ and $\L$ is of characteristic $0$.

Now consider the map $\Theta:\Pfour\rDashto\Omega\subseteq \Pfive$
introduced in fact (VI) of the proof of \ref{symmoore}. Recall from
[Au] (see also [ADHPR1], Remark 37) that if $L_1$ is
an $\alpha$-plane in $\Omega$, then $\Theta^{-1}(L_1)$ is
a $(1,5)$-polarized abelian surface union the base locus of $\Theta$. 
Moreover, if $L_1$ is linked to a $\beta$-plane $L_2$ via a hyperplane
section of $\Omega$, then $\Theta^{-1}(L_2)$ must be a degree $15$
non-minimal $(2,10)$-polarized abelian surface 
(see [ADHPR1], pg. 898). In particular, there is a three-dimensional
family of such $\HHH_5$ and $\iota$-invariant surfaces. These
surfaces are embedded in $\Pfour$ by a linear system of the type given in the
hypotheses of this Proposition, and since these surfaces are linearly
normal, the embedding is given by a complete linear system
$|H^0(\pi^*\L\otimes\omega_{\tilde A}^{-1})|$, where 
$\L$ is of characteristic zero with respect to some decomposition. 

Now the moduli space of $(2,10)$-polarized abelian surfaces along with
a choice of symplectic basis is the Siegel upper half-space ${\cal H}_2$.
This is three dimensional. There is a universal family $\A\rightarrow
{\cal H}_2$ with a zero section, and a line bundle $\L$ on $\A$ such
that for any $Z\in{\cal H}_2$, $\L|_{\A_Z}$ is the line bundle of type $(2,10)$
and of characteristic zero with respect to the decomposition on $\A_Z$
induced by the choice of symplectic basis (see [LB], Chapter 8, (7.1)).
Furthermore, we can blow up the submanifold $\bigcup_{Z\in{\cal H}_2} 
2K(\L|_{\A_Z})\subseteq \A$ to obtain a three-dimensional family
of non-minimal $(2,10)$-polarized abelian surfaces 
$\tilde A\rightarrow {\cal H}_2$,
along with the line bundle 
$\pi^*\L\otimes\omega_{\tilde\A/{\cal H}_2}^{-1}$.
Here $\pi:\tilde \A\rightarrow\A$ denotes the blow-up map. 
This defines the universal
family of polarized non-minimal abelian surfaces of precisely the sort
we are interested, and every such surface appears in this family.
Now the point is that this is a three-dimensional family, and further for
each $Z\in{\cal H}_2$, there are only a finite number of bases
of $H^0(\pi^*\L\otimes\omega_{\tilde\A/{\cal H}_2}^{-1}|_{\A_Z})$ for
which the induced map $f:\tilde\A_Z\rDashto\Pfour$ is equivariant
with respect to translation by elements of $2K(\L|_{\A_Z})$ on $\tilde\A_Z$ and
the Schr\"odinger representation of $\HHH_5$ on $\Pfour$. Indeed,
all such bases are related by an element of $\SL_2(\boldz_5)$ acting on
$H^0(\pi^*\L\otimes\omega_{\tilde\A/{\cal H}_2}^{-1}|_{\A_Z})$. 
Since we already have a
three-dimensional family of such abelian surfaces embedded in $\Pfour$, these
two families must coincide. Thus for a general $Z\in{\cal H}_2$, 
$|H^0(\pi^*\L\otimes\omega_{\tilde\A/{\cal H}_2}^{-1}|_{\A_Z})|$ 
is very ample, while the embedded surface
$f(\tilde\A_Z)$ is of the form $\Theta^{-1}(L_2)$ for some $\beta$-plane
$L_2\in\Omega\subset\Pfive$. 
Now if $y\in f(\tilde\A_Z)$, and $H$ is a hyperplane in $\Pfive$
tangent to $\Omega$ at $\Theta(y)$, then $\Theta^{-1}(H)$ is the 
Horrocks-Mumford quintic $\{\det M_5'(x,y)=0\}$, and clearly
$f(\tilde\A_Z)\subseteq \Theta^{-1}(H)$, since $L_2\subseteq H$.\Box

\section{1.6} {Moduli of $(1,6)$-polarized abelian surfaces}

We will show that the general abelian surface $A\subset\P^5$ with a
$(1,6)$ polarization is determined by the cubics containing it, and
this will allow us to define a rational map $\Theta_6$, essentially
taking $A$ to the set of cubics containing it.  In fact, we will find out in
addition that through a general point $y$ of $\Pfive$ passes exactly
one $(1,6)$-polarized abelian surface! The strategy for determining
this abelian surface is to consider the cubics passing through the
$\HHH_6$-orbit of a given $y\in\Pfive$.  These cubics (or a specific
subspace of these cubics) will also contain the unique abelian surface
passing through $y$.

We first need to discuss the representation theory of $\HHH_6$ acting
on the vector space $H^0(\O_{\Pfive}(3))$, the space of cubic
forms on $\Pfive$. There are no Heisenberg invariant cubic forms
on $\Pfive$; however if $\HHH'\subseteq\HHH_6$ is the subgroup
generated by $\sigma^2$ and $\tau^2$, then there are $\HHH'$-invariant
cubic forms. We denote the space of such forms by 
$H^0(\O_{\Pfive}(3))^{\HHH'}$.
It is easy to see that 
$H^0(\O_{\Pfive}(3))^{\HHH'}$ has as a basis $f_0,\ldots,f_3,\sigma f_0,
\ldots,\sigma f_3$, where
$$\eqalign{f_0&=x_0^3+x_2^3+x_4^3\cr
f_1&=x_1^2x_4+x_3^2x_0+x_5^2x_2\cr
f_2&=x_1x_2x_3+x_3x_4x_5+x_5x_0x_1\cr
f_3&=x_0x_2x_4.\cr}$$ 
Furthermore,
$H^0(\O_{\Pfive}(3))^{\HHH'}$ is a representation of $\HHH_6$, which
splits up into four isomorphic representations
$$H^0(\O_{\Pfive}(3))^{\HHH'}
\cong \bigoplus_{i=0}^3 \langle f_i,\sigma f_i \rangle.$$
We will identify $H^0(\O_{\Pfive}(3))^{\HHH'}$
with $V_0\otimes W$, where $V_0$ is a two dimensional 
representation of $\HHH_6$ and $W$ is a four
dimensional complex vector space with basis $e_0,e_1,e_2,e_3$ so that
$V_0\otimes\langle e_i\rangle=\langle f_i,\sigma f_i\rangle$. 

The importance of $H^0(\O_{\Pfive}(3))^{\HHH'}$ comes from 
the following key lemma:

\lemma{invar3} Let $A\subseteq\Pfive$ be a $(1,6)$-polarized 
abelian surface which is invariant
under the action of $\HHH^e_6=\HHH_6\rtimes\langle\iota\rangle$. 
Then $\dim H^0(\I_A(3))^{\HHH'} \ge 4$. 

\proof: This is a very similar argument to that in [HM], bottom of
page $76$. We consider the restriction map
$$H^0(\O_{\Pfive}(3))^{\HHH'}\longrightarrow H^0(\L^3)^{\HHH'},$$
where $\L$ is the line bundle on $A$ inducing the embedding. Let
$A'=A/2K(\L)$, and $\pi:A\longrightarrow A'$ the quotient map. Then
$2K(\L)$ acts on $\L^3$, and there exists a bundle $\M$ on $A'$ such
that $\L^3=\pi^*(\M)$. Furthermore, $H^0(\L^3)^{\HHH'}=H^0(\M)$. Now
$c_1(\M)^2=c_1(\L^3)^2/\deg\pi=12$, so $\M$ is a line bundle inducing
a polarization of type $(1,6)$ on $A'$. Thus $\iota$ acts on
$H^0(\M)={\CC}^6$ in the usual way so that it has two eigenspaces, one
of dimension four and one of dimension two.  Now
$H^0(\O_{\Pfive}(3))^{\HHH'}$ is $\iota$-invariant, so that it must
map to one of these two eigenspaces. Hence the kernel of the above
restriction map is at least four dimensional. \Box

\remark{not.proj.normal}
Note that the Riemann-Roch theorem
tells us that we should only expect $H^0(\I_A(3))$ to
be two dimensional. Thus $A\subset\Pfive$ is not cubically normal. 

Eventually we will show that for a general abelian surface
$A\subseteq\Pfive$, $\dim H^0(\I_A(3))^{\HHH'}=4$. This will allow
us to define a rational map 
$$\Theta_6:
\A_{6}^{lev}\rDashto\Gr(2,W)$$ by taking 
$A$ to the two-dimensional subspace
$V\subseteq W$ such that $V_0\otimes V=H^0(\I_A(3))^{\HHH'}$. 
(See \ref{rat.map6} for the precise construction.)

\definition{rat.map} Define a rational map 
$$\phi:\Pfive\rDashto\Gr(2,W)$$ by taking a point $y\in\Pfive$ to a
subspace $V\subseteq W$ such that $V_0\otimes V$ is the largest
$\HHH_6$-subrepresentation of $H^0(\O_{\Pfive}(3))^{\HHH'}$ vanishing
at $y$.
\hfill\break\indent 
Equivalently, we first identify $\Gr(2,W)$ with
$\Gr(2,\dual{W})$ by identifying a subspace $V\subseteq W$ with its
annihilator $V^0\subseteq \dual{W}$, and then we can  define the map
$\phi:\Pfive \rDashto \Gr(2,\dual W)$ in coordinates, using the dual basis
to $e_0,\ldots,e_3$, as follows: a point $y \in \Pfive$ is mapped to
the subspace spanned by $(f_0(y),\ldots,f_3(y))$ and $((\sigma
f_0)(y),\ldots,(\sigma f_3)(y))\in \dual W$. 
\hfill\break\indent
We note also that $V_0\otimes V$ is
the subspace of $H^0(\O_{\Pfive}(3))^{\HHH'}$ vanishing on
the $\HHH_6$-orbit of $y\in\Pfive$. 

We will use Pl\"ucker coordinates for $\Gr(2,\dual W)$: the
coordinates for a $2$-plane spanned by $(x_0,\ldots,x_3)$ and
$(y_0,\ldots,y_3)$ will be $p_{ij}=x_iy_j-x_jy_i$.  These satisfy the
Pl\"ucker relation $p_{03}p_{12}-p_{02}p_{13}+p_{01}p_{23}=0,$ which
gives the equation of the Pl\"ucker embedding 
$\Gr(2,\dual W)\subseteq \Pfive$.  Thus, more
explicitly, $\phi$ is defined in Pl\"ucker coordinates by
$$y\in\Pfive\mapsto 
(f_0(y)(\sigma f_1)(y)-f_1(y)(\sigma f_0)(y),\ldots,
f_2(y)(\sigma f_3)(y)-f_3(y)(\sigma f_2)(y))\in\Gr(2,\dual{W})
\subset\Pfive.$$

\remark{detM3} Another way to obtain the $\HHH_6$-subrepresentation of
$H^0(\O_{\Pfive}(3))^{\HHH'}$ vanishing at a general point
$y\in\Pfive$ is to take the four cubics $\det(M_3(x,y))$,
$\det(M_3(\sigma^3(x),y))$, $\det(M_3(\tau^3(x),y))$, and
$\det(M_3(\sigma^3\tau^3(x),y))$; see \ref{prelim} and [GP1] for their
explicit form. As observed in [GP1], Remark 2.13, these four cubics do
not always span a four-dimensional space, e.g. when $y$ is contained
in an elliptic normal curve.

\lemma{ideal.of.cub} Let $A\subseteq\Pfive$ be a general $\HHH_6$-invariant
abelian surface. Then the map $\phi$ in \ref{rat.map} is defined 
at the general point of $A$, $\dim H^0(\I_A(3))^{\HHH'}=4,$ 
and the cubics in $H^0(\I_A(3))^{\HHH'}$ cut
out a scheme of dimension $\le 2$.

\proof: First, let $X(\Gamma_6)
\subseteq\Pfive$ be the hexagon of [GP1], \S 3. Then $\phi$
is defined on an open subset of $\Sec(X(\Gamma_6))=X(\partial C(6,4))$,
with notation as in [GP1].  Indeed, $\phi$ is
defined for instance at the point $(1:1:2:1:0:0)\in \Sec(X(\Gamma_6))$, 
as one checks by explicit evaluation. Since $X(\Gamma_6)\subset\Pfive$ 
is a degeneration of a general $\HHH_6$-invariant elliptic normal curve,
$\phi$ must be defined on a non-empty open subset of the secant variety 
of a general such curve. Thus $\phi$ is defined generically on a 
general translation scroll, and thus also on the the general 
$(1,6)$-polarized abelian surface by [GP1], Theorem 3.1.

Let $A\subseteq\Pfive$ be such an abelian surface. Let $y\in A$ be a
general point of the abelian surface, let $V\subseteq \dual{W}$ be the
subspace corresponding to $\phi(y)$, and let $V^0$ be its
annihilator. Then $V_0\otimes V^0$ is a four-dimensional space of
cubics whose zero-locus contains the orbit of $y$ under $\HHH_6$, and
no other cubic of $H^0(\O_{\Pfive}(3))^{\HHH'}$ vanishes along this
orbit.  Hence $V_0\otimes V^0\supseteq H^0(\I_A(3))^{\HHH'}$. On the
other hand $\dim H^0(\I_A(3))^{\HHH'}\ge 4$, by \ref{invar3}, so
equality holds.

Finally, suppose that $A\subset\Pfive$ is a general translation
scroll, so that $A\subseteq \Sec(E)$ for some non-singular elliptic
normal curve $E\subset\Pfive$. As before, $\dim
H^0(\I_A(3))^{\HHH'}=4$, and by [GP1], Theorem 5.2, two of these
cubics cut out $\Sec(E)$, which is irreducible. It follows that the
cubics of $H^0(\I_A(3))^{\HHH'}$ must cut out a surface, and the same
holds for the general $(1,6)$-polarized abelian surface, again by the
degeneration argument.\Box

\definition{rat.map6} Let $U\subseteq\Pfive$ be the largest open
set on which $\phi$ is defined.
We define a rational map 
$$\Theta_6:\A_{6}^{lev}\rDashto\Gr(2,\dual{W})$$ 
by sending the general Heisenberg invariant $(1,6)$-polarized
abelian surface $A\subseteq\Pfive$ to 
$\phi(A\cap U)$.\hfill\break\noindent 
We note that $\phi(A\cap U)$ is the single point in
$\Gr(2,\dual{W})$ whose annihilator corresponds to 
$H^0(\I_A(3))^{\HHH'}$.

Now in the diagram
\diagram[tight,height=2em,midshaft]
&&\X\subseteq \Pfive\times \Gr(2,\dual W)&&\\
&\ldTo^{\pi_1}&&\rdTo^{\pi_2}&\\
\Pfive&&&&\Gr(2,\dual W)\\
\enddiagram
let $\X$ be the closure of the graph of $\phi$, 
so that if $\pi_1$ and $\pi_2$ are the
projections to $\Pfive$ and $\Gr(2,\dual W)$, respectively, then
$\pi_1:\X\longrightarrow\Pfive$ is birational and $\pi_2=\phi\circ\pi_1$
on the open set of $\X$ where $\phi\circ\pi_1$ is defined.

We now define another subvariety $\C\subseteq \Pfive\times\Gr(2,\dual
W)$ as follows: If $V\subseteq\dual W$ is a two-dimensional subspace,
and $V^0\subseteq W$ is its annihilator, then $V_0\otimes V^0$ is a
subrepresentation of $V_0\otimes W$, which, by the above discussion,
may be identified with a four dimensional space of cubics.  Let
$$\C\subseteq \Pfive\times\Gr(2,\dual W)$$ be the universal family
defined by these cubics: i.e., if $V\in\Gr(2,\dual W)$, then $\C_V$ is
the scheme of zeros in $\Pfive$ of the ideal generated by the cubics
in $V_0\otimes V^0$.

Finally, let $Q\subseteq\Gr(2,\dual W)$ be the quadric defined by the
hyperplane section $\set{p_{03}+p_{12}=0}$ in the Pl\"ucker embedding.

Note first that $\X\subseteq\C$. Indeed, if $y\in\Pfive$ is a point
where $\phi$ is defined, and $\phi(y)=V\subset\dual W$, then the
cubics in $V_0\otimes V^0$ necessarily vanish at $y$.

We will now show the following

\theorem{descrip6} 
\item{\rm (1)} $\pi_2(\X)=Q$, the non-singular hyperplane section of
the Pl\"ucker embedding of
$\Gr(2,\dual W)$ given by $\set{p_{03}+p_{12}=0}$.
\item{\rm (2)} $\Theta_6:\A_{6}^{lev}\rDashto Q$ is a birational map, and 
$\pi_2 : \X \longrightarrow Q$ is birational to a twist of the
universal family over $\A^{lev}_{6}$. In particular, $\A^{lev}_{6}$
is a rational threefold.

\proof: (1) First note, by direct calculation, that
$$f_0(\sigma f_3)-f_3(\sigma f_0)+f_1(\sigma f_2)-f_2(\sigma f_1)=0,$$
which means exactly $\pi_2(\X)\subseteq Q$. To show that
$\pi_2(\X)=Q$, we will find a fibre of $\pi_2:\X\longrightarrow Q$ which
is non-empty and has dimension $\le 2$. Then since $\dim(\X)=5$ and
$\dim(Q)=3$, $\pi_2$ must be surjective, and have generic fibre
dimension $2$.

To show that ${\pi_2}_{\mid \X}$ has non-empty fibres of dimension
$2$, let $A\subseteq \Pfive$ be a general Heisenberg invariant
$(1,6)$-polarized abelian surface. As observed in \ref{rat.map6},
$\phi(A\cap U)$ consists of exactly one point, say, $V
\in\Gr(2,\dual{W})$. If $\tilde A\subseteq\X$ is the proper transform
of $A$ in $\X$, then $\pi_2(\tilde A)=V$. Furthermore, $\X_V\subseteq
\C_V$ and by \ref{ideal.of.cub}, $V_0\otimes V^0=H^0(\I_A(3))^{\HHH'}$
and $\C_V$ is a scheme of dimension $\le 2$. Hence $\X_V$ is non-empty
of dimension $2$.

(2) To prove the first part, it is enough to show that for general
$V\in Q$, the fibre $\C_V$ of $\pi_2:\C\longrightarrow Q$ is a scheme
whose two-dimensional component is of degree 12. Then there cannot be
two distinct abelian surfaces $A,A'\subseteq \C_{V}$, and so
$\Theta_6$ must be generically $1$ to $1$.

To check this claim about $\C_V$, we choose one specific $V\in Q$
spanned by $(1,0,0,0)$ and $(0,-1,1,0)$ in $\dual{W}$, so that
$V^0\subseteq W$ is spanned by $(0,1,1,0)$ and $(0,0,0,1)$.  Then
$\C_V$ is the subscheme of $\Pfive$ defined by the equations
$$\set{f_1+f_2=f_3=\sigma f_1+\sigma f_2=\sigma f_3=0},$$
or explicitly
$$\displaylines{ x_0x_2x_4=x_1x_3x_5=
x_1^2x_4+x_3^2x_0+x_5^2x_2+x_1x_2x_3+x_3x_4x_5+x_5x_0x_1=\cr
=x_1x_4^2+x_3x_0^2+x_5x_2^2+x_2x_3x_4+x_4x_5x_0+x_0x_1x_2=0.}$$ The
first two equations yield the degree $9$ (Stanley-Reisner) threefold
$X(\partial C(6,4))$ (see [GP1], Proposition 4.1) consisting of the
union of linear subspaces $L_{ij}$, $i\in\set{0,2,4},
j\in\set{1,3,5}$, with $L_{ij}=\set{x_i=x_j=0}$.  If $i\not\equiv
j+3\mod 6$, then $L_{ij}\cap\C_V$ is given by the equations
$$\set{x_i=x_j=x_k^2x_{k+3}+x_{k-2}x_{k-1}x_k+x_kx_{k+1}x_{k+2}=
x_kx_{k+3}^2+x_{k+1}x_{k+2}x_{k+3}+x_{k-3}x_{k-2}x_{k-1}=0},$$
where $2k-3=i+j \mod 6$, and $k\in\set{1,3,5}$.
This is a quadric surface 
$$\set{x_i=x_j=x_kx_{k+3}+x_{k-2}x_{k-1}=0},$$ or
$$\set{x_i=x_j=x_kx_{k+3}+x_kx_{k+1}=0},$$ 
together with a line, depending on the value of $i$ and $j$. 
If $i\equiv j+3\mod 6$,
then $L_{ij}\cap \C_V$ is given by the equations
$$x_i=x_j=\sum_{k\in\set{0,2,4}\atop k\not= i} x_k^2x_{k+3}=
\sum_{k\in\set{0,2,4}\atop k\not= i} x_kx_{k+3}^2=0,$$
which is easily seen to be a curve. Thus $\C_V$ is a union of 6
quadric surfaces and a number of (possibly embedded) curves. This
verifies the claim, and thus by the above discussion the general fibre
of $\pi_2:\X\longrightarrow Q$ is a $(1,6)$-polarized Heisenberg
invariant abelian surface. If $V=\Theta_6(A)$, then $\X_V=A$. This
shows the last part of (2) in the theorem.\Box

\remark{ninelines} (1) A more careful analysis of the specific set of
cubics studied in the proof of \ref{descrip6}, (2) shows that  in fact
the scheme $\C_V$ they define consists of the union of the six quadric 
surfaces $S,\sigma(S),\ldots,\sigma^5(S)$, where
$$S=\set{x_0=x_1=x_2x_5+x_3x_4=0}\subset\P^5,$$ and the $\HHH_6$-orbit
of the involution $\iota$ $(-1)$-eigenspace $\P^1_-$.  
This orbit consists of exactly nine lines. 
Conversely, if $A$ is a $\HHH_6^e$-invariant $(1,6)$ abelian
surface, $\P^1_-\cap A$ consists of four points, namely
the odd two-torsion points of $A$, and so $\P^1_-$ (and all its Heisenberg
translates) are contained in any cubic hypersurface containing
$A\subset\P^5$. This shows that for a general $(1,6)$-polarized
abelian surface $A\subset\Pfive$, 
the cubics containing it cut out the union of $A$
with the $\HHH_6$-orbit of $\Pone_-$.\hfill\break\indent 
(2) Furthermore, one may directly calculate the ideal of $\bigcup_{i=0}^5
\sigma^i(S)$, and one finds that the ideal is generated by the four
cubics $f_1+f_2,\sigma f_1+\sigma f_2, f_3,\sigma f_3$, and six
additional quartics. However, we have seen that for an arbitrary
$(1,6)$-polarized abelian surface $A$, $\dim H^0(\I_A(3))\ge 4$.  Thus we see
that for the general such $A$, $\dim H^0(\I_A(3))=4$, and that $\I_A$ is
generated by quartics.

\medskip
We next study the complete intersections of type $(3,3)$ arising in this
construction. They turn out to be  partial smoothings of the degenerate
``Calabi-Yau'' threefolds described in \ref{recallsec} and \ref{join}.

\definition{CI6} For a point $p\in\P(W)$ corresponding to a one dimensional
subspace $T\subseteq W$, let $V_{6,p}\subseteq\Pfive$ be the complete
intersection of type $(3,3)$ determined by the cubic hypersurfaces in
$V_0\otimes T\subset H^0(\O_{\P^5}(3))^{\HHH'}$. 
Notice that, by construction, $V_{6,p}\subseteq\Pfive$
is $\HHH_6$-invariant.

\theorem{CY6} For general $p\in\P(W)$ corresponding to $T\subseteq W$
we have 
\item{\rm (1)} $V_{6,p}$ contains a pencil of
$(1,6)$-polarized abelian surfaces, parametrized by a pencil of
two-dimensional subspaces contained in the annihilator of $T$.
\item{\rm (2)} $V_{6,p}$ is an irreducible threefold whose singular locus
consists of $72$ ordinary double points.  
These $72$ ordinary double points are the base locus of
the pencil in {\rm (1)}.  
\item{\rm (3)} There is a small resolution
$V^1_{6,p}\longrightarrow V_{6,p}$ of the ordinary double points, with
$V^1_{6,p}$ a Calabi-Yau threefold, and such that there is a map
$\pi_1:V^1_{6,p} \longrightarrow \Pone$ whose fibres form the pencil of
abelian surfaces of part {\rm (1)}.  
\item{\rm (4)} $\chi(V^1_{6,p})=0$ and
$h^{1,1}(V^1_{6,p})=h^{1,2}(V^1_{6,p})=6$.

\proof: (1) The choice of $V_{6,p}$ is given by $p\in\P(W)$, with
$p=\P(T)$, while the set of $(1,6)$-polarized abelian surfaces
contained in $V_{6,p}$ is parametrized by $V\in
Q\subseteq\Gr(2,\dual{W})$ such that $T\subseteq V^0$. Now as is
well-known, the lines in a non-singular hyperplane section of
$\Gr(2,\dual{W})$ correspond to pencils of lines in $\P(W)$, and each
point of $\P(W)$ is the center of exactly one of these pencils. Hence
for $p\in \P(W)$, with $p=\P(T)$, the set $\setdef{V\in Q}{V^0\supseteq T}$
is a line in $Q$. Each such $V$ gives rise to an abelian surface in
$V_{6,p}$, and hence we obtain a pencil of abelian surfaces in
$V_{6,p}$.

(2) That $V_{6,p}$ is irreducible for general $p\in\P(W)$ follows, for
example, from the fact that for a Heisenberg invariant
elliptic normal curve $E\subseteq\Pfive$, $\Sec(E)$ is of the form
$V_{6,p}$ for some $p$. This is part of [GP1], Remark 2.13 and Theorem 5.2.

To prove that in general $V_{6,p}$ has only $72$ ordinary double
points seems difficult without resorting to computational means. First
note that if $f_1$ and $f_2$ are two cubic hypersurfaces in $\P^5$
containing a $(1,6)$-polarized abelian surface $A$, and defining a
threefold $X=\set{f_1=f_2=0}$, then $f_1$ and $f_2$ induce sections of
$(\I_A/\I_A^2)(3)$, which are linearly dependent precisely on the set
$\Sing(X)\cap A$. Now $c_2((\I_A/\I_A^2)(3))=72$, so we would expect
to have exactly $72$ singularities of $X$ on $A$, counted with
multiplicities. On the other hand, $72$ ordinary double points can be
easily accounted for. Namely, in the pencil $V_0\otimes T\subset
H^0(\O_{\P^5}(3))^{\HHH'}$ defining $V_{6,p}$, there are exactly $4$
cubic hypersurfaces which are determinants of $3\times 3$-matrices
$M_3(x,y)$, for certain values of the parameter $y$. Such a cubic
hypersurface is singular along the locus where its $2\times 2$-minors
vanish, which is an elliptic normal curve in $\P^5$. In
particular, $V_{6,p}$ has $18$ ordinary double points at the points of
intersection of such an elliptic curve with (any) other cubic in the
pencil $V_0\otimes T$.  Using {\it Macaulay/Macaulay2\/}, one can easily
find examples of threefolds $V_{6,p}$ which do have precisely $72$
ordinary double points.  Thus $V_{6,p}$ for a general $p\in\P(W)$
possesses exactly $72$ ordinary double points.

Now, blowing up a smooth abelian surface $A\subseteq V_{6,P}$ in the pencil,
we obtain a small resolution $V^2_{6,p}\longrightarrow V_{6,p}$. Flopping
simultaneously the $72$ exceptional $\Pone$'s gives a small resolution
$V^1_{6,p}\longrightarrow V_{6,p}$ which is a Calabi-Yau threefold containing
a pencil of minimal abelian surfaces. Thus by \ref{ABpencil}, $V^1_{6,p}$ has
an abelian surface fibration. In particular, the base points of the 
pencil of abelian surfaces on $V_{6,p}$ consist precisely of the $72$ 
nodes. This completes the proof of (2) and (3).

For (4), the Euler characteristic computation follows immediately from
the fact that the Euler characteristic of a non-singular $(3,3)$
complete intersection in $\Pfive$ is $-144$. See \ref{defect6} below 
for the Hodge numbers.  \Box

\remark{defect6}
We sketch here a method of computing the Hodge numbers of $V^1_{6,p}$. 
An effective method of computing the Hodge numbers of nodal quintics in $\P^4$
is well-known (see [Scho], [We]); we essentially generalize this 
to nodal type $(3,3)$ and $(2,2,2,2)$ complete intersections in $\P^5$ and
$\P^7$, respectively.

Let $X$ be a complete intersection Calabi-Yau threefold in $\P^{n+3}$
with only ordinary double points as singularities, and assume that $X$ is
given by equations $f_1=\cdots=f_n=0$, with all equations being of the same
degree $d$. We assume also that there exists a projective small resolution
$\widetilde{X}\longrightarrow X$, for instance obtained by blowing up
a smooth Weil divisor passing through all the nodes, which is the case for
$V_{6,p}$ in \ref{CY6}. 

Let $S$ be the homogeneous coordinate ring of $X$, that is $S\cong
\CC[x_0,\ldots,x_{n+3}]/(f_1,\ldots,f_d)$. Let $T^1\cong
\Ext^1(\Omega^1_X,\O_X)$ be the tangent space to the deformation space
of $X$. It is easy to see that there is a natural isomorphism
$$T^1\cong \left( {S^n\over 
\setdef{(\partial f_1/\partial x_i,\ldots,\partial f_n/
\partial x_i)}{0\le i \le n+3}}\right)_d,$$
where the subscript $d$ refers to the degree $d$ homogeneous part.
Let $T^1_{loc}$ be the tangent space of the deformation space of the
germ of the singular locus of $X$. There is a natural map
$$T^1\longrightarrow T^1_{loc},$$ and we wish to compute the dimension
of the kernel of this map, which is the tangent space of the deformation
space of $\widetilde{X}$.  Thus, we need to know when
$(g_1,\ldots,g_n)\in S^n$ representing an element of $T^1$ yields a
trivial deformation on the ordinary double points. This is precisely
when the matrix
$$\pmatrix{g_1&\partial f_1/\partial x_0&&\partial f_1/\partial
x_{n+3}\cr \vdots&\vdots&\cdots&\vdots\cr g_n&\partial f_n/\partial
x_0&&\partial f_n/\partial x_{n+3}\cr}$$ 
has rank $<n$ at the singular
points of $X$. Thus the $n\times n$ minors of this matrix must be
contained in the ideal of the singular locus of $X$. This is a
calculation which can be easily performed using {\it Macaulay/Macaulay2\/}. 
One thus determines the dimension of the kernel of $T^1\longrightarrow
T^1_{loc}$, and this dimension is precisely $h^{1,2}(\widetilde{X})$. From the
knowledge of the topological Euler characteristic of $\widetilde{X}$, 
we can then compute the Picard number of $\widetilde{X}$.

We also note that even if $X\subseteq \P^m$ is not a complete intersection,
but the homogeneous ideal of $X$ is generated by elements $f_1,\ldots,f_n$
all of the same degree $d$, then $T^1$ can be computed as a subspace of
$$\left( {S^n\over 
\setdef{(\partial f_1/\partial x_i,\ldots,\partial f_n/
\partial x_i)}{0\le i \le m}}\right)_d.$$
See for example [Tei]. This makes it possible to compute Hodge
numbers even in this case. We omit the details, but will make use
of this on occasion to compute Hodge numbers of other Calabi-Yau threefolds.

\remark{birat16} Because the Picard number of $V^1_{6,p}$ is fairly
large, it would be difficult to understand all details of the geometry
of $V_{6,p}$! However, it is already interesting to identify several different
minimal models, by focusing only on the subgroup of $\Pic(V^1_{6,p})$
generated by $H$, the pull-back of a hyperplane in $\Pfive$, and $A$,
the class of a $(1,6)$-polarized abelian surface, 
that is a fibre of the abelian
surface fibration $\pi_1:V^1_{6,p}\longrightarrow\Pone$. There are
several curves of interest in $V^1_{6,p}$ also. Let $[e]$ be the class
of one of the $72$ exceptional lines. (Caution: these lines are not all in the
same homology class, but their homology classes cannot be
distinguished using intersection numbers with only $H$ and $A$.) By
\ref{ninelines} we have $\P^1_-\subset V_{6,p}$; let
$[l]$ be the class of $\Pone_-\subset V_{6,p}$. It is not difficult
to check that, for general $p$, $\Pone_-$ is 
disjoint from the singular locus of
$V_{6,p}$, and hence the proper transform of $\Pone_-$
is disjoint from all of the exceptional curves. 
With these preparations, in $V^1_{6,p}$,
we find the following intersection numbers:
$$\eqalign{H^3=9,\quad&H^2A=12,\quad A^2=0,\cr H\cdot e=0,\quad&H\cdot
l=1,\quad A\cdot e=1,\quad A\cdot l=4.\cr}$$ 
Here $A\cdot l=4$ since, by \ref{ninelines}, $\Pone_-$ intersects 
$A\subseteq\P^5$ in the four odd two-torsion points of $A$.

As in the proof of \ref{CY6}, by flopping simultaneously the $72$
exceptional curves we obtain the model $V^2_{6,p}$. There is still a
map $V^2_{6,p}\longrightarrow V_{6,p}$, which is obtained by blowing
up a general $A\subseteq V_{6,p}$. In $V^2_{6,p}$, we now have the
intersection numbers
$$\eqalign{H^3=9,\quad&H^2A=12,\quad HA^2=0,\quad A^3=-72,\cr
H\cdot e=0,\quad&H\cdot l=1,\quad A\cdot e=-1,\quad A\cdot l=4.\cr}$$
Now, since the ideal of a general abelian surface $A$ is generated by
cubics and quartics, the linear system $|4H-A|$
is base-point free on $V^2_{6,p}$. However, $(4H-A)\cdot l=0$, and hence
$|4H-A|$ induces a contraction $V^2_{6,p}\longrightarrow \overline{V^2_{6,p}}$,
contracting, by \ref{ninelines}, (2), only $l$ and its Heisenberg translates,
thus a total of $9$ lines. Furthermore, the normal bundle of $l$
can be easily computed to be $\O_{\Pone}(-1)\oplus\O_{\Pone}(-1)$.
These $9$ lines can be simultaneously flopped to obtain a model 
$V^3_{6,p}$. In this model,
$$\eqalign{H^3=0,\quad&H^2A=-24,\quad HA^2=-144,\quad A^3=-648,\cr
H\cdot e=0,\quad&H\cdot l=-1,\quad A\cdot e=-1,\quad A\cdot l=-4.\cr}$$ 
Now in $V^2_{6,p}$, the linear system $|3H-A|$ had the
union of $l$ and its Heisenberg translates as base locus, by
\ref{ninelines}.  Thus in $V^3_{6,p}$, $|3H-A|$ is base-point free,
and $\dim |3H-A|\ge 1$. Hence two distinct elements of $|3H-A|$ are in
fact disjoint, and thus $|3H-A|$ gives rise to a fibration
$\pi_2:V_{6,p}^3\longrightarrow\Pone$ whose (general) fibres are then
non-singular surfaces. A calculation shows that $(3H-A)\cdot
c_2(V^3_{6,p})=0$, and hence $\pi_2$ is an abelian surface
fibration. Let $A'$ be the proper transform in
$V_{6,p}\subseteq\Pfive$ of a non-singular fibre of $\pi_2$.  Then
$A'\subseteq\P^5$ is a non-minimal abelian surface, containing $9$
exceptional lines, precisely $\Pone_-$ and its Heisenberg
translates. In addition, we compute from the above tables that $\deg
A'=15$. Thus the linear system embedding $A'$ in $\Pfive$ is
$|L-\sum_{i=1}^9E_i|$, where $E_1,\ldots, E_9$ are the exceptional
curves and $L$ is a polarization on the minimal model $A_{min}'$
either of type $(1,12)$ or $(2,6)$.

To see that this polarization is in fact of type $(2,6)$, first note
by inspection that all cubics in ${H^0(\O_{\Pfive}(3))}^{\HHH'}$ are
in fact invariant under the Heisenberg involution $\iota$. Thus
$\iota$ acts on the abelian surface $A'$. The fixed points of this
action are exactly $A'\cap\Pthree_+$ and $A'\cap\Pone_-=\Pone_-$.
Intersection theory shows that $\# A'\cap\Pthree_+=15$. On the other
hand $\iota$ descends to the involution $x\mapsto -x$ on the minimal
model $A'_{min}$ of $A'$. In particular, if $H$ is a hyperplane in
$\Pfive$ containing $\Pthree_+$, then $H\cap A'$ descends to a
symmetric divisor $D$ on $A'_{min}$ vanishing with multiplicity $1$ at
all $16$ two-torsion points of $A'_{min}$. Therefore by [LB],
Proposition 4.7.5, $L$ must be of type $(2,6)$ rather than $(1,12)$.

\remark{} The general cubic hypersurface containing a $(1,6)$-polarized abelian
surface is smooth. An interesting question is whether it is
rational. See also [Has] for a detailed discussion of such rationality
issues.

\section {1.7} {Moduli of $(1,7)$-polarized abelian surfaces.}

The rationality of the moduli space $\A^{lev}_7$ of $(1,7)$-polarized abelian
surfaces was proved by Manolache and Schreyer [MS], who show that this
moduli space is birational to a special Fano threefold of type
$V_{22}$ of index $1$ and genus $12$. Their approach is based on a
detailed description of the minimal free resolution of the ideal sheaf
of a $(1,7)$-polarized abelian surface in $\P^6$.

We will give in the following a somewhat different proof of their
result, by using the fact that the general abelian surface $A$ with a
$(1,7)$ polarization is determined by the Pfaffian cubics containing
it. This will focus attention on certain Pfaffian Calabi-Yau threefold
containing a pencil of $(1,7)$-polarized abelian surfaces which were
first discovered by Aure and Ranestad (unpublished).  We will also
show in [GP3] the existence on such Pfaffian Calabi-Yau threefolds of
a second pencil of $(1,14)$-polarized abelian surfaces.

Let $N(\HHH_{7})$ denote the normalizer of the Heisenberg group
$\HHH_{7}$ inside $\SL(V)$, where as in \ref{prelim} the inclusion
$\HHH_{7}\hookrightarrow \SL(V)$ is via the Schr\"odinger
representation.  We have a sequence of inclusions
$$\langle I\rangle \subseteq Z(\HHH_{7})=\mu_7 \subseteq \HHH_{7}
\subseteq N(\HHH_{7}),$$ and as is well-known,
$N(\HHH_{7})/\HHH_{7}\cong \SL_2(\boldz_{7})$, and in fact
$N(\HHH_{7})$ is a semi-direct product $\HHH_{7}\rtimes
\SL_2(\boldz_{7})$ (see [HM], \S 1 for an identical discussion for the
group $\HHH_5$). Therefore the Schr\"odinger representation of
$\HHH_7$ induces a seven dimensional representation
$$\rho_7:\SL_2(\boldz_7)\longrightarrow \SL(V).$$
In terms of generators and relations (see for example [BeMe]), one has
$$\PSL_2(\boldz_7)=\langle S,T\mid S^7=1, (ST)^3=T^2=1, (S^2TS^4T)^3=1\rangle$$
where
$$\hbox{$S=\pmatrix{1&1\cr 0&1\cr}$ and $T=\pmatrix{0&-1\cr 1&0\cr}$,}$$
while the representation $\rho_7$ is given projectively by
$$\rho_7(S)={(\xi^{-{ij\over 2}}\delta_{ij})}_{i,j\in\ZZ_7}, \quad
\rho_7(T)={1\over\sqrt{-7}}{(\xi^{ij})}_{i,j\in\ZZ_7},$$
see [Ta] and [Si] for details. Here $\xi$ is a fixed primitive 
$7^\th$ root of unity.

The center of $\SL_2(\boldz_7)$ is generated by $T^2$, and
$\rho_7(T^2)=-\iota$.  Thus the representation $\rho_7$ is
reducible. In fact, if $V_+$ and $V_-$ denote the
positive and negative eigenspaces, respectively, of the Heisenberg
involution $\iota$ acting on $V$, then $V_+$ and $V_-$ are both
invariant under $\rho_7$, and $\rho_7$ splits as $\rho_+\oplus
\rho_-$, where $\rho_{\pm}$ is the representation of $\SL_2(\boldz_7)$
acting on $V_{\pm}$. Note that $\rho_-$ is trivial on the center of
$\SL_2(\boldz_7)$, so in fact it descends to give an irreducible
representation
$$\rho_-:\PSL_2(\boldz_{7})\longrightarrow \GL(V_-).$$

For the reader's convenience, we reproduce from the Atlas
of finite groups [CNPW] the character table for 
$\PSL_2(\boldz_{7})$:
$$
\vbox{\tabskip=0pt \offinterlineskip
\halign to 320pt
{\hfil#\hfil\quad\vrule&\strut\quad\hfil#\hfil\quad&\quad\hfil#\hfil\quad
&\quad\hfil#\hfil\quad&\quad\hfil#\hfil\quad&\quad\hfil#\hfil\quad
&\quad\hfil#\hfil\quad\cr
Size of&&&&&&\cr
Conjugacy class&1&21&42&56&24&24\cr
Representative&$I$&$S$&$T^2ST^{-2}ST^2$&$ST$&$T$&$T^{-1}$\cr
\noalign{\hrule}
Characters&&&&&&\cr
$\chi_1$&1&1&1&1&1&1\cr
$\chi_2$&3&-1&1&0&$\beta$&$\bar\beta$\cr
$\chi_3$&3&-1&1&0&$\bar\beta$&$\beta$\cr
$\chi_4$&6&2&0&0&-1&-1\cr
$\chi_5$&7&-1&-1&1&0&0\cr
$\chi_6$&8&0&0&-1&1&1\cr}}
$$ 
where $\beta={1\over 2}(-1+\sqrt{-7})$. We'll denote in the sequel
the corresponding representations by their characters.  The
representation $\rho_-$ is irreducible and has character $\chi_2$. The
polynomial invariants of this representation are classical and they
have first been determined by F. Klein, see [Kl1], [Kl2]. It turns out that
there are no invariants of degree $<4$, and the  quartic
$$f_4=x_{1}^{3}x_{2}-x_{2}^{3}x_{3}-x_{3}^{3}x_{1},$$
is the unique invariant in degree $4$. The smooth quartic curve
defined by this invariant
$\overline{X(7)}=\set{f_4=0}\subset\Ptwo_-$ is in fact
an isomorphic image of the modular curve of level $7$,
and has $\PSL_2(\boldz_{7})$ as its full automorphism group.
(See for instance [GP1], [Kl1], and [Ve] for details.)
The other primary invariants of this representation 
are a sextic $f_6$, which is the determinant of the Hessian matrix
of $f_4$, and $f_{14}$ a polynomial of degree fourteen which
is obtained as the determinant of a bordered Hessian of $f_4$.

We will also need a number basic facts concerning apolarity 
and polars of hypersurfaces, and we recall them briefly here.

Let $V$ be an $n$-dimensional vector space, and fix a basis for it. 
For a point $a$ with coordinates $a=(a_1,\ldots,a_n)\in V$, and a homogeneous
polynomial $F\in S^d(\dual{V})$ of degree $d$, one defines
$$P_a(F):={1\over d}\sum_{i=1}^n a_i{\partial F\over \partial x_i},$$
where $x_1,\ldots,x_n\in\dual{V}$ form the dual of the chosen basis. 
It is easy to see that the previous definition is independent 
of the choice of basis, and so if we further set
$$P_{a_1\cdot\cdots\cdot a_k}(F):=P_{a_1}(P_{a_2}(\cdots P_{a_k}(F))\cdots),$$
and then extend by linearity $P_{\Phi}(F)$ to all $\Phi\in S^k(V)$, we have
defined a pairing between $S(\dual{V})$ and $S({V})$, called
the {\it apolarity} pairing (see [DK] for a modern account and for details).
The resulting form  $P_{\Phi}(F)$ is called the {\it polar of $F$ with respect to $\Phi$}. 
Two homogeneous forms $\Phi$ and $F$ are called {\it apolar} if
$P_{\Phi}(F)=0=P_{F}(\Phi)$ (cf. [Sal] who says that the term was coined by Reye).  
One says that $\Phi\in S^{d-k}(V)$ is a {\it $k^\th$ antipolar} 
of a hyperplane $H$ with respect to $F$ if $P_{\Phi}(F)=H^k$. 

Finally, if $l_1,\ldots,l_s\in \dual{V}$ are linear forms such that 
$F=l_1^d+ \ldots +l_s^d$, then  $P_{\Phi}(F)=0$ for all
$\Phi\in I_\Gamma\subset S(V)$, where $I_\Gamma$ is the homogeneous
ideal of $\Gamma=\set{H_1,H_2,\ldots, H_s}\subset\P(\dual{V})$, the collection
of hyperplanes $H=\set{l_i=0}$. Conversely, if $P_{\Phi}(F)=0$ for all
$\Phi\in I_\Gamma$, with $\Gamma=\set{H_1,H_2,\ldots, H_s}$ a collection of points in 
the dual space, then $F$ is a sum of powers $F=l_1^d+ \ldots +l_s^d$, for suitably 
rescaled linear forms. One says that $l_1, \ldots,l_s$ 
(or more precisely that the corresponding points  $H_1, H_2,\ldots, H_s$ in 
$\P(\dual{V})$) form a {\it polar $s$-polyhedron} to $F$ 
if $F=\lambda_1 l_1^d+\ldots +\lambda_s l_s^d$, for
suitable scalars $\lambda_i$ (see  [DK] and [RS] for modern
accounts of apolarity).  

\medskip
We start by looking at the equations of an abelian surface in $\P^6$:

\proposition{ideal17} Let $A\subseteq\Psix$ be a $(1,7)$-polarized
abelian surface. Then $h^0(\I_A(2))=0$, $h^0(\I_A(3))=21$, and the
ideal of $A$ is generated by these $21$ cubics.

\proof: The first claim is proved by a direct argument in [MS], Lemma
2.3. Alternatively, [La] or [MS], Lemma 2.4, show that a
$(1,7)$-polarized abelian surface $A\subseteq\Psix$ is projectively
normal, and hence Riemann-Roch gives $h^0(\I_A(2))=0$ and
$h^0(\I_A(3))=21$.

Projective normality, together with Kodaira vanishing, imply also that
$\I_A$ is $4$-regular in the Castelnuovo-Mumford sense.  The
obstruction to being $3$-regular is that $h^3(\I_A)=h^2(\O_A)=1$.
However, $h^i(\I_A(3-i))=0$ for $i\ne 3$, and the comultiplication map
$$H^3(\I_A(-1))\longrightarrow H^0(\O_{\P^6}(1))^*\otimes H^3(\I_A),$$
is dual to the natural multiplication $H^0(\O_{\P^6}(1))\otimes H^0(\O_A)
\longrightarrow H^0(\O_A(1))$ and thus is an isomorphism. Therefore
we may apply [EPW], Lemma 8.8, to see that $\I_A$
is generated by $21$ cubics, and in fact to determine
all Betti numbers in the minimal resolution of this ideal sheaf. 
An alternative however partly incomplete argument may be found
in [MS], Corollary 2.4.\Box

Recall from [GP1], Corollary 2.8 that the $7\times 7$-matrix
$$M'_7(x,y)={(x_{{(i+j)}\over 2}y_{{(i-j)}\over 2})}_{i,j\in\boldz_7}$$ 
has rank at most $4$ on an embedded $\HHH_7$ invariant
$(1,7)$-polarized abelian surface in $\P^6$.
On the other hand, for any  parameter point 
$y=(0:y_1:y_2:y_3:-y_3:-y_2:-y_1)\in\Ptwo_-$, 
the matrix $M'_7$ is alternating. We will denote in the sequel by 
$I_3(y)\subset \CC[x_0,\ldots,x_6]$ the homogeneous ideal
generated by the $6\times 6$-Pfaffians of the alternating matrix $M_7'(x,y)$,
and by $V_{7,y}\subset\P^6$ the  closed subscheme defined by this ideal.

\proposition{CY7} 
\item{\rm (1)} For $y\in \overline{X(7)}=\set{f_4=0}\subset\Ptwo_-$, the scheme 
$V_{7,y}$ is the secant variety of an elliptic normal curve in $\P^6$
(the level $7$ structure elliptic curve corresponding to the point
$y$ on the modular curve $\overline{X(7)}$).
\item{\rm (2)} For a general $y\in\Ptwo_-$, the scheme 
$V_{7,y}$ is a projectively Gorenstein irreducible threefold of
degree $14$ and sectional genus $15$.

\proof: If $V_{7,y}$ is of the expected codimension three, 
then it will be of the degree and genus stated, and will be
projectively Gorenstein, as any such Pfaffian subscheme has these properties.
On the other hand, [GP1], Theorem 5.4, shows that for $y\in\P^2_-$ the origin of an 
Heisenberg invariant elliptic normal curve $E\subset\P^6$, 
we have $V_{7,y}=\Sec(E)$. Thus for general $y\in\P^2_-$, $V_{7,y}$ is irreducible
and of the expected codimension.\Box

\proposition{KleinQuartic}
\item{\rm (1)} For all $y\in \P^2_-\setminus \overline{X(7)}$, 
the scheme $V_{7,y}$ meets $\P^2_-$ along a conic $C_y$ and the point
$y$.
\item{\rm (2)} The conic $C_y$ is defined by the second polar $P_{y^2}(f_4)$
with respect to $y$ of the Klein modular quartic curve
$\overline{X(7)}=\set{f_4=0}\subset\P^2_-$. Furthermore, $C_y=C_{y'}$
if and only if $y=y'\in \Ptwo_-$, and $C_y$ is a singular conic if and
only if $y\in {\rm Hess}(\overline{X(7)})=
\set{f_6=0}\subset\P^2_-$. 
\item{\rm (3)} The point $y$ lies on $C_y$ if and only if $y\in \overline{X(7)}$. 
Moreover, if $y\in \overline{X(7)}$, then the conic $C_y$ touches
$\overline{X(7)}$ at the point $y$ (with multiplicity $2$, if $y$ is not
a flex of $\overline{X(7)}$).

\proof: Let $I'$ be the bihomogeneous ideal in
$\CC[x_1,x_2,x_3]\otimes \CC[y_1,y_2,y_3]$ generated by the $6\times
6$-Pfaffians of the alternating matrix $M'_7(x,y)$, where this time
both sets of coordinates $y=(0:y_1:y_2:y_3:-y_3:-y_2:-y_1)$ and
$x=(0:x_1:x_2:x_3:-x_3:-x_2:-x_1)$ are chosen in $\Ptwo_-$.

We wish to show that the ideal $I'$ takes the form $J\cdot
P_{y^2}(f_4)$, where $J$ is the bihomogeneous ideal of the diagonal
$\Delta$ in $\Ptwo_-\times\Ptwo_-$, namely
$$J:=(x_1y_2-x_2y_1, x_1y_3-x_3y_1, x_2y_3-x_3y_2),$$
and $P_{y^2}(f_4)$ is the second polar of the Klein quartic curve
$\overline{X(7)}\subseteq \Ptwo_-$ with respect to $y$, namely
$$P_{y^2}(f_4)={1\over 2}(
y_1y_2x_1^2-y_2y_3x_2^2-y_1y_3x_3^2+y_1^2x_1x_2-y_2^2x_2x_3-y_3^2x_1x_3).$$
While this can be shown by direct computation, it is easier to do so via
representation theory.

One may check by inspection that the bihomogeneous ideal $I'$ 
is equivariant with respect to the diagonal action of 
$\PSL_2(\boldz_7)$ on $\Ptwo_-\times\Ptwo_-$
given by the representation $\rho_-$ introduced at the beginning of
this section. In addition, the zero locus in  $\Ptwo_-\times\Ptwo_-$ of
$I'$ contains the diagonal $\Delta$, since the matrix $M_7'(y,y)$
has rank at most $4$ for every $y\in\Ptwo_-$. Thus $I'\subseteq J$,
and so the bihomogeneous $(3,3)$-part
${(I')}_{(3,3)}\subseteq J\cdot H^0(\O_{\Ptwo_-\times\Ptwo_-}(2,2))$.
Now
$$\eqalign{
H^0(\O_{\Ptwo_-\times\Ptwo_-}(1,1))&\cong\chi_3\oplus\chi_4\cr
H^0(\O_{\Ptwo_-\times\Ptwo_-}(2,2))&\cong\chi_1\oplus 2\chi_4\oplus \chi_5
\oplus 2\chi_6\cr
H^0(\O_{\Ptwo_-\times\Ptwo_-}(3,3))&\cong\chi_1\oplus 2\chi_2\oplus\chi_3
\oplus 5\chi_4\oplus 4\chi_5\oplus 4\chi_6\cr}$$
as can be easily calculated from the given character table of
$\PSL_2(\boldz_7)$.  In particular, as
$\PSL_2(\boldz_7)$-representations, $J_{(1,1)}\cong\chi_3$, and since
$\dim{(I')}_{(3,3)}\le 7$, because $I'$ is generated by $7$ Pfaffians
there is then no choice but for ${(I')}_{(3,3)}\cong\chi_3\otimes
\chi_1=\chi_3$. This means that ${(I')}_{(3,3)}=J_{(1,1)}\cdot \langle
f \rangle$, where $f$ is the unique $\PSL_2(\boldz_7)$-invariant of
bidegree $(2,2)$.  This invariant is in fact $P_{y^2}(f_4)$. Thus for
a fixed $y\in\Ptwo_-$, the scheme $V_{7,y}$ meets $\Ptwo_-$ in the
point $y$ and the conic $C_y=\set{P_{y^2}(f_4)=0}$. This proves part
(1).

To finish the proof of (2) and (3) note that $P_{y^2}(f_4)(y,y)=
y_{1}^{3}y_{2}-y_{2}^{3}y_{3}-y_{3}^{3}y_{1}$ which vanishes if and
only if $y\in \overline{X(7)}\subset\Ptwo_-$.  And in case $y\in
\overline{X(7)}\subset\Ptwo_-$, the conic $C_y$ is tangent to
$\overline{X(7)}$ at the point $y$ since
$$T_y(C_y) =T_y(\overline{X(7)})=
\setdef{ x\in \Ptwo_-}{(3y_1^2y_2-y_3^3)x_1+(y_1^3-3y_2^2y_3)x_2
-(3y_1y_3^2+y_2^3)x_3=0}.$$
Miele shows in his thesis [Mie] that the intersection multiplicity of
$C_y$ and $\overline{X(7)}$ at $y$ is two when $y\in\overline{X(7)}$
is not one of the flexes, and that in this case all of the other
intersection points of $\overline{X(7)}$ and $C_y$ are in fact simple.

Finally, the conic $C_y$ is singular if and only if
$$\rank{\left({\partial^2 f_4(y) \over \partial y_i\partial y_j}\right)}<3,$$
that is $y\in {\rm Hess}(\overline{X(7)})=\set{f_6=0}$ 
by the definition of the Hessian locus.\Box

\proposition{descript7}
Let $A\subset\P^6$ be a general Heisenberg invariant
$(1,7)$-polarized abelian surface, and let $A\cap\P^2_-=
\set{p_1, p_2,\ldots p_6}$ be the odd $2$-torsion points of 
$A$. Then: 
\item{\rm (1)} The points $p_i$ form a polar hexagon to the Klein
quartic curve $\overline{X(7)}$ in the dual space $(\Ptwo_-)^\ast$.
\item{\rm (2)} The surface $A$ is contained in $V_{7,p_i}$, 
for all $i=1,\ldots,6$. Moreover, $21$ cubic 
Pfaffians defining three of the six $V_{7,p_i}$'s 
generate the homogeneous ideal $I_A$ of $A$.

\proof: It follows from [GP1], Corollary 2.8, that for $y=p_i$ one
of the odd $2$-torsion points of $A$, the matrix $M'_7(x,y)$ is alternating
and of rank at most $4$ along the surface $A$, hence $A$ is 
contained in the Pfaffian scheme $V_{7,p_i}$, for all $i=1,\ldots,6$. 

None of the points $p_i\in A\cap\P^2_-$ lie on the quartic
$\overline{X(7)}\subset\P^2_-$ since otherwise, by \ref{KleinQuartic} (3), 
the corresponding scheme $V_{7,p_i}$ would be the secant variety 
of an elliptic normal curve in $\P^6$, which in turn does not 
contain smooth abelian surfaces. Thus, from \ref{KleinQuartic}, 
we deduce that $p_i\not\in C_{p_i}$, for all $i=1,\ldots,6$. 
In particular, for all $i\ne j$, the point $p_j$ lies necessarily 
on the conic $C_{p_i}$, and the point $p_i$ lies on the conic $C_{p_j}$. 
The four remaining odd $2$-torsion points of $A$, namely $\setdef{p_k}
{k\ne i,j}$, are then well determined as the four points of 
intersection of the two (smooth) conics $C_{p_i}$ and 
$C_{p_j}$. Now, [DK], Theorem 6.14.2, gives a criterion for 
when $6$ points in the plane form a polar hexagon for
a plane quartic curve, and the above description of 
$\set{p_1,\ldots,p_6}$ fits that criterion precisely. This proves
part (1).

We next show that $21$ cubics given by the submaximal Pfaffians
of  three of the skew-symmetric matrices $M'_7(x,p_i)$, $i=1,\ldots,6$, 
are linearly independent, and thus generate the ideal $I_A$. 
Indeed, obviously any $14$ such Pfaffians coming from two 
different matrices $M'_7(x,p_i)$ and $M'_7(x,p_j)$, 
are linearly independent. To show that  a third 
set of submaximal Pfaffians coming from $M'_7(x,p_k)$, for some $k\ne i,j$,
would be linearly dependent on the first two it is
enough to check this for a degenerate abelian surface $A\subset\P^7$.
This may be checked directly for the Stanley-Reisner degeneration
in [GP1], Proposition 4.4. By \ref{ideal17}, this proves (2).\Box

\remark{} One may check that for a general $(1,7)$-polarized abelian
surface $A$, the $21$ cubic Pfaffians defining any three of the six
associated $V_{7,p_i}$'s generate the homogeneous ideal $I_A$.

We have now obtained the same description of the moduli space
as in [MS], Theorem 4.9:

\corollary{moduli7} The moduli space of $(1,7)$-polarized abelian
surfaces with canonical level structure is birational to the space
$\VSP(\overline{X(7)},6)$ of polar hexagons to the Klein quartic curve
$\overline{X(7)}\subset (\P^2_-)^\ast$. This is a smooth special Fano
threefold of type $V_{22}$ of index $1$ and genus $12$. In particular,
$\A^{lev}_7$ is rational.

\proof: By \ref{descript7} (2), the polar hexagon
of odd two-torsion points of a general $(1,7)$-polarized
abelian surface $A$ with level structure uniquely determines 
$A$. Mukai [Muk1], [Muk2] has shown that the space of 
polar hexagons to a general plane quartic curve is a
a smooth Fano threefold of index $1$ and genus $12$. 
On the other hand, the analysis of [MS], Theorem 4.4  and (4.5), 
shows that the variety $\VSP(\overline{X(7)},6)$ of polar hexagons to 
the Klein quartic $\overline{X(7)}\subset (\P^2_-)^\ast$ is
general in Mukai's sense and thus also a smooth Fano threefold. 
In particular, being of the same dimension as the moduli space $\A_7^{lev}$, 
the general polar hexagon to $\overline{X(7)}\subset (\P^2_-)^\ast$
is the set of odd two-torsion points of some $(1,7)$-polarized
abelian surface with canonical level structure.

A discussion of the properties (smoothness, type of Fano threefold) of
$\VSP(\overline{X(7)}, 6)$  can be found in [MS], see
also [Schr] and [Muk2]. The rationality of such a Fano threefold
appears to have been known to Mukai [Muk1], and Ivskovskih [Isk].  
See also [MS], Theorem 4.10, for the sketch of a proof. \Box

\remark{fibration7} For a general $y\in \Ptwo_-$, the smooth
conic $C_{y}$ parametrizes a pencil of $(1,7)$-polarized abelian 
surfaces contained in $V_{7,y}$. The conic $C_y$, is in fact
also a conic in the anticanonical embedding of the special 
Fano threefold $\VSP(\overline{X(7)},6)$. See [Schr]
for a description of the conics on a Fano threefold
of index $1$ and genus $12$.

\proposition{hodge7} 
Let $y\in\Ptwo_-$ be a general point. Then
\item{\rm (1)} The threefold $V_{7,y}$ has as singularities
$49$ ordinary double points, which occur at the $\HHH_7$-orbit of the point $y$.
\item{\rm (2)} There is a small resolution
$V^1_{7,p}\longrightarrow V_{7,p}$ of the ordinary double points such that
$V^1_{7,p}$ is a Calabi-Yau threefold, and such that there is a map
$\pi_1:V^1_{7,p} \longrightarrow \Pone$ whose fibres form the pencil of
$(1,7)$-polarized abelian surfaces in \ref{fibration7}.
\item{\rm (3)} $\chi(V^1_{7,y})=0$ and $h^{1,1}(V^1_{7,y})=
h^{1,2}(V^1_{7,y})=2$.

\proof: For (1), we know of no better proof than that given
in Rodland's thesis [Rod]. One calculates the tangent cone
of $V_{7,y}$ at $y$ and finds in general an ordinary double point;
thus the $\HHH_7$-orbit of the point $y$ accounts for $49$ ordinary double
points. On the other hand, a {\it Macaulay/Macaulay2\/} calculation shows that
one can find $y\in\Ptwo_-$ for which $V_{7,y}$ has only
$49$ singular points. Thus for $y\in\Ptwo_-$ general, $V_{7,y}$
has precisely $49$ ordinary double points.

The small resolution in (2) is obtained, as in \ref{CY6}, by blowing up
a smooth $(1,7)$-polarized abelian surface contained in $V_{7,y}$
to obtain a small resolution $V^2_{7,y}\longrightarrow V_{7,y}$,
and then by flopping the $49$ resulting exceptional curves.

Part (3) follows from the fact that the general nonsingular Pfaffian
Calabi-Yau has Euler characteristic $\chi=-98$, while 
the calculation of the Hodge numbers $h^{1,1}$ and $h^{1,2}$ is 
done via the techniques of \ref{defect6}.\Box

\remark{pfaffkahler} We now discuss the K\"ahler cone of
various models of $V_{7,y}$. First note
that $H_2(V^1_{7,y},\boldz)$ contains two classes of interest: $e$,
the class of an exceptional curve in the small resolution, and $c$ 
the class of the conic $C_y$ contained in $V_{7,y}\cap\Ptwo_-$. 
It is then clear that in $V^1_{7,y}$, we have
$$\eqalign{H^3=14,\quad&H^2A=14,\quad A^2=0,\cr
H\cdot e=0,\quad&H\cdot c=2,\quad A\cdot e=1,\quad A\cdot c=5.\cr}$$
This in fact shows that $\Pic(V^1_{7,y})/{\rm Torsion}$ is generated by $H$ and $A$.
Indeed, if not, first note since $A\cdot e=1$, $A$ must be primitive
in $\Pic(V^1_{7,y})/{\rm Torsion}$, so $A$ and $aA+bH$ form a basis for
$\Pic(V^1_{7,y})/{\rm Torsion}$, where $a,b\in{\mathbf Q}$. But since
$(aA+bH)\cdot e=a$, $a\in\boldz$, and since 
$H^3=14$ is not divisible by a cube, we must have $b\in\boldz$.

We next consider $V^2_{7,y}$, the model obtained from
$V^1_{7,y}$ by flopping the above $49$ exceptional curves. 
One sees easily that
$$\eqalign{H^3=14,\quad&H^2A=14,\quad HA^2=0,\quad A^3=-49,\cr
H\cdot e=0,\quad&H\cdot c=2,\quad A\cdot e=-1,\quad A\cdot c=5.\cr}$$
We shall show in [GP3] that the K\"ahler cone of $V^2_{7,y}$
is spanned by $H$ and $5H-2A$, and that $|5H-2A|$ contracts the
$\HHH_7$-orbit of the conic $C_y$. Flopping these curves will
then yield a third model $V^3_{7,y}$ with intersection table
$$\eqalign{H^3=-378,\quad&H^2A=-966,\quad HA^2=-2450,\quad A^3=-6174,\cr
H\cdot e=0,\quad&H\cdot c=-2,\quad A\cdot e=-1,\quad A\cdot c=-5.\cr}$$
The K\"ahler cone of $V^3_{7,y}$ will be seen to be spanned
by $5H-2A$ and $7H-3A$. Finally, in this model $|7H-3A|$ is
a base-point free pencil of abelian surfaces with a polarization
of type $(1,14)$.

Finally we will see in [GP3] that the linear system $|5H-2A|$
maps the threefold $V^2_{7,y}$ into $\P^{13}$ 
as a codimension $7$ linear section
of the Pl\"ucker embedding of the Grasmmanian $\Gr(2,7)$.
This should potentially explain some of the numerical
similarities in the mirror symmetry computations in [Rod] and [Tjo].
  
\section {1.8} {Moduli of $(1,8)$-polarized abelian surfaces.}

We start by determining the quadratic equations of a $(1,8)$-polarized
abelian surface in $\P^7$.

Much as in the $(1,6)$ case, let $\HHH'$ be the subgroup of the
Heisenberg group $\HHH_8$ generated by $\sigma^4$ and $\tau^4$. We can
easily compute the space of $\HHH'$-invariant quadrics,
${H^0(\O_{\Pseven}(2))}^{\HHH'}$, and we see that it splits up into
three 4-dimensional isomorphic $\HHH_8$-representations
$${H^0(\O_{\Pseven}(2))}^{\HHH'}\cong \bigoplus_{i=0}^2
\langle f_i,\sigma f_i, \sigma^2 f_i, \sigma^3 f_i\rangle,$$
where 
$$f_0 = x_0^2+x_4^2,\qquad
f_1 = x_1x_7+x_3x_5,\quad {\rm and}\quad
f_2 = x_2x_6.$$

\remark{invar2} For a point $y\in\Ptwo_-$, we can, as in the
$(1,6)$ case, ask what is the largest 
$\HHH_8$-subrepresentation of ${H^0(\O_{\Pseven}(2))}^{\HHH'}$ 
vanishing at $y$, or equivalently, what is the subspace of 
${H^0(\O_{\Pseven}(2))}^{\HHH'}$ vanishing on the $\HHH_8$-orbit
of $y$. To determine this, for
$y=(0:y_1:y_2:y_3:0:-y_3:-y_2:-y_1)\in\Ptwo_-$, 
consider the $4\times 3$-matrix 
$$(\sigma^i(f_j(y)))_{{0\le i\le 3}\atop{0\le j\le 2}}=
\pmatrix{0&-y_1^2-y_3^2&-y_2^2\cr
y_1^2+y_3^2&0&-y_1y_3\cr
2y_2^2&2y_3y_1&0\cr
y_1^2+y_3^2&0&-y_1y_3\cr}.$$
Note that for any $y\in\Ptwo_-$, this matrix has rank $2$, and its
kernel is spanned by $(y_1y_3,-y_2^2,y_1^2+y_3^2)$. Thus the representation
of $\HHH_8$ spanned by 
$$f=y_1y_3f_0-y_2^2f_1+(y_1^2+y_3^2)f_2,\quad 
\sigma(f),\quad \sigma^2(f),\quad  {\rm and}\quad \sigma^3(f)$$ 
is the subspace of ${H^0(\O_{\Pseven}(2))}^{\HHH'}$ vanishing along
the $\HHH_8$ orbit of $y$. We note that this representation is also
spanned by the $4\times 4$ Pfaffians of the matrices $M_4(x,y)$,
$M_4(\sigma^4(x),y)$, $M_4(\tau^4(x),y)$ and
$M_4(\sigma^4\tau^4(x),y)$, though we will not need this fact in what
follows.

\lemma{dimH08} If $A\subseteq\Pseven$ is an abelian surface invariant
under the Schr\"odinger representation  of $\HHH_8$, then 
$\dim {H^0(\I_A(2))}^{\HHH'}=4$. 

\proof: This is similar to \ref{invar3}. First consider the
restriction map
$${H^0(\O_{\Pseven}(2))}^{\HHH'}\longrightarrow 
{H^0(\L^{\otimes 2})}^{\HHH'},$$ 
where $\L$ is the line bundle inducing the embedding
of $A$ in $\Pseven$.  Let $A'=A/4K(\L)$, and let $\pi:A\longrightarrow
A'$ be the quotient map. Then $4K(\L)$ acts on $\L^{\otimes 2}$, and
there exists a bundle $\M$ on $A'$ such that $\L^{\otimes 2}=\pi^*\M$,
and ${H^0(\L^{\otimes 2})}^{\HHH'}\cong H^0(\M)$. Now
$c_1(\M)^2=c_1(\L^{\otimes 2})^2/\deg(\pi) =16$, so $\dim
H^0(\M)=8$. Thus the dimension of the kernel of the restriction map is
at least $4$, so $\dim {H^0(\I_A(2))}^{\HHH'}\ge 4$.  Equality then
follows from \ref{invar2}, since $A\cap \Ptwo_-$ is non-empty.  \Box

Note that $\HHH'$ acts on $\Ptwo_-$, and the quotient $\Pzz$ is
easily seen to be isomorphic to $\Ptwo$.  From [GP1], \S 6, we have a morphism
$$\Theta_8: \A_{(1,8)}^{lev}\longrightarrow\Pzz\cong\Ptwo,$$ 
that essentially associates to an $(1,8)$-polarized abelian surface with
canonical level structure the class of its odd $2$-torsion
points. Note that this map is defined on all of $\A_{(1,8)}^{lev}$,
even for abelian surfaces where the $(1,8)$-polarization is not very
ample, since one still has a concept of where the odd $2$-torsion
points are mapped.

\theorem{rat.map8} 
The map $\Theta_8: \A_{(1,8)}^{lev}\longrightarrow \Pzz$ is
dominant.

\proof: Let $Z\subseteq\P^2_-$ be the inverse image of $\im(\Theta_8)$
under the projection
$\P^2_-\longrightarrow\P^2_-/\boldz_2\times\boldz_2$; we wish to show
that $\overline Z=\Ptwo_-$. Let $\A\subseteq \P^7_{\Delta}$ be a
Heisenberg invariant degeneration of a $(1,8)$-polarized abelian
surface with canonical level structure to a translation scroll
$S_{\tau,E}$ of an $\HHH_8$ invariant elliptic normal curve
$E\subset\P^7$, as given by [GP1],Theorem 3.1; see also the discussion
after \ref{recallsec}.  Then $S_{E,\tau}\cap\Ptwo_-\subseteq\overline
Z$, and so $\Sec(E)\cap\Ptwo_-\subseteq\overline Z$ for a general
Heisenberg invariant elliptic normal curve $E$ in $\Pseven$.

What is $\Sec(E)\cap \Ptwo_-$? By the proof of [GP1], Lemma 6.1 (b), a
secant line $\langle x,y \rangle$ of $E\subset\P^7$ intersects
$\Ptwo_-$ if and only if $x=-y$ on $E$. Thus $\Sec(E)\cap \Ptwo_-$ is
the projection of $E$ from $\Pfour_+$ to $\Ptwo_-$. Since $E\cap
\Pfour_+$ consists of the four two-torsion points of $E$ (see [LB],
Corollary 4.7.6), this projection is a two-to-one cover of a conic in
$\P^2_-$.

Thus to prove that $\overline Z=\P^2_-$, it is enough to show that as
$E$ varies, $\Sec(E)\cap \P^2_-$ sweeps out $\Ptwo_-$. To do this, it
is enough to find two (possibly degenerate) Heisenberg invariant
elliptic normal curves $E_1$ and $E_2$ such that $\Sec(E_1)\cap
\Ptwo_-$ and $\Sec(E_2)\cap \Ptwo_-$ do not coincide. We can take
$E_1$ to be the octagon $X(\Gamma_8)$, see [GP1], Theorem 3.2, and
$E_2$ to be any non-singular Heisenberg invariant elliptic normal
curve in $\P^7$. Then using the equations of [GP1], Proposition 4.1 and the
fact proved there that $\Sec(E_1) =X(\partial C(8,4))$, one finds
that, with coordinates $(x_1:x_2:x_3)$ on $\Ptwo_-$,
$\Sec(E_1)\cap\Ptwo_-$= $\set{x_1x_3=0}$, a reducible conic.  On the
other hand, $\Sec(E_2)\cap\Ptwo_-$ is an irreducible conic, so these
two curves cannot coincide.  This shows that $\Theta_8$ is dominant,
and concludes the proof. \Box

\definition{V8y} For a fixed point $y\in\Ptwo_-$, let $V_{8,y}$ denote the scheme in
$\Pseven$ defined by the quadrics in ${H^0(\O_{\Pseven}(2))}^{\HHH'}$
vanishing on the $\HHH_8$-orbit of $y$.

\theorem{singV8y} For general $y\in\Ptwo_-$, $V_{8,y}\subset\Pseven$ is a 
$(2,2,2,2)$-complete intersection, which is singular precisely at the
$64$ points which form the $\HHH_8$-orbit of $y$, and each of these
singular points is an ordinary double point.

\proof: This can be checked computationally in much the same way as in
\ref{CY6}, (2).  However, in this case the calculation can be carried
out by hand, and we do this here.

It is easy to see that for general $y$, $V_{8,y}\subset\P^7$ is a
threefold. Being a complete intersection, $V_{8,y}$ is singular at
$x\in V_{8,y}$ if and only if there is a quadric $Q\in
\P(H^0(\I_{V_{8,y}}(2)))$ which is singular at $x$. Thus, we need to
identify all singular quadrics in the web
$\P(H^0(\I_{V_{8,y}}(2)))=\Pthree$. Now $H^0(\I_{V_{8,y}}(2))$ is an
$\HHH_8$-representation, but since the subgroup $\HHH'$ acts
trivially, it is in fact an $\HHH_8/\HHH'\cong\HHH_4$
representation. Using coordinates $z=(z_0:z_1:z_2:z_3)$ on
$\Pthree=\P(H^0(\I_{V_{8,y}}(2)))$ so that $z$ corresponds to the
quadric $\sum_{i=0}^3 z_i\sigma^i(f)$, the action of $\HHH_4$ on
$\Pthree$ becomes the standard Schr\"odinger representation of
$\HHH_4$.  We will continue to write $\sigma$ and $\tau$ for the
generators of $\HHH_4$; this should not create any confusion.

Put
$$w=(w_0:w_1:w_2:w_3)=(2y_1y_3:-y_2^2:y_1^2+y_3^2:-y_2^2),$$
so, with notation as in \ref{invar2},
$$f={1\over 2}w_0f_0+w_1f_1+w_2f_2.$$
The point $w$ is considered fixed.

To find the singular quadrics in the web, we compute the Hessian of
the quadric $\sum_{i=0}^3 z_i\sigma^i(f)$ with respect to the
variables $x_0,\ldots,x_7$. It is convenient to order these variables
in two blocks as $x_0$, $x_2$, $x_4$, $x_6$, and $x_1$, $x_3$, $x_5$,
$x_7$, for then the Hessian is a block matrix
$$R(z,w)=\pmatrix{A&0\cr 0&B\cr}$$ 
with $A=(z_{i+j}w_{i-j})_{i,j\in\boldz_4}$ and
$B=\sigma(A)=(z_{i+j-1}w_{i-j})_{i,j\in\boldz_4}$, where $\sigma$ acts
only on the variables $z_0,\ldots,z_3$. Keeping in mind that
$w_3=w_1$, one can further write
$$A=\pmatrix{ A_1&A_2\cr A_2&A_1}$$
with
$$A_1=\pmatrix{z_0w_0&z_1w_1\cr z_1w_1&z_2w_0},\quad
A_2=\pmatrix{z_2w_2&z_3w_1\cr z_3w_1&z_0w_2\cr}.$$
Similarly,
$$B=\pmatrix{B_1&B_2\cr B_2&B_1\cr}$$ 
with $B_i=\sigma(A_i)$.
Then
$$\det R(z,w)=\det(A_1+A_2)\cdot\det(A_1-A_2)\cdot\det(B_1+B_2)\cdot
\det(B_1-B_2).$$
If $Q$ is the quadric in $\Pthree(z_0,z_1,z_2,z_3)$ given by the equation
$\det(A_1+A_2)=0$, then the locus $\set{\det R(z,w)=0}\subseteq
\Pthree(z_0,z_1,z_2,z_3)$ is given by
$$Q\cup\sigma(Q)\cup\tau(Q)\cup\sigma\tau(Q).$$
In other words, this is the locus of singular quadrics in the web
$|H^0(\I_{V_{8,y}}(2))|$. 

We note the following facts about the geometry 
of this discriminant locus, which hold for general $y\in\Ptwo_-$:

\lemma{rank-behavior} For general $y\in\Ptwo_-$, 
\item{\rm (1)} $Q$ is a quadric cone, singular at the point $p=(0:1:0:-1)\in\P^3$,
which is precisely where $A_1+A_2=0$.
\item{\rm (2)} $p\not\in \sigma(Q)\cup\tau(Q)\cup\sigma\tau(Q)$.
\item{\rm (3)} $\rank R(z,w)=8$ if $z\not\in Q\cup\sigma(Q)\cup\tau(Q)\cup
\sigma\tau(Q)$.
\item{\rm (4)} $\rank R(z,w)=7$ if $z$ is contained in precisely
one of $Q$, $\sigma(Q)$, $\tau(Q)$ and $\sigma\tau(Q)$, and
$z\not\in\set{p,\sigma(p),\tau(p),\sigma\tau(p)}$.
\item{\rm (5)} $\rank R(z,w)\ge 6$ if 
$z\in\set{p,\sigma(p),\tau(p),\sigma\tau(p)}$,
or $z$ is contained in precisely two of the quadrics $Q$, $\sigma(Q)$, 
$\tau(Q)$, $\sigma\tau(Q)$.
\item{\rm (6)} $\rank R(z,w)\ge 5$ if $z$ is contained in precisely three of
$Q,\sigma(Q),\tau(Q),\sigma\tau(Q)$.
\item{\rm (7)} $\rank R(z,w)\ge 5$ for all $z\in\Pthree$.
\item{\rm (8)} There are precisely $16$ points
$z$ in $\Pthree$ such that $\rank R(z,w)=5$, and they form is the
$\HHH_4$-orbit of the point $(w_3:w_2:w_1:w_0)$.

\proof:
(1) and (2) are easily checked. For (3)--(6), one notes that
$$\rank R(z,w)\ge\rank (A_1+A_2)+\rank (A_1-A_2)+\rank (B_1+B_2)
+\rank (B_1-B_2),$$
from which (3)--(6) follow, using (2) to note that if at least two
of these four matrices drop rank, none have rank $0$.
For (7), one needs to check that
$Q\cap\sigma(Q)\cap\tau(Q)\cap\sigma\tau(Q)=\emptyset$ for general
$y\in\Ptwo_-$. It is sufficient to check this at one $y$,
say $(y_1:y_2:y_3)=(0:1:0)$, in which case $w=(0:-1:0:-1)$ and
$Q=\set{(z_1+z_3)^2=0}$, in which case $Q\cap\sigma(Q)\cap
\tau(Q)\cap\sigma\tau(Q)=\emptyset$.

For (8), one can check by hand that
$$\eqalign{&(4w_1^4-2w_0^3w_2-2w_0w_2^3)(z_0+z_2)^2
+(w_0+w_2)^2\det(A_1+A_2)\cr
+&4w_1^2\sigma(\det(A_1+A_2))+(w_0-w_2)^2\tau(\det(A_1+A_2))=0.\cr}
\leqno{(6.1)}$$
Thus the net of quadrics spanned by $Q$, $\sigma(Q)$ and $\tau(Q)$
contains the doubled plane $(z_0+z_2)^2=0$, and hence $Q\cap
\sigma(Q)\cap\tau(Q)$ consists of at most $4$ distinct points, counted
doubly. The same then holds true by symmetry for any other
intersection of three of the four quadrics, and hence there are at most $16$ 
points contained in the intersection of three of the four quadrics.
On the other hand, the point $(w_3:w_2:w_1:w_0)$ is seen
to be in $Q\cap\tau(Q)\cap\sigma\tau(Q)$, and the $\HHH_4$-orbit
of $(w_3:w_2:w_1:w_0)$ consists of $16$ points
(for general $y$), all contained in the intersection of $3$
of the four quadrics. Thus this accounts for all such points.
One checks that $\rank R(z,w)=5$ for $z=(w_3:w_2:w_1:w_0)$.\Box

\noindent
{\sl Proof of \ref{singV8y} continued}:
Having understood the discriminant locus of the web, we now have to
understand the loci of vertices of singular quadrics in the web.  For
a given point $z\in\Pthree$, the corresponding quadric $\sum_{i=0}^3
z_i\sigma^i(f)$ has vertex $\P(\ker R(z,w))\subseteq \Pseven$.

Suppose $z\in Q$, and $z\not\in \sigma(Q)\cup\sigma\tau(Q)$. Then 
$\rank(B)=4$, and thus 
$$\P(\ker R(z,w))\subseteq \set{x_1=x_3=x_5=x_7=0}.$$
Thus 
$$\eqalign{V_{8,y}\cap\P(\ker R(z,w))\subseteq
\{x_1=x_3=x_5=x_7&={1\over 2} w_0(x_0^2+x_4^2)+w_2x_2x_6=\cr
&=w_1(x_0x_6+x_2x_4)=\cr
&={1\over 2}w_0(x_6^2+x_2^2)+w_2x_0x_4=\cr
&=w_1(x_6x_4+x_0x_2)=0\}\cr}$$
the latter of which is easily seen to be empty for general $y$.

If $z\in Q$, and $z\not\in\tau(Q)\cup\sigma\tau(Q)$, then  $\rank(A_1-A_2)=2$,
$\rank(B_1-B_2)=2$, and it follows that
$$\P(\ker R(z,w))\subseteq \set{x_0-x_4=x_1-x_5=x_2-x_6=x_3-x_7=0}.$$
Thus
$$\eqalign{V_{8,y}\cap\P(\ker R(z,w))\subseteq
\{x_0-x_4=x_1-x_5&=x_2-x_6=x_3-x_7=\cr
&=w_0x_0^2+2w_1x_3x_1+w_2x_2^2=\cr
&=w_0x_3^2+2w_1x_2x_0+w_2x_1^2=\cr
&=w_0x_2^2+2w_1x_1x_3+w_2x_0^2=\cr
&= w_0x_1^2+2w_1x_0x_2+w_2x_3^2=0\}\cr}$$
which again is seen to be empty for general $w$.

Finally, if $z\in Q$, and $z\not\in\sigma(Q)\cup\tau(Q)$, then $\rank(B_1+B_2)=2$,
$\rank(A_1-A_2)=2$, and 
$$\P(\ker R(z,w))\subseteq\set{x_0-x_4=x_1+x_5=x_2-x_6=x_3+x_7=0},$$
and
$$\eqalign{V_{8,y}\cap\P(\ker R(z,w))\subseteq
\{x_0-x_4=x_1+x_5&=x_2-x_6=x_3+x_7=\cr
&=w_0x_0^2-2w_1x_1x_3+w_2x_2^2=\cr
&=w_0x_3^2+2w_1x_0x_2-w_2x_1^2=\cr
&=w_0x_2^2+2w_1x_3x_1+w_2x_0^2=\cr
&=w_0x_1^2+2w_1x_2x_0-w_2x_3^2=0\}\cr}$$
which is again empty.

By the $\HHH_4$ symmetry, it follows that for general $y$ and
any point $z\in\Pthree$ such that $\rank R(z,w)\ge 6$, the singular
locus of the quadric $\set{\sum z_i\sigma^i(f)=0}$ is disjoint from $V_{8,y}$.
Thus the only contribution to the singularities of $V_{8,y}$ comes
from the quadric 
$$\sum_{i=0}^3 w_{3-i}\sigma^i(f)=0$$
and its Heisenberg translates.
Now $\sigma(\P(\ker R((w_3:w_2:w_1:w_0), w)))$ is easily seen to be $\Ptwo_-$,
and $V_{8,y}\cap\Ptwo_-=\set{y,\sigma^4(y),\tau^4(y),\sigma^4\tau^4(y)}$
for general $y$. Thus we see the singular locus of $V_{8,y}$, for
general $y\in\Ptwo_-$, is precisely the $\HHH_8$-orbit of $y$.

To figure out the nature of the singularities of $V_{8,y}$,
we only need now to observe that for general $y\in\Ptwo_-$,
\ref{dimH08} implies there is an abelian surface $A\subseteq
V_{8,y}$, and then $f,\sigma(f),\sigma^2(f),\sigma^4(f)$ yield four
sections of $(\I_A/\I_A^2)(2)$ which are linearly dependent precisely on
$\Sing(V_{8,y})\cap A$. Now a Chern class calculation as in \ref{CY6} shows that one
expects $64$ such points. However, any such singular point which is not
an ordinary double point counts with some non-trivial
multiplicity. Since we have identified precisely $64$ distinct
singular points, these must all be ordinary double points. \Box

\theorem{fibre8} The fibre of $\Theta_8$
over a general point $y\in\Pzz$ corresponds to a pencil of abelian
surfaces contained in the singular Calabi-Yau complete intersection
$V_{8,y}\subset\P^7$. In particular $\Theta_8$ gives
$\A_{(1,8)}^{lev}$ birationally the structure of a $\Pone$-bundle over
an open set of $\Ptwo_-/\boldz_2\times
\boldz_2\cong\Ptwo$.

\proof: By \ref{singV8y}, there is an open set $U\subseteq
\Ptwo_-/\boldz_2\times\boldz_2\cong\P^2$ such that $V_{8,y}$ is a singular
Calabi-Yau threefold with 64 ordinary double points for $y\in U$.
For $y\in U\cap\im(\Theta_8)$ (which is non-empty by \ref{rat.map8},)
there is a $(1,8)$-polarized Heisenberg invariant abelian surface
$A$ for which $y$ is the class of an odd two-torsion point.
In particular, $A\subseteq V_{8,y}$, and $A$ contains
$\Sing(V_{8,y})$, which is the $\HHH_8$-orbit of $y$. Thus
we obtain a projective small resolution $V^2_{8,y}\longrightarrow V_{8,y}$
which is a Calabi-Yau threefold,
by blowing-up $A$. After flopping the $64$ exceptional curves, we obtain
a small resolution $V^1_{8,y}\longrightarrow V_{8,y}$ which, by
\ref{ABpencil}, contains a base-point free pencil of abelian surfaces. On
the other hand the fibre of $\Theta_8$ is one-dimensional,
so we obtain in this way a one-parameter family of
abelian surfaces in $V^1_{8,y}$. If this one parameter family
of abelian surfaces is connected, then it must coincide with 
the pencil we have already constructed. However, if it were
not connected, then $V_{8,y}$ would contain at least two distinct
pencils of abelian surfaces. Let $A,A'\subseteq V_{8,y}$ be 
abelian surfaces in these two pencils. Then $A\cap A'$ necessarily
contains a curve. Indeed, otherwise their proper transforms in
$V^1_{8,y}$ meet only at points, hence have an empty intersection
as $V^1_{8,y}$ is non-singular. But this is only possible
if $A$ and $A'$ belong to the same pencil.

Let $C$ be the one-dimensional component of $A\cap A'$. If $A$, $A'$
are general, then their Neron-Severi groups are generated by $H$, so
$C$ is numerically equivalent to $nH$ for some $n$. Furthermore,
since $A$ and $A'$ are Heisenberg invariant, so is $C$.
Now we must have $n\le 3$. Indeed, by Riemann-Roch, $\dim
H^0(\I_A(3))/H^0(\I_{V_{8,y}}(3))\ge 16$. If $f\in H^0(\I_A(3))$
is a cubic not vanishing on $V_{8,y}$, then certainly for general choice
of $A'$ in the second pencil, $f$ does not vanish on $A'$, so
$A\cap A'$ is contained in a divisor of type $3H$.

Finally, to rule out this possibility, we note there are no Heisenberg
invariant curves on $A$ numerically equivalent to $nH$ for $n\le 3$;
this follows from [LB], Ex. (4), pg. 179.
Hence there is only one component of the fibre of $\Theta_8$,
and it is contained in a $\Pone$, which concludes the proof.
\Box

The above theorem allows us to conclude that $\A_{(1,8)}^{lev}$ is
uniruled, but it does not show it is rational or unirational. This is
because the $\Pone$-bundle may not be the projectivization of a rank 2
vector bundle. To determine rationality, one needs to know the open
set of $\Ptwo$ over which the $\Pone$-bundle structure is
defined. This can be done through a careful analysis of the
discriminant locus of the family of $(2,2,2,2)$-complete intersection
Calabi-Yau threefolds of type $V_{8,y}$. This is quite a tedious
exercise, so we will only sketch the results below.

\theorem{mod8} $\A_{(1,8)}^{lev}$ is birational to a conic bundle over
$\Ptwo$ with discriminant locus contained in the plane quartic
$$D=\set{2w_1^4-w_0^3w_2-w_0w_2^3=0}.$$
In particular $\A_{(1,8)}^{lev}$ is rational.

\proof: The basic idea is that the $\Pone$-bundle structure can only
break down over those points $y\in\Pzz$ for which $V_{8,y}$ is degenerate,
i.e. has worse than $64$ ordinary double points. However, some of these 
degenerations may also contain pencils of abelian surfaces, so a more careful
analysis of the discriminant locus is required. To determine the 
discriminant locus of the family of Calabi-Yau threefolds, one looks in detail
at the proof of \ref{singV8y} and determines precisely where each
step breaks down.

Examining facts (1)--(8) in \ref{rank-behavior} above, one sees that (1) breaks down on
the locus
$$L=\set{w_1(w_0+w_2)(w_0-w_2)=0},$$
which is a union of three lines.
Outside of these three lines, (2)--(6) continue to hold. However,
for (7) to hold, we need to know that
$$Q\cap\sigma(Q)\cap\tau(Q)\cap\sigma\tau(Q)=\emptyset.$$
Via a straightforward calculation, one can see that this occurs precisely on the 
curve 
$$D=\set{2w_1^4-w_0^3w_2-w_0w_2^3=0}.$$ 
In fact, by equation (6.1) in the proof of \ref{rank-behavior},
one even sees that on this curve the four quadrics $Q$, $\sigma(Q)$, 
$\tau(Q)$ and $\sigma\tau(Q)$ only span a pencil and hence intersect
along an elliptic curve. Now (8) holds off of $C$ and $L$.

Finally, in the last part of the proof of \ref{singV8y}, one finds that
$$V_{8,y}\cap\set{x_1=x_3=x_5=x_7=0}\not=\emptyset$$
if and only if 
$$w_0w_1w_2(w_0^2+w_2^2)=0,$$
while 
$$V_{8,y}\cap\set{x_0-x_4=x_1-x_5=x_2-x_6=x_3-x_7=0}\not=\emptyset$$
if and only if
$$(w_0-w_2)((w_0+w_2)^4-(2w_1)^4)=0,$$
while
$$V_{8,y}\cap\set{x_0-x_4=x_1+x_5=x_2-x_6=x_3+x_7=0}\not=\emptyset$$
if and only if
$$(w_0+w_2)((w_0-w_2)^4+(2w_1)^4)=0.$$
Putting this all together, one finds that $V_{8,y}$ might have
worse singularities than $64$ ordinary double points only over the locus
$$\Delta:=D\cup 
\set{w_0w_1w_2(w_0^2+w_2^2)(w_0^2-w_2^2)((w_0-w_2)^8-(2w_1)^8)=0},$$
which is a union of the smooth quartic curve $D$ and fifteen lines.

To aid the further analysis, we need to bring in the additional
$\SL_2(\boldz_8)$ symmetry present. Recall there is an exact sequence
(see for example [LB], Exercise 6.14)
$$\exact{\H(8)}{N(\H(8))}{\SL_2(\boldz_8)}$$
where $N(\H(8))$ is the normalizer of the Heisenberg group
$\H(8)\subseteq \GL(H^0(\O_A(1)))$
via the Schr\"odinger representation. Letting $\xi=e^{2\pi i/16}$
be a fixed primitive $16^\th$ root of unity, let $S$ and 
$T$ be the $8\times 8$ matrices
$$\eqalign{(S_{ij})&={(\xi^{-(i-j)^2})}_{0\le i,j\le 7}\cr
(T_{ij})&={(\xi^{-2ij})}_{0\le i,j\le 7}.\cr}$$ 
It is easy to check
that $S,T\in N(\H(8))$, and give via conjugation an action on
$\H(8)/{\CC}^*=\boldz_8\times\boldz_8$ defined by the matrices
$\smallmath\pmatrix{1&1\cr 0&1\cr}$ and 
$\smallmath\pmatrix{\hfill 0&1\cr-1&0\cr}$, respectively. 
In particular, $S$ and $T$ along with
$\H(8)$ generate the normalizer $N(\H(8))$, though we do not need
this fact. What is important for us is that for $\alpha\in
N(\H(8))$, we have $V_{8,\alpha(y)}=\alpha(V_{8,y})$.  Furthermore,
one sees easily that $S$ and $T$ act on the components of the
discriminant locus $\Delta$ and there are three orbits of this action:
the first being $D$, which is an $\SL_2(\boldz_8)$ invariant, the
second being $L=\set{w_1(w_0^2-w_2^2)=0}$, and the third being the
remaining $12$ lines in $\Delta$.  Thus to understand the generic
degeneration of $V_{8,y}$ along each component, it is enough to study
the components $D$, $\set{w_1=0}$, and $\set{w_0=0}$.

Looking at the equations for $V_{8,y}$ when $w_1=0$, it is immediate
that $V_{8,y}$ is the join $\Join(E_1,E_2)$ of two elliptic normal
quartic curves as in \ref{join}, and these threefolds still have a
natural pencil of abelian surfaces by \ref{newabeldegen}.  Over
$w_0=0$, one finds a Calabi-Yau with $72$ ordinary double points, and
these too still possess a pencil of abelian surfaces. Thus it is only
over the quartic $D$ that the $\Pone$-bundle structure may be lost.
It then follows from [Bea] that $\A_8^{lev}$ is rational.\Box

The following summarizes information about the threefold $V_{8,y}$:

\theorem{defect8} For general $y\in\Ptwo_-$, let $V^2_{8,y}\longrightarrow
V_{8,y}$ be the small resolution obtained by blowing up a smooth $(1,8)$-polarized
abelian surface $A\subseteq V_{8,y}$, and let $V^1_{8,y}$ be the small
resolution of $V_{8,y}$ obtained by flopping the $64$ exceptional curves
on $V^2_{8,y}$. Then
\item{\rm (1)} There exists an abelian surface fibration $\pi_1:
V^1_{8,y}\longrightarrow\Pone$;
\item{\rm (2)} $\chi(V^1_{8,y})=0$ and 
$h^{1,1}(V^1_{8,y})=h^{1,2}(V^1_{8,y})=2$.

\proof: (1) has already been proved in \ref{fibre8}.
(2) follows from calculations similar to those in \ref{defect6},
so in particular we may see $\P^2_-/\boldz_2\times\boldz_2$ 
as a compactification of the moduli space of $V_{8,y}$.\Box

\remark{kaehler8} 
The structure of the birational models of $V_{8,y}$ is quite
interesting and can be completely determined.  Let $H$ denote the pullback
of a hyperplane section of $V_{8,y}$ to $V^1_{8,y}$. Let $A$ be the
class of a $(1,8)$-polarized abelian surface in $V^1_{8,y}$.  Classes of curves
in $V^1_{8,y}$ include
\item{$\bullet$}$[l]$, the class of a line in $V_{8,y}$
disjoint from the singular locus and  contained in a translation scroll
fibre of the fibration $V^1_{8,y}\longrightarrow\Pone$, and 
\item{$\bullet$}$[e]$, the
class of an exceptional curve of the small resolution of $V_{8,y}$.

\noindent 
In general $V^1_{8,y}\longrightarrow\Pone$ has a translation
scroll fibre because, as in the proof of \ref{rat.map8}, there is a
translation scroll containing a general point $y\in\Ptwo_-$.  Then in
the model $V^1_{8,y}$,
$$\eqalign{H^3=16,\quad&H^2A=16,\quad A^2=0,\cr
H\cdot e=0,\quad&H\cdot l=1,\quad A\cdot e=1,\quad A\cdot l=0.\cr}$$
Thus in particular, $\Pic(V^1_{8,y})/{\rm Torsion}$ 
is generated by $H$ and $A$.
In the model $V^2_{8,y}$,
$$\eqalign{H^3=16,\quad&H^2A=16,\quad HA^2=0,\quad A^3=-64,\cr
H\cdot e=0,\quad&H\cdot l=1,\quad A\cdot e=-1,\quad A\cdot l=0.\cr}$$

The K\"ahler cone of $V^1_{8,y}$ is spanned by $H$ and $A$. We will
see shortly that in $V^2_{8,y}$ the linear system $|4H-2A|$ is a
base-point free pencil of abelian surfaces, by \ref{pencil2} below, so
the K\"ahler cone of $V^2_{8,y}$ is spanned by $H$ and
$2H-A$. Furthermore, it is then clear that $V^1_{8,y}$ and $V^2_{8,y}$
are the only minimal models of $V_{8,y}$. Computational evidence
leads us to speculate that the family of Calabi-Yau threefolds $V_{8,y}$
could possibly be self-mirror!

\bigskip
To determine the K\"ahler cone of the model $V^2_{8,y}$ in \ref{kaehler8}
we first need to determine the generators of the ideal of a 
$(1,8)$-polarized abelian surface.

\definition{minors8} If $z\in \Pseven$, let $W_{8,z}$ be the vanishing locus
of the three by three minors of the matrix $M_4(x,z)$; see \ref{prelim}
for notation.

Recall now the degenerations $X_8^{\lambda}$, $\lambda\in \CC^*$, of
$(1,8)$-polarized abelian surfaces with canonical level structure
defined in [GP1], \S 4.  For fixed $\lambda\in \CC^*$, the surface
$X_8^{\lambda}$ is the union of eight quadric surfaces
$$X_8^{\lambda}:=\bigcup_{i\in\boldz_8} Q_i^{\lambda},$$
with 
$$Q_i^{\lambda}:=
\{\hbox{$x_ix_{i+2}+\lambda x_{i-1}x_{i+3}=x_j=0,$ for $j\in 
\boldz_8\setminus \set{i,i+2,i-1,i+3}$}\}.$$

\theorem{idealx8} 
\item{\rm (1)} The ideal homogeneous $I(X_8^{\lambda})$ is generated by
the quadrics
$$\setdef{\sigma^i(f_1+\lambda f_2)}{0\le i\le 3}$$
and the cubics
$$\setdef{x_{i-2}x_ix_{i+2}}{i\in\boldz_8}\cup
\setdef{x_{i-1}x_ix_{i+1}}{i\in\boldz_8}.$$
\item{\rm (2)} For all $i\in\boldz_8$, $x_ix_{i+2}x_{i+5}\in I(X_8^{\lambda})$.
\item{\rm (3)} The family $\setdef{X_8^{\lambda}}{\lambda\in \CC^*}
\subseteq \CC^*\times\Pseven$ extends uniquely to a flat family
$$\setdef{X_8^{\lambda}}{\lambda\in\Pone\setminus\set{0}}
\subseteq (\Pone\setminus \set{0})\times\P^{n-1},$$ 
and  $X_8^{\infty}$ is the face variety
of the triangulation of the torus
\drawnoname{92}{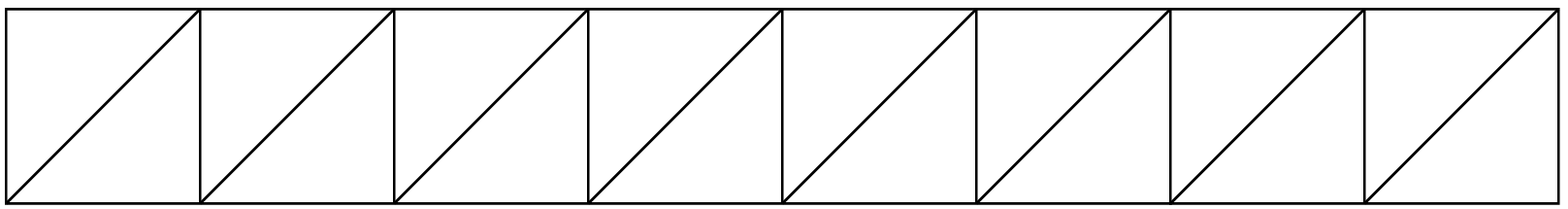}{
\swtext{.037}{1.000}{\swpad{2pt}{2pt}{$2$}}
\swtext{.151}{1.000}{\swpad{2pt}{2pt}{$3$}}
\swtext{.264}{1.000}{\swpad{2pt}{2pt}{$4$}}
\swtext{.378}{1.000}{\swpad{2pt}{2pt}{$5$}}
\swtext{.492}{1.000}{\swpad{2pt}{2pt}{$6$}}
\swtext{.606}{1.000}{\swpad{2pt}{2pt}{$7$}}
\swtext{.720}{1.000}{\swpad{2pt}{2pt}{$0$}}
\swtext{.833}{1.000}{\swpad{2pt}{2pt}{$1$}}
\swtext{.948}{1.000}{\swpad{2pt}{2pt}{$2$}}
\swtext{.037}{-0.169}{\swpad{2pt}{2pt}{$0$}}
\swtext{.151}{-0.169}{\swpad{2pt}{2pt}{$1$}}
\swtext{.264}{-0.169}{\swpad{2pt}{2pt}{$2$}}
\swtext{.378}{-0.169}{\swpad{2pt}{2pt}{$3$}}
\swtext{.492}{-0.169}{\swpad{2pt}{2pt}{$4$}}
\swtext{.606}{-0.169}{\swpad{2pt}{2pt}{$5$}}
\swtext{.720}{-0.169}{\swpad{2pt}{2pt}{$6$}}
\swtext{.833}{-0.169}{\swpad{2pt}{2pt}{$7$}}
\swtext{.948}{-0.169}{\swpad{2pt}{2pt}{$0$}}
}
In particular, $X_8^{\lambda}$ has the same Hilbert polynomial as a
smooth $(1,8)$-polarized abelian surface.

\proof: This uses very similar combinatorics to the proof of [GP1],
Theorem 4.6, where we neglected to mention in that theorem that part (a) held
only for $n\ge 10$. We omit the details.\Box

\theorem{ideal18} Let $A$ be a general $(1,8)$-polarized Heisenberg
invariant abelian surface in $\Pseven$. Then the embedding is
projectively normal and the homogeneous ideal of $A$ is generated by
the quadrics of ${H^0(\I_A(2))}^{\HHH'}$, and the $3\times 3$ minors
of a matrix $M_4(x,y)$, for $y\in A$ a general point.

\proof: This is a standard degeneration argument. Define $S_8
\subseteq \P^2_-\times\Pseven$, with coordinates $y_1,y_2,y_3$ on
the first factor and $z_0,\ldots,z_7$ on the second 
as the locus cut out by the equations
$$S_8=\setdef{\sigma^if(y_1,y_2,y_3,z_0,\ldots,z_7)=0}{i=0,\ldots,3}\subset
\P^2_-\times\Pseven,$$ 
where $f:=y_1y_3f_0-y_2^2f_1+(y_1^2+y_3^2)f_2$. Though one can show that
this scheme is irreducible by using the analysis of \ref{mod8}, we
will avoid this, so if necessary, replace $S_8$ with the irreducible
component of $S_8$ which dominates $\P^2_-$. Note that if $\dim
V_{8,y}=3$, then $\set{y}\times V_{8,y}\subseteq S_8$.  Next consider
the scheme $\SS\subseteq S_8\times \Pseven$ defined in
$S_8\times\Pseven$ by the equations
$\sigma^if(y_1,y_2,y_3,x_0,\ldots,x_7)=0$, $i=0,\ldots,3$, and the
$3\times 3$ minors of a matrix $M_4(x,z)$. A general fibre of the
projection $\SS\longrightarrow S_8$ over a point $(y,z)\in S_8$ is
obviously contained in the Calabi-Yau threefold $V_{8,y}$, and
moreover, by \ref{fibre8}, there is a unique $(1,8)$-polarized
Heisenberg invariant abelian surface $A\subseteq V_{8,y}$ containing
the point $z\in V_{8,y}$.  By [GP1], Corollary 2.7, the abelian
surface $A$ is contained in the fibre $\SS_{(y,z)}$. We show that
$A=\SS_{(y,z)}$ for the general $(y,z)$, by checking equality for one
special choice of $(y,z)$, namely $y=(0:1:\sqrt{-\lambda})$, with
$\lambda\not=0$, and $z=(1:1:0:0:0:0:0:0)\in V_{8,y}$.  For this
choice $\dim V_{8,y}=3$, and hence $(y,z)\in S_8$.  Moreover, for this
choice, the ideal of the fiber $\SS_{(y,z)}$ is generated by the
quadrics $\sigma^i(f_1+\lambda f_2)$, and the $3\times 3$ minors of
the matrix
$$\pmatrix{x_0&0&0&x_7\cr
x_1&x_2&0&0\cr
0&x_3&x_4&0\cr
0&0&x_5&x_6\cr}.$$
Therefore, by \ref{idealx8}, 
we deduce that $\SS_{(y,z)}=X_8^{\lambda}$, and so 
$\SS_{(y,z)}$ has the same Hilbert function as a $(1,8)$-polarized
abelian surface and is projectively normal. 
It follows that for  general $(y,z)$, we have $A=\SS_{(y,z)}$, 
which concludes the proof. \Box

\proposition{pencil2} For general $y\in\P^2_-$, $|4H-2A|$ is a
base-point free pencil on $V^2_{8,y}$, inducing a second abelian
surface fibration $\pi_2:V_{8,y}^2\longrightarrow\Pone$. A fibre $A'$
of $\pi_2$ is a $(2,8)$ or $(4,4)$-polarized abelian surface, mapped
to $\Pseven$ as a surface of degree $32$.

\proof: By [GP1], Corollary 2.7, $A\subset W_{8,z}$ for all $z\in
A\subset\Pseven$, so in particular the quartic hypersurface
$Q_z=\set{\det M_4(x,z)=0}$ vanishes doubly along the abelian surface
$A$, and hence defines an element of $|4H-2A|$ on $V^2_{8,y}$, for all
$z\in A$.  One checks for a special value of $y$, and therefore also
for the general $y$, that the divisors $Q_z$'s in $|4H-2A|$ span at
least a pencil. For example, one may take
$y=(0:1:\sqrt{-\lambda})\in\Ptwo_-$ for some $\lambda\in \CC^*$ and
$A=X_8^{\lambda}$, so that $V_{8,y}$ is given by the equations
$\set{\sigma^i(f_1+\lambda f_2)=0}$, and then take
$z_1=(1:0:\cdots:0)$, $z_2=(0:1:0:\cdots:0)$. Then
$Q_{z_1}=\set{x_0x_2x_4x_6=0}$, while $Q_{z_2}=\set{x_1x_3x_5x_7=0}$,
and one can check by hand that $V_{8,y}\cap Q_{z_1}\cap
Q_{z_2}=X_8^{\lambda}$, at least set-theoretically. Thus $Q_{z_1}$ and
$Q_{z_2}$ in particular define two independent elements of $|4H-2A|$,
and hence this linear system is at least a pencil. In addition, this
shows that $|4H-2A|$ has base locus contained in $A$. Since $A$ could
have been chosen to be any member of $|A|$, the base locus of
$|4H-2A|$ must be supported on the base locus of $|A|$. In the model
$V^2_{8,y}$, this base locus is the union of the $49$ exceptional
$\Pone$'s of numerical class $[e]$.  On the other hand, from the
intersection tables in \ref{kaehler8}, we compute on $V^2_{8,y}$ that
$(4H-2A)\cdot e=2$, so $4H-2A$ is nef.

Since $(4H-2A)^3= (4H-2A)^2H=(4H-2A)^2A=0$ it follows by [Og] that
$|4H-2A|$ is base-point free and induces a map to $\Pone$, with fibres
either isomorphic to a K3 or an abelian surface. Since moreover
$|2H-A|$ is empty as the quadrics containing $A$ are precisely those
containing $V_{8,y}$, we see that an element $A'\in |4H-2A|$ is
irreducible. It is in fact an abelian surface: by [Og], it is enough
to note that $(4H-2A)\cdot c_2(V_{8,y}^2)=0$, and this can be verified
easily.

Fix now a general point $z\in A\subset\Pseven$ and define
$\psi_z:\Pseven\rDashto \P^{15}$ to be the rational map induced by the
$3\times 3$-minors of the symmetric matrix $M_4(x,z)$. We regard here
$\psi_z$ as the restriction to $\Pseven$ of the Cremona involution
$\Phi:\P^{15}\rDashto\P^{15}$ which associates to each $4\times 4$
matrix its adjoint matrix. $\Phi$ is not defined on the locus $X_2$ of
rank $\le 2$ matrices (the secant variety to the Segre embedding of
$\Pthree\times\Pthree$ into $\P^{15}$), it contracts the locus $X_3$
of rank $\le 3$ matrices (a quartic hypersurface defined by the
determinant) to the locus $X_1$ of rank $\le 1$ matrices (the Segre
embedding of $\Pthree\times\Pthree$ into $\P^{15}$), and it is
one-to-one outside of the locus $X_3$.

In particular $\psi_z$ is birational onto its image and an isomorphism
outside the quartic hypersurface $Q_z$.  
Since $\O_{A'}(3H-A)=\O_{A'}(H+D)$, where
$D$ is a 2-torsion element, possibly zero, we also deduce that
$\psi_z$ maps any abelian surface $A'$ in $|4H-2A|$ via a linear system
induced by the restriction of $|H+D|$.  
Let $A^\prime_z$ be the surface in $|4H-2A|$
corresponding to the quartic $Q_z$. Then $\psi_z$ maps $A^\prime_z$
into the intersection of the Segre embedding of 
$\Pthree\times\Pthree\subset\P^{15}$
and $\psi_z(V_{8,y})$; in other words the restriction of $H+D$ to
$A^\prime_z$ decomposes as the sum 
$L_1+L_2$, where each $L_i$ induces the map on one of the two
$\Pthree$ components. Interchanging $L_1$ and $L_2$ amounts to
transposing $M_4(x,z)$, and since $M_4(x,z)^t=M_4(x,\iota(z))$, we
deduce that $L_1$ and $L_2$ are numerically of the same type, and
so on $A'$, we have
$H\equiv 2L$, where $L$ induces a polarization of type $(1,4)$ or
$(2,2)$. Hence $H$ is of type $(2,8)$ or $(4,4)$ on $A'$.  \Box

\remark{pol28}
A straightforward {\it Macaulay/Macaulay2\/} 
computation shows that $|L|$ induces a $2:1$ map of
$A^\prime_z$ onto a quartic surface in $\P^3$ which is singular along
two skew lines.  Thus, by [BLvS], $L$ must induce a polarization of
type $(1,4)$, and thus the abelian surfaces in the pencil $|4H-2A|$ are
in fact $(2,8)$-polarized.

\remark{RemarkCY8} By [GP1], Corollary 2.7, if $A\subseteq\Pseven$ 
is a Heisenberg invariant $(1,8)$-polarized abelian surface and $z\in A$, 
then $A\subseteq W_{8,z}$.
The expected codimension of the variety determined by
the $3\times 3$ minors of a $4\times 4$ matrix of linear forms is $4$.
On the other hand, by [GP1], Theorem 5.2, for special
values of $z$, $W_{8,z}$ is the secant scroll of an elliptic normal curve
in $\P^7$, so in general $W_{8,z}\subset\P^7$ is a threefold of degree $20$,
which is a partial smoothing of a degenerate ``Calabi-Yau'' threefold, by
by \ref{recallsec}. It turns out that $W_{8,z}$ is in fact a singular
model of Calabi-Yau threefold.
\smallskip

We mention without proof interesting facts about $W_{8,z}$:

$W_{8,z}\subset\P^7$ is invariant under the subgroup of $\HHH_8$
generated by $\sigma$ and $\tau^2$ and has, in general,
$32$ ordinary double points. The $(1,8)$-polarized abelian surface 
$A\subset W_{8,z}$ moves in a pencil whose base locus contains
the singular of the threefold. 

By blowing up $W_{8,y}$ along an abelian surface $A$ we obtain a small
resolution $W^2_{8,y}$ which is a Calabi-Yau threefold.  A Calabi-Yau
threefold in $\P^7$ defined by the $3\times 3$ minors of a general
$4\times 4$ matrix of linear forms, has Hodge numbers $h^{1,1}=2$,
$h^{1,2}=34$, while for $W^2_{8,y}$ one has $h^{1,1}=4$,
$h^{1,2}=4$. Let $H$ denote the pullback of a hyperplane section of
$W_{8,y}$ to $W^2_{8,y}$. Classes of curves in $W^2_{8,y}$ includes
$[c]$, the class of the conic in which $W_{8,y}$ meets $\Ptwo_-$ (see
the proof of \ref{rat.map8} for the fact that $W_{8,y}$ intersects
$\Ptwo_-$ in a conic), and $[e]$, the class of an exceptional $\P^1$
of the small resolution of $W_{8,y}$. Then on the model $W^2_{8,y}$,
$$\eqalign{H^3=20,\quad&H^2A=16,\quad HA^2=0,\quad A^3=-32,\cr H\cdot
e=0,\quad&H\cdot c=2,\quad A\cdot e=-1,\quad A\cdot c=4.\cr}$$ 
The linear system $|2H-A|$ is at least three dimensional, having a
subsystem defined by the quadrics containing $A$.  This gives a
morphism $\Psi:W^2_{8,y} \longrightarrow\Pthree$, whose image is in
fact a smooth quadric surface and whose fibres are elliptic
curves. $\Psi$ maps the exceptional lines of class $e$ onto the
rulings of the quadric $Q$, while the conics of class $c$ are
contracted. Each ruling of $Q$ induces on $W^2_{8,y}$ a K3-surface
fibration. A detailed analysis of the geometry of $W_{8,y}$ seems
difficult due to the fairly large rank of the Picard group, which
is generated over ${\bf Q}$
by $H$, $A$ and the fibres of the above two K3-fibrations.

\remark{CY24} It is interesting to compare the structures described
above with those of a $(2,4)$-polarized abelian surface
$A\subseteq\Pseven$, as studied in [Ba]. While numerically such a 
surface $A$ looks similar to a $(1,8)$-polarized abelian surface, 
in fact it is cut out by six quadrics. Any four of these 
quadrics defined a complete
intersection Calabi-Yau $X$ in $\Pseven$ containing $A$, again in
general with $64$ ordinary double points.  A Calabi-Yau small
resolution $\widetilde{X}\longrightarrow X$ exists, but now
$h^{1,1}(\widetilde X)=h^{1,2}(\widetilde X)=10$. We will say nothing
more about this threefold.

\section {1.10} {Moduli of $(1,10)$-polarized abelian surfaces.}

We consider first the representation theory of $\HHH_{10}$ on
the space of quadrics $H^0(\O_{\P^9}(2))$. 
There are four types of five dimensional irreducible
representations appearing in $H^0(\O_{\P^9}(2))$. Indeed, we can write
$H^0(\O_{\P^9}(2))\cong 3V_1\oplus 3V_2\oplus 3V_3 \oplus 2V_4$. The
three representations of type $V_1$ are given by the span of
$f,\sigma(f),\ldots,\sigma^4(f)$, where $f=x_0^2+x_5^2$,
$x_1x_9+x_6x_4$, or $x_2x_8+x_7x_3$. Three representations of type
$V_2$ are given similarly with $f=x_0x_5$, $x_9x_6 +x_4x_1$, or
$x_8x_7+x_3x_2$. Three representations of type $V_3$ are given by
$f=x_0^2-x_5^2$, $x_1x_9-x_6x_4$ or $x_2x_8-x_7x_3$, and finally the
two representations of type $V_4$ are given by $f=x_9x_6-x_4x_1$ or
$x_8x_7-x_3x_2$.

Note that for $y\in\P^3_-$, the matrix $M_5(x,y)$ is skew-symmetric,
while for $y\in\P^5_+$, $M_5(x,y)$ is symmetric. For a general
parameter point $y\in \P^3_-$, the $4\times 4$ Pfaffians of $M_5(x,y)$
cut out a variety $G_y\subset\P^9$ which can be identified with a
Pl\"ucker embedding of $\Gr(2,5)$ in $\P^9$. The corresponding
varieties defined by the Pfaffians of $M_{5}(\sigma^5(x),y)$,
$M_{5}(\sigma^5\tau^5(x),y)$ and $M_{5}(\tau^5(x),y)$ are then
$\sigma^5(G_y)$, $\sigma^5\tau^5(G_y)$ and $\tau^5(G_y)$,
respectively.

Note that the subgroup $\langle \sigma^5,\tau^5\rangle$ of $\HHH_{10}$ acts on
$\Pthree_-$, and we denote the quotient of this action by $\Ptzz$.  By
[GP1], Theorem 6.2, we have a map
$$\Theta_{10}:
\A_{(1,10)}^{lev}\longrightarrow \Ptzz$$
which essentially maps an abelian surface to the orbit of its
odd $2$-torsion points. The map $\Theta_{10}$ is birational onto its image. 
Since both spaces are three dimensional, we obtain the following result:

\theorem{rat.map10} $\A_{10}^{lev}$ is birationally equivalent to
$\P^3_-/\boldz_2\times\boldz_2$.\Box

\remark{quotPtzz} The quotient $\Ptzz$, being a quotient of $\Pthree$
by a finite abelian group, is rational by [Miya]. One can also compute
explicitly the ring of invariants of this action, and this ring is generated by
quadrics and quartics. There are $11$ independent invariants of degree
$4$, and these can be used to map $\P^3_-/\boldz_2\times\boldz_2$
into $\P^{10}$, where one finds the image to be a singular Fano
threefold of index $1$ and genus $9$.

We also have

\theorem{dimH010} Let $y\in\P^3_-$ be a general point. Then 
\item{\rm (1)} The twenty quadrics defining the four Pl\"ucker embedded
Grassmannians $G_y$, $\sigma^5(G_y)$, $\sigma^5\tau^5(G_y)$, $\tau^5(G_y)$ 
span only a $15$ dimensional space of quadrics $R_y\subset
H^0(\O_{\P^9}(2))$, which as a representation of
$\HHH_{10}$ is of type $V_1\oplus V_2 \oplus V_3$.  
\item{\rm (2)} The subspace $R_y$ generates the homogeneous ideal of the 
$(1,10)$-polarized abelian surface corresponding to the image via
$\Theta_{10}$ of $y$ in $\Ptzz$.

\proof: By [GP1], Corollary 2.7, if $A$ is a $\HHH_{10}$ invariant
abelian surface in $\P^9$ and  $y\in A\cap\P^3_-$ is an odd
$2$-torsion point, then $A\subseteq G_y\cap \sigma^5(G_y)
\cap \tau^5(G_y)\cap \sigma^5\tau^5(G_y)$, so all quadrics in the subspace
$R_y$ vanish on $A$. In fact, Theorem 6.2 of [GP1] states that $A=G_y\cap
\sigma^5(G_y)\cap\tau^5(G_y)$, and moreover the corresponding space of $15$
quadrics generates the homogeneous ideal of $A$ in $\P^9$. Thus we must have
$\dim R_y=15$. To see how $R_y$ decomposes as an $\HHH_{10}$ representation, 
it is sufficient to check the isomorphism 
$R_y\cong V_1\oplus V_2\oplus V_3$ for
a special value of $y$. One may use the point 
$y=(0:1:1:0:0:0:0:0:-1:-1)\in\P^3_-$, which in the notation
of [GP1], \S 4, yields the degenerate surface $X_{10}^1$, whose
homogeneous ideal, by [GP1] Theorem 3.6, is generated by the quadrics
$$\setdef{x_ix_{i+2}+x_{i-1}x_{i+3}}{i\in\boldz_{10}}\cup
\setdef{x_ix_{i+5}}{i\in\boldz_{10}}.$$
This clearly decomposes as $V_1\oplus V_2\oplus V_3$. \Box

\theorem{CY10} 
\item{\rm (1)} For general $y\in\Pthree_-$,
$V_{10,y}=G_y\cap \tau^5(G_y)\subset\P^9$ is a Calabi-Yau threefold
with $50$ ordinary double points, which is invariant under the
subgroup generated by $\sigma^2$ and $\tau$.  The singular locus of
$V_{10,y}$ is the $\langle\sigma^2,\tau\rangle$-orbit of
$\set{\sigma^5(y)}$. If $A\subseteq V_{10,y}$ is the $\HHH_{10}$-invariant
$(1,10)$-polarized
abelian surface which has $y$ as an odd $2$-torsion point, then the
linear system $|2H-A|$ induces a $2:1$ cover from $V_{10,y}$ to the
symmetric HM-quintic threefold $X_{5,y'}\subset\P^4$, where
$y'=(y_1y_2+y_3y_4:-y_2y_3:y_1y_4:y_1y_4:-y_2y_3)\in\Ptwo_+
\subseteq\Pfour$ (see \ref{symmoore}).  
\item{\rm (2)} If $y\in A\cap \Pfive_+$, let
$W_{10,y}\subset\P^9$ be the variety defined by the 
$3\times 3$ minors of the matrix
$M_5(x,y)$. For general $y$, this is a Calabi-Yau threefold of degree
$35$ with $25$ ordinary double points. Furthermore, $W_{10,y}$ is not simply
connected, but has an unbranched double cover birational to 
$X_{5,y''}\subset\P^4$,
with
$y''=(2(y_3y_4-y_1y_2):y_0y_1-y_4y_5:y_2y_5-y_0y_3:y_2y_5-y_0y_3:y_0y_1-y_4y_5)
\in\Ptwo_+\subseteq\Pfour$.

\proof: (1) To understand the singularity structure of $V_{10,y}$,
we will first understand the image of the map induced by the
linear system $|2H-A|$ in $\Pfour$, and show that this map expresses 
$V_{10,y}$ as a partial resolution of a double cover of $X_{5,y'}\subset\P^4$ 
branched over its singular locus.

First we note that for general $y\in\Pthree_-$, $V_{10,y}$ is an
irreducible threefold of degree $25$. To see this, we show it is true
for special choice of $y$. Choose $y\in\Pthree_-$ so that
$y_i=y_{5+i}$. Then $M_5(x,y)=(x'_{i+j}y_{i-j})$ where
$x_{i+j}'=x_{i+j}+x_{i+j+5}$. Thus $G_y$ is just a cone over an
elliptic normal curve $E$ in the $\Pfour$ given by
$\set{x_i-x_{i+5}=0}$, with vertex the $\Pfour$ given by
$\set{x_i+x_{i+5}=0}$.  $G_{\tau^5(y)}$ is a similar cone over
$\tau^5(E)$, and then $V_{10,y}=G_y\cap G_{\tau^5(y)}$ is just the linear
join of $E$ and $\tau^5(E)$, which is an irreducible threefold
of degree $25$, by \ref{join}. For general $y$, $G_y$ and $G_{\tau^5(y)}$
have dimension $6$ and degree $5$ in $\P^9$, so $G_y\cap G_{\tau^5(y)}$ is
expected to be of dimension $3$ and degree $25$. Since this holds for
special $y$, it holds for general $y$.

Next we show that $V_{10,y}$ contains a pencil of abelian surfaces.
Let $l$ be the line joining $y$ and $\tau^5(y)$, so that $\sigma^5(l)$
is the line joining $\sigma^5(y)$ and $\tau^5\sigma^5(y)$. Let
$z\in\sigma^5(l)$, such that $z\not\in\set{\sigma^5(y)+\sigma^5\tau^5(y),
\sigma^5(y)-\sigma^5\tau^5(y)}$. Then there is a linear transformation
$T\in \PGL(\P^9)$ of the form $T=\diag(a,b,\ldots,a,b)$ such that 
$T^2(z)=\sigma^5(y)$.
Note that $\tau^5 T=T\tau^5$, $\sigma^5T=T^{-1}\sigma^5$ in $\PGL(\P^9)$,
and that for any $w\in\Pthree_-$, we have $T(G_w)=G_{T^{-1}(w)}$.
(The significance of these commutation relations will
be explained in more detail in \S 1 of [GP3].) 
Denote by $A_{T(z)}\subset\P^9$ the abelian surface
determined by $T(z)\in\Pthree_-$, i.e.
$$A_{T(z)}=G_{T(z)}\cap G_{\sigma^5T(z)}\cap
G_{\tau^5T(z)}\cap G_{\sigma^5\tau^5T(z)}.$$
Then
$$\eqalign{T(A_{T(z)})&=G_z\cap G_{\sigma^5T^2(z)}\cap G_{\tau^5(z)}
\cap G_{\sigma^5\tau^5T^2(z)}\cr
&=G_z\cap G_{\tau^5(z)}\cap (G_y\cap G_{\tau^5(y)})\cr
&\subseteq V_{10,y}.\cr}$$
As $z$ varies on $\sigma^5(l)$, the point $T(z)$ also varies on this line,
and thus we obtain a pencil of distinct $(1,10)$-polarized 
abelian surfaces in $V_{10,y}$, filling up $V_{10,y}$ by irreducibility. 
Note that in particular, $\sigma^5(z)=T(\sigma^5(T(z)))
\in T(A_{T(z)})$ moves along the line $l$. Thus $l\subseteq V_{10,y}$,
and each abelian surface in the pencil intersects the line $l$.

Now we will study the linear system $|2H-A|$ on $V_{10,y}$. To do so,
we must first understand the equations of $V_{10,y}$ more deeply.
The Pfaffian of the matrix obtained by deleting the first row and
column of $M_5(x,y)$ is
$$\eqalign{p=&-x_0x_5(y_1y_2+y_3y_4)+(x_1x_4+x_6x_9)y_2y_3
-(x_2x_3+x_7x_8)y_1y_4\cr
&+(x_3x_7y_1^2+x_2x_8y_4^2)-(x_4x_6y_2^2+x_1x_9y_3^2)
+(x_5^2y_1y_3+x_0^2y_2y_4).\cr}$$
The ideal of $V_{10,y}$ is then generated by
$$\setdef{\sigma^{2i}(p)}{0\le i\le 4}\cup 
\setdef{\tau^5\sigma^{2i}(p)}{0\le i\le 4}.$$
On the other hand, the ideal of $A$ is generated by these ten
quadrics plus the additional quadrics 
$$\setdef{\sigma^{2i+5}(p)}{0\le i\le 4}.$$
Thus in particular, if we set
$$f_i:={1\over 2}\sigma^{2i}(\sigma^5(p+\tau^5(p))+p+\tau^5(p))$$
we see that
$$\eqalign{f_i=&(x_{3+2i}x_{7+2i}+x_{2+2i}x_{8+2i})(y_1^2+y_4^2)\cr
&+(x_{5+2i}^2+x_{2i}^2)(y_1y_3+y_2y_4)
-(x_{4+2i}x_{6+2i}+x_{1+2i}x_{9+2i})(y_2^2+y_3^2),}$$
and that $f_0,\ldots,f_4$ are linearly independent quadrics which
are still linearly independent when restricted to $V_{10,y}$. Furthermore,
each $f_i$ vanishes on $A$, and in fact $f_0,\ldots,f_4$ cut out
$A$ on $V_{10,y}$. Thus $f_0,\ldots,f_4$ define a four-dimensional
linear subsystem of $|H^0(\O_{V_{10,y}}(2))|$ whose
base locus is precisely the abelian surface $A$.

We now use $f_0,\ldots,f_4$ to define a rational map 
$f:V_{10,y}\rDashto \Pfour$. We describe next its image. 
The map $f$ is induced by the linear system $|2H-A|$ on 
$V^2_{10,y}$, where the model
$V^2_{10,y}\longrightarrow V_{10,y}$ is obtained by blowing up 
$V_{10,y}$ along $A$. We also denote by $f$ the induced 
morphism $V^2_{10,y}\rightarrow\Pfour$.

Let $z$ be a general point on the line $\sigma^5(l)$ joining
$\sigma^5(y)$ and $\sigma^5\tau^5(y)$, and $T$ as before. Let us consider
the restriction of $f$ to the surface $A'=T(A_{T(z)})$. 
There are two alternatives.
First if $A\cap A'$ is an isolated set of points, then the number of
these points is divisible by 50, since $A$ and $A'$ are invariant
under the action of $\tau$ and $\sigma^2$. Since $\sigma^5(y)
\in A\cap A'$, we have at least 50 such points, and the degree
of $f:A'\rDashto\Pfour$ is $4\cdot 20-a\cdot 50$, for some $a\ge 1$,
which is non-negative only if $a=1$. Thus $A\cap A'$ can only consist
of one $\langle\sigma^2,\tau\rangle$-orbit of $\sigma^5(y)$, and
$A$ must intersect $A'$ transversally. Thus the proper transform
of $A'$ in $V^2_{10,y}$, which we denote by $\tilde A'$, 
is $A'$ blown up in $50$ points.

The second alternative is that $A\cap A'$ contains a curve $C$. 
In this case, since $A$ is general, the curve $C$ can only 
be numerically equivalent to $nH$ on $A'$, for some $n$. 
Since $A$ is cut out by quadrics, $n\le 2$. On the
other hand, since $A$ and $A'$ are invariant with respect to $\sigma^2$ and
$\tau$, so is $C$. However, by [LB], Ex. (4), pg. 179, this is impossible.
Thus this case does not occur!

Next observe that for $x\in V_{10,y}$,
$f(x)=f(\tau^5(x))$, since each $f_i$ is $\tau^5$-invariant. Thus
$f$ factors via $V_{10,y}\longrightarrow V_{10,y}/\langle\tau^5
\rangle\rDashto\Pfour$.
Now $\tau^5$ acts on
$A'$. The $(1,10)$ polarization
$\L$ on $A'$ descends to a $(1,5)$ polarization $\L'$ on
$A'/\tau^5$, and $f_0,\ldots,f_4$ descend to sections of ${\L'}^{\otimes 2}$. 
The map $f:A'/\tau^5\rDashto\Pfour$ lifts 
then to a morphism $\tilde f:\tilde A'
/\tau^5\rightarrow\Pfour$.
It is easy to see that $\tilde f$ satisfies all the hypotheses
of \ref{nonmin15}. Thus $f$ maps
$T(A_{T(z)})$, for $z$ general, to a Heisenberg invariant
abelian surface $A''\subseteq\Pfour$ of degree $15$,
and $A''$ is contained in a quintic $X_{5,y'}\subset\P^4$ for any $y'\in
A''$. If we find a point $y'$ which is contained in $f(T(A_{T(z)}))$
for all $z\in\sigma^5(l)$, then $f(V_{10,y})\subseteq X_{5,y'}$.
A simple calculation shows that 
$$f(l)=(y_1y_2+y_3y_4:-y_2y_3:y_1y_4:y_1y_4:-y_2y_3)=:y'.$$
Since $T(A_{T(z)})$ intersects $l$ for each $z$, $y'$ is the desired
point, and so $f(V_{10,y})\subseteq X_{5,y'}$. To show 
that $f$ maps two-to-one onto
$X_{5,y'}$, we observe that on $V^2_{10,y}$, $H$ and $A$
are Cartier divisors, with $H^3=25$, $H^2A=20$, $HA^2=0$ and $A^3=-50$,
(the latter holds since
$(2H-A)^2\cdot A=30$, as $f$ maps $A$ two-to-one to a degree
15 surface in $\Pfour$). From this we see that $(2H-A)^3=10$. 
Thus $f(V_{10,y})=X_{5,y'}$ and $f$ maps two-to-one generically.

To compute the branch locus of the double cover $f:V^2_{10,y}\rightarrow
X_{5,y'}$, we need to identify
the fixed locus of the involution $\tau^5$ acting on $V_{10,y}$. This is
easily done: the fixed locus of $\tau^5$ acting on $\P^9$ consists of
$L_0=\set{x_0=x_2=x_4=x_6=x_8=0}$ and $L_1=\set{x_1=x_3=x_5=x_7=x_9=0}$.
Note we can write
$M_5(x,y)=(x_{6(i+j)}y_{6(i-j)}+x_{6(i+j)+5}y_{6(i-j)+5})_{0\le i,j\le 4}$,
making it clear that
$$M_5(x,y)|_{L_0}=-M_5(\tau^5(x),y)|_{L_0}
=(x_{6(i+j)+5}y_{6(i-j)+5})_{0\le i,j\le 4}$$
and
$$M_5(x,y)|_{L_1}=M_5(\tau^5(x),y)|_{L_1}
=(x_{6(i+j)}y_{6(i-j)})_{0\le i,j\le 4}.$$
Each of these are skew-symmetric Moore matrices on $\Pfour$. Thus
by [Hu1], [ADHPR2], Proposition 4.2, each of these matrices drops rank on 
an elliptic normal curve. Thus $(L_0\cup L_1)\cap V_{10,y}=E_0\cup E_1$, 
the disjoint union of two elliptic normal curves. 
This is the fixed locus of $\tau^5$ acting
on $V_{10,y}$. It is now clear that $f$ must map $E_0$ and $E_1$
to the two singular curves of $X_{5,y'}$ of \ref{symmoore}, (3). 

We can now understand from this the singularity structure of $V^2_{10,y}$.
First note that the exceptional locus of $\pi_2:V^2_{10,y}\rightarrow 
V_{10,y}$ consists of 50 $\Pone$'s. $A$ is a Cartier divisor
on $V^2_{10,y}$, so in particular $V^2_{10,y}$ is non-singular
along any non-singular element of the linear system $|A|$. 
But the exceptional locus of $\pi_2$ is contained in the base locus
of $|A|$, hence $V^2_{10,y}$ is non-singular along the exceptional
locus. In particular, since $\pi_2$ is
the small resolution of 50 ordinary double points and
$\omega_{V_{10,y}}$ is trivial, so is $\omega_{V^2_{10,y}}$.

Let $V^2_{10,y}\mapright{s_1} \hat X_{5,y'}\mapright{s_2}
X_{5,y'}$ be the Stein factorization of $f:V^2_{10,y}\rightarrow X_{5,y'}$.
Then $s_2$ is branched precisely over $\Sing(X_{5,y'})$, and is a double
cover. From the description of the singularities of $X_{5,y'}$ of
\ref{symmoore}, (5), one sees that $s_2(\Sing(\hat X_{5,y'}))$
is the Heisenberg orbit of $y'$, and shows
that $\hat X_{5,y'}$ has one ordinary double point sitting over each point of
the orbit of $y'$.

Next consider $s_1:V^2_{10,y}\rightarrow\hat X_{5,y'}$. We have already
seen that $f(l)=y'$; thus $s_1$ must contract $l$ and its
$\langle \sigma^2,\tau\rangle$-orbit (25 lines) to the 25 ordinary
double points of $\hat X_{5,y'}$. Since $\omega_{V^2_{10,y}}$ and 
$\omega_{\hat X_{5,y'}}$ are trivial, $s_1$ must be crepant, and the only
possibility then is that $s_1$ only contracts these $25$ lines, 
while $V^2_{10,y}$ is non-singular. Thus $V_{10,y}$ itself 
has $50$ ordinary double points
obtained by contracting the exceptional locus of $\pi_2$.

(2) Fix any point $y\in\P^5_+$. Consider first $\P(\Sym^2({\CC}^5))
=\Proj {\CC}[\setdef{x_{ij}}{0\le i,j\le 4, x_{ij}=x_{ji}}]$. The generic 
symmetric matrix $M=(x_{ij})$ has rank $2$ precisely on $\Sym^2(\Pfour)
\subseteq\P(\Sym^2({\CC}^5))$. One computes that $\Sym^2(\Pfour)$ has degree
35 in this embedding and thus $W_{10,y}$ can be described as $i(\P^9)\cap
\Sym^2(\Pfour)$, where $i:\P^9\longrightarrow\P(\Sym^2({\CC}^5))$ is the 
linear map given
by $x_{ij}=x_{i+j}y_{i-j}+x_{i+j+5}y_{i-j+5}$. Thus, if the intersection
is transversal, $W_{10,y}$ has dimension $3$ and degree $35$.

To understand $W_{10,y}$ further, consider it embedded in 
$\Sym^2(\Pfour)$ via the map $i$. Let $\pi:\Pfour\times\Pfour
\longrightarrow \Sym^2(\Pfour)$ be the quotient map, and let $z_0,\ldots,z_4,
w_0,\ldots,w_4$ be coordinates on the first and second $\Pfour$'s respectively.
The map $\pi$ is given by $x_{ij}=z_iw_j+z_jw_i$. Let us determine the equations
of $\pi^{-1}(W_{10,y})$, which is a double cover of $W_{10,y}$. 
One checks easily
that the equations of $i(\P^9)$ are
$$\setdef{(y_3y_4-y_1y_2)x_{i,i}+(y_0y_1-y_4y_5)x_{i-1,i+1}
+(y_2y_5-y_0y_3)x_{i-2,i+2}=0}{0\le i\le 4}.$$
Hence $\pi^{-1}(W_{10,y})$ is given by the five bilinear equations
$$\eqalign{
\{(y_3y_4-y_1y_2)(2z_iw_i)&+(y_0y_1-y_4y_5)(z_{i-1}w_{i+1}+z_{i+1}w_{i-1})\cr
&+(y_2y_5-y_0y_3)(z_{i-2}w_{i+2}+z_{i+2}w_{i-2})=0\}\cr}$$
If we write
$$y_0''=2(y_3y_4-y_1y_2),\quad y_1''=y_4''=y_0y_1-y_4y_5,
\quad y_2''=y_3''=y_2y_5-y_0y_3,$$
then we can rewrite the above set of equations as
$$L(z,y'')w=0,$$
where $L(z,y'')$ is the $5\times 5$ matrix 
$$L(z,y'')=(z_{2i-j}y''_{i-j}).$$
If $p_1,p_2:\Pfour\times\Pfour\longrightarrow\Pfour$ are the first and second
projections, then $X'_{5,y''}=p_1(\pi^{-1}(W_{10,y}))$ is given by 
the equation $\set{\det L(z,y'')=0}$. 
Thus it is clear that $W_{10,y}$ is of dimension three and hence of
degree $35$. Note that this argument shows that a generic intersection of
$\Sym^2(\Pfour)$ and $\P^9$ is a non-singular threefold (since the singular
locus of $\Sym^2(\Pfour)$ has codimension $4$) and has an unbranched covering
which is a Calabi-Yau threefold. We conclude that the generic section
of $\Sym^2(\Pfour)$ is a non-singular Calabi-Yau threefold, and hence also
$\omega_{W_{10,y}}=\O_{W_{10,y}}$.

Now note that
$${}^tL(z,y'')w=M'_5(w,y'')z$$
given that $y''\in\Ptwo_+$. Thus, if $\tilde X_{5,y''}\subseteq
\Pfour(z)\times\Pfour(w)$ is given by the equation
$${}^tL(z,y'')w=0,$$
then $p_1:\tilde X_{5,y''}\rightarrow X'_{5,y''}$ is a birational
map and $p_2(\tilde X_{5,y''})=X_{5,y''}$, where $X_{5,y''}$ is the
(symmetric) Horrocks-Mumford quintic given by $y''$. 
Then $p_2:\tilde X_{5,y''}
\rightarrow X_{5,y''}$ is also birational. Thus we see that
$\pi^{-1}(W_{10,y})$ is birational to $X_{5,y''}$. In addition, by
\ref{symmoore}, (4), $\tilde X_{5,y''}$ has $50$ ordinary double points.
Now $p_1:\pi^{-1}(W_{10,y})\rightarrow X'_{5,y''}$ is also
a crepant partial resolution, and fails to be an isomorphism
where $\rank L(z,y'')\le 3$. However,
$\rank L(z,y'')=\rank {}^tL(z,y'')$, so $p_1:\tilde X_{5,y''}
\rightarrow X'_{5,y''}$ and $p_1:\pi^{-1}(W_{10,y})\rightarrow X'_{5,y''}$
fail to be isomorphisms on the same set, which by \ref{symmoore}, (6),
is a union of two degree $5$ elliptic curves $E_1\cup E_2$. Furthermore,
$p_1^{-1}(E_i)\subseteq\tilde X_{5,y''}$ and $p_1^{-1}(E_i)\subseteq
\pi^{-1}(W_{10,y})$ are both $\Pone$-bundles. It follows from this
description that while $\tilde X_{5,y''}$ and $\pi^{-1}(W_{10,y})$
are not necessarily isomorphic over $X'_{5,y''}$, they must have the same
singularity structure, both being partial crepant resolutions of curves of
$cA_1$ singularities with $50$ $cA_2$ points. 
In particular, $\pi^{-1}(W_{10,y})$
has $50$ ordinary double points.

Finally, it is easy to see that for a general choice of $y$, 
$\pi^{-1}(W_{10,y})$
is disjoint from the diagonal in $\Pfour\times\Pfour$, so that
$\pi^{-1}(W_{10,y})\rightarrow W_{10,y}$ is an unbranched covering.
To see this, one can check for one point, e.g. a point of the
form $y=(y_0,\ldots,y_5,\ldots)=(0,y_1,y_2,y_3,y_4,0,\ldots)$, in which
case $y''=(2(y_3y_4-y_1y_2),0,\ldots,0)$, and it is particularly
easy to see that $\pi^{-1}(W_{10,y})$ is disjoint from the
diagonal.

Thus, for general $y$, $W_{10,y}$ must be a Calabi-Yau threefold
with $25$ ordinary double points.
\Box

\remark{k10} We discuss the structure of the K\"ahler cones of minimal models
of $V_{10,y}$.

First consider $V^1_{10,y}$, obtained by
flopping the 50 exceptional curves of $V^2_{10,y}\rightarrow V_{10,y}$.
Then $|A|$ is a base-point free pencil on $V^1_{10,y}$. 
Let $H$ be the pullback of a hyperplane section. Let $e$
be the class of an exceptional curve of the small resolution, and let
$l$ be the proper transform of the line $l$ joining $y$ and $\tau^5(y)$
as in the proof of \ref{CY10}. Then
$$\eqalign{H^3=25,\quad&H^2A=20,\quad A^2=0,\cr
H\cdot l=1,\quad&H\cdot e=0,\quad A\cdot l=2,\quad A\cdot e=1.\cr}$$
One can compute that $h^{1,1}(V^1_{10,y})=h^{1,2}(V^1_{10,y})=2,$ and
then it is clear that $\Pic(V^1_{10,y})/{\rm Torsion}$ is generated by
$H$ and $A$, and the K\"ahler cone of $V^1_{10,y}$ is spanned by $H$ and
$A$.

In $V^2_{10,y}$,
$$\eqalign{H^3=25,\quad&H^2A=20,\quad HA^2=0,\quad A^3=-50,\cr
H\cdot l=1,\quad&H\cdot e=0,\quad A\cdot l=2,\quad A\cdot e=-1.\cr}$$
The K\"ahler cone of $V^2_{10,y}$ is spanned by $H$ and $2H-A$; indeed,
it follows from the proof of \ref{CY10} that $|2H-A|$ on $V^2_{10,y}$
is base-point free and induces the map $V^2_{10,y}\longrightarrow X_{5,y'}$,
and the Stein factorization of this map contracts the
$\langle \sigma^2,\tau^2\rangle$ orbit of $l$. Thus these $25$ lines can
be flopped, to obtain a new model $V^3_{10,y}$. In this model,
$$\eqalign{H^3=0,\quad&H^2A=-30,\quad HA^2=-100,\quad A^3=-250\cr
H\cdot l=-1,\quad&H\cdot e=0,\quad A\cdot l=-2,\quad A\cdot e=-1.\cr}$$
Now $|2H-A|$ induces the morphism $f:V^3_{10,y}\longrightarrow X_{5,y'}$,
and $X_{5,y'}$ contains a pencil $|A'|$ of minimal $(1,5)$-polarized
abelian surfaces. Furthermore, if $X^1_{5,y'}\subseteq\Pfour\times |A'|$
defined by 
$$X^1_{5,y'}=\setdef{(x,S)}{x\in S\in|A'|}$$
is a partial resolution of $X^1_{5,y'}$ in which the pencil
$|A'|$ is base-point free, then $f$ factors through $\tilde f:V^3_{10,y}
\rightarrow X^1_{5,y'}$. It is then easy to see that if $S$
is a general element of this pencil, then $\tilde f^{-1}(S)$ is of class
$10H-6A$ on $V^3_{10,y}$. Also,
$\tilde f^{-1}(S)\longrightarrow S$ is a two-to-one unbranched cover.
Thus either $\tilde f^{-1}(S)$ is irreducible, or $\tilde f^{-1}(S)$
splits as a union of two $(1,5)$ abelian surfaces. However it is not
difficult to see that the divisor $5H-3A$ is effective. For example,
choosing a general point $z\in A$, $\det M_5(x,z)$ is a quintic
vanishing triply along $A$. Since a pencil of abelian surfaces on a Calabi-Yau
threefold cannot have a multiple fibre, we conclude that $|5H-3A|$
must be a pencil of abelian surfaces, and $\tilde f^{-1}(S)$ splits.
Hence $|5H-3A|$ yields a base-point free pencil of $(1,5)$-polarized
abelian surface, and the K\"ahler cone of $V^3_{10,y}$ is spanned by
$2H-A$ and $5H-3A$.

We will not address the K\"ahler cone structure of $W_{10,y}$
here. Because $h^{1,1}\ge 3$ for a small resolution of this
Calabi-Yau threefold, we in fact expect $W_{10,y}$ to inherit a rather
complicated K\"ahler cone structure from its double cover, a Horrocks-Mumford
quintic. See [Fry] for an analysis of the moving cone of the general
Horrocks-Mumford quintic.

\references
\item{[Ad]}
Adler, A., ``On the automorphism group of a certain cubic threefold'',
{\it Amer. J. Math.}  {\bf 100}, (1978), no. 6, 1275--1280.
\item{[Au]} Aure, A.,
 ``Surfaces on quintic threefolds associated to the Horrocks-Mumford bundle'',
in {\it Lecture Notes in Math.} {\bf 1399}, (1989), 1--9, Springer.
\item{[ADHPR1]} Aure, A.B., Decker, W., Hulek, K., Popescu, S., Ranestad, K.,
``The Geometry of Bielliptic Surfaces in $\Pfour$'',
{\it Internat. J. Math.} {\bf 4}, (1993) 873--902.
\item{[ADHPR2]} Aure, A.B., Decker, W., Hulek, K., Popescu, S., Ranestad, K.,
``Syzygies of abelian and bielliptic surfaces in ${\bf P}\sp 4$'',
{\it Internat. J. Math.} {\bf 8}, (1997), no. 7, 849--919.
\item{[Ba]} Barth, W., ``Abelian surfaces with $(1,2)$ polarization'',
{\it Adv. Studies in Pure Math.} {\bf 10}, Alg. Geom. Sendai, (1985) 41--84.
\item{[BM]} Barth, W., Moore, R., 
``Geometry in the space of Horrocks-Mumford surfaces'', 
{\it Topology}  {\bf 28}, (1989), no. 2, 231--245.
\item{[BHM]} Barth, W., Hulek, K., Moore, R., 
``Degenerations of Horrocks-Mumford surfaces'', 
{\it Math. Ann.} {\bf 277}, (1987), no. 4, 735--755. 
\item{[BS]} Bayer, D., Stillman, M.,
``Macaulay: A system for computation in
        algebraic geometry and commutative algebra
Source and object code available for Unix and Macintosh
        computers''. Contact the authors, or download from
        {\tt ftp://math.columbia.edu} via anonymous ftp.
\item{[Bea]} Beauville, A.,
``Vari\'et\'es de Prym et jacobiennes interm\'ediaires'',
{\it Ann. Sci. \'Ecole Norm. Sup.} (4) {\bf 10}, (1977), no. 3, 309--391.
\item{[BeMe]} Behr, H., Mennicke, J., ``A presentation of the groups 
${\rm PSL}(2,\,p)$'', {\it Canad. J. Math.}, {\bf 20}, (1968), 1432--1438.
\item{[BL]} Birkenhake, Ch., Lange, H., ``Moduli spaces of
abelian surfaces with isogeny'', {\it Geometry and analysis} (Bombay, 1992), 
225--243, {\it Tata Inst. Fund. Res.,} Bombay, 1995.
\item{[BLvS]} Birkenhake, Ch., Lange, H., van Straten, D., ``Abelian 
surfaces of type $(1,4)$'', {\it Math. Ann.} {\bf 285}, (1989) 625--646.
\item{[Bor1]} Borcea, C., ``On desingularized Horrocks-Mumford quintics'',
{\it J. Reine Angew. Math.} {\bf 421}, (1991), 23--41. 
\item{[Bor2]} Borcea, C.,
``Nodal quintic threefolds and nodal octic surfaces'', 
{\it Proc. Amer. Math. Soc.} {\bf 109}, (1990), no. 3, 627--635. 
\item{[Bori]} Borisov, L., ``A Finiteness Theorem for ${\rm Sp}(4,\ZZ)$'',
Algebraic geometry, 9. {\it J. Math. Sci.} (New York) {\bf 94} (1999), no. 1,
1073--1099, and {\tt math.AG/9510002}.
\item{[Cle]} Clemens, C. H., ``Double solids'', 
{\it Adv. in Math.} {\bf 47} (1983), no. 2, 107--230. 
\item{[CNPW]} Conway, J. H., Curtis, R. T., Norton, S. P., Parker, R. A., 
Wilson, R. A.,  {\it Atlas of finite groups. 
Maximal subgroups and ordinary characters for
simple groups. With computational assistance from J. G. Thackray}, 
Oxford University Press, Oxford, 1985.
\item{[DHS]} Debarre, O., Hulek, K., and Spandaw, J.,
``Very ample linear systems on abelian varieties,'' {\it Math. Ann.}
{\bf 300}, (1994) 181-202.
\item{[DK]}Dolgachev, I., Kanev, V., 
``Polar covariants of plane cubics and quartics'', 
{\it Adv. Math.} {\bf 98}, (1993), no. 2, 216--301. 
\item{[DO]} Dolgachev, I., Ortland, D., ``Points sets in
projective spaces and theta functions'', {\it Ast\'erisque}
{\bf 165}, (1988).
\item{[Do]} Donagi, R.,  ``The unirationality of ${\cal A}_5$'', 
{\it Ann. of Math.} (2) {\bf 119} (1984), no. 2, 269--307. 
\item{[EPW]} Eisenbud, D., Popescu, S., Walter, Ch.,
``Enriques Surfaces and other Non-Pfaffian Subcanonical Subschemes of
Codimension $3$'', {\it Communications in Algebra}, special volume for
Robin Hartshorne's $60^{\rm th}$ birthday, to appear,
{\tt math.AG/9906171}.
\item{[Fry]} Fryers, M., {\it The Moveable Cones of Calabi--Yau Threefolds,}
Thesis, University of Cambridge 1998. 
\item{[vG]} van Geemen, B., ``The moduli space of curves of genus $3$
with level $2$ structure is rational'', unpublished preprint.
\item{[GS]} Grayson, D., Stillman, M., 
``Macaulay 2: A computer program designed to support
computations in algebraic geometry and computer algebra.''
Source and object code available from
{\tt http://www.math.uiuc.edu/Macaulay2/}.
\item{[Gri1]} Gritsenko, V., ``Irrationality of the moduli spaces of 
polarized abelian surfaces. With an appendix by the author and K. Hulek''
in Abelian varieties (Egloffstein, 1993), 63--84, de Gruyter, Berlin, 1995.
\item{[Gri2]} Gritsenko, V., ``Irrationality of the moduli spaces of 
polarized abelian surfaces'', {\it Internat. Math. Res. Notices} (1994), 
no. {\bf 6}, 235 ff., approx. 9 pp. (electronic).
\item{[GP1]} Gross, M., Popescu, S., ``Equations of $(1,d)$-polarized
abelian surfaces,'' {\it Math. Ann.} {\bf 310}, (1998) 333--377.
\item{[GP2]} Gross, M., Popescu, S., ``The moduli space of $(1,11)$-polarized
abelian surfaces is unirational,'', {\it Compositio Math.} (2000), to appear,
preprint {\tt math.AG/9902017}.
\item{[GP3]} Gross, M., Popescu, S., ``Calabi-Yau threefolds and moduli of
abelian surfaces II'', in preparation.
\item{[Ha]} Hartshorne, R., {\it  Algebraic Geometry,} Springer 1977.
\item{[Has]} Hassett, B.,  
``Some rational cubic fourfolds'', {\it J. Algebraic Geom.}  
{\bf 8}, (1999), no. 1, 103--114. 
\item{[HM]} Horrocks, G., Mumford, D., `` A Rank 2 vector bundle on  $\Pfour$
with 15,000 symmetries'', {\it Topology}, {\bf 12}, (1973) 63--81.
\item{[Hu1]} Hulek, K., ``Projective geometry of elliptic curves'',
{\it Ast\'erisque} {\bf 137} (1986).
\item{[Hu2]} Hulek, K.,
``Geometry of the Horrocks-Mumford bundle'',
in ``Algebraic Geometry, Bowdoin 1985'',
{\it Proc. Symp. Pure Math.}, vol {\bf 46}, (2), (1987), 69--85.
\item{[HKW]} Hulek, K., Kahn, C., Weintraub, S., {\it Moduli Spaces of Abelian
Surfaces: Compactification, Degenerations, and Theta Functions},
Walter de Gruyter 1993.
\item{[HS1]} Hulek, K., Sankaran, G. K., 
``The Kodaira dimension of certain moduli spaces of abelian surfaces'', 
{\it Compositio Math.} {\bf 90}, (1994), no. 1, 1--35.
\item{[HS2]} Hulek, K., Sankaran, G. K., 
``The geometry of Siegel modular varieties'',
preprint {\tt math.AG/9810153}.
\item{[I]} Igusa, J., 
``Arithmetic variety of moduli for genus two'', 
{\it Ann. of Math.} (2) {\bf 72}, (1960), 612--649. 
\item{[Isk]} Iskovskikh, V. A., 
``Algebraic threefolds (a brief survey)'', (Russian) 
{\it Algebra}, 46--78, Moskov. Gos. Univ., Moscow, 1982.
\item{[Iy]} Iyer, J. N.,
Projective normality of abelian surfaces given by primitive line bundles,
{\it Manuscripta Math.} {\bf 98}, (1999), no. 2, 139--153.
\item{[Kat]} Katsylo, P., ``Rationality of the moduli 
variety of curves of genus $3$'', 
{\it Comment. Math. Helv.} {\bf 71} (1996), no. 4, 507--524. 
\item{[Kl1]} Klein, F., ``\"Uber transformationen siebenter Ordnung 
der elliptischen Funktionen'' (1878/79), Abhandlung {\bf LXXXIV}, in
{\it Gesammelte Werke}, Bd. {\bf III}, Springer, Berlin 1924.
\item{[Kl2]} Klein, F., ``\"Uber die ellip\-ti\-schen Nor\-mal\-kur\-ven der
$n$-ten Ord\-nung'' (1885), Abhandlung {\bf XC}, in
{\it Gesammelte Werke}, Bd. {\bf III}, Springer, Berlin 1924.
\item{[La]} Lazarsfeld, R., ``Projectivit\'e normale  
des surfaces abeliennes'', (written by O. Debarre), 
preprint Europroj {\bf 14}, Nice.
\item{[LB]} Lange, H., Birkenhake, Ch., {\it Complex abelian varieties},
Springer-Verlag 1992.
\item{[MS]} Manolache, N., Schreyer, F.-O.,
``Moduli of $(1,7)$-polarized abelian surfaces via syzygies'',
preprint {\tt math.AG/9812121}.
\item{[Mie]} Miele, M., ``Klassifikation der 
Durchscnitte Heisenberg-invarianter
Systeme von Quadriken in $\P^6$'', Thesis, Erlangen 1995.
\item{[Mi]} Miyata, T., ``Invariants of Certain Groups I,'' {\it Nagoya
Math. J.,} {\bf 41} (1971), 69--73.
\item{[Moo]} Moore, R.,
``Heisenberg-invariant quintic 3-folds, and sections
of the Horrocks-Mumford bundle'', preprint, Canberra 1985.
\item{[MM]} Mori, S., Mukai, S., ``The uniruledness of the moduli 
space of curves of genus $11$'', in  {\it Algebraic geometry} 
(Tokyo/Kyoto, 1982), 334--353,  Lecture Notes in Math. {\bf 1016}, 
Springer, 1983. 
\item{[Muk1]} Mukai, S., ``Biregular classification of Fano 
$3$-folds and Fano manifolds of coindex $3$'', 
{\it Proc. Nat. Acad. Sci. U.S.A.}  {\bf 86}, (1989), no. 9, 3000--3002. 
\item{[Muk2]} Mukai, S., ``Fano $3$-folds'', in
{\it Complex projective geometry} (Trieste, 1989/Bergen, 1989), 
255--263, London Math. Soc. Lecture Note Ser. {\bf 179}, 
Cambridge Univ. Press, Cambridge, 1992. 
\item{[Mu1]} Mumford, D., ``On the equations defining abelian varieties'',
{\it Inv. Math.} {\bf 1} (1966) 287--354.
\item{[Mu2]} Mumford, D., {\it Abelian varieties}, Oxford 
University Press 1974.
\item{[Mu3]} Mumford, D., {\it Tata Lectures on Theta I}, 
 Progress in Math. {\bf 28}, Birkh\"auser, 1983.
\item{[O'G]} O'Grady, K. G., ``On the Kodaira dimension of moduli
spaces of abelian surfaces'',
{\it Compositio Math.} {\bf 72} (1989), no. 2, 121--163.
\item{[Og]} Oguiso, K., ``On algebraic fiber space structures on a 
Calabi-Yau $3$-fold'',  With an appendix by Noboru Nakayama. 
{\it Internat. J. Math.} {\bf 4} (1993), no. 3, 439--465. 
\item{[RS]} Ranestad, K., Schreyer, F.-O., ``Varieties of sums of powers'',
{\it J. Reine Angew. Math.} (2000), to appear, preprint {\tt math.AG/9801110}.
\item{[Rod]} Rodland, E., ``The Pffafian Calabi-Yau, its mirror,
and their link to the Grassmannian $G(2,7)$,'' preprint 
{\tt math.AG/9801902}.
\item{[Sal]} Salmon, G., {\it Modern Higher Algebra} {\bf 4}, Hodges, Figgis,
and Co., Dublin 1885.
\item{[Scho]} Schoen, C., ``Algebraic cycles on certain desingularized nodal
hypersurfaces'', {\it Math. Ann.}, {\bf 270}, (1985), 17--27.
\item{[Scho2]} Schoen, C., ``On the geometry of a special determinantal
hypersurface associated to the Horrocks-Mumford vector bundle,''
{\it J. Reine Angew. Math.} {\bf 364} (1986), 85--111.
\item{[Schr]} Schreyer, F.-O., ``Geometry and algebra of prime Fano
3-folds of index 1 and genus 12'', preprint {\tt math.AG/9911044}.
\item{[SR]} Semple, G., Roth, L., {\it Algebraic Geometry}, Chelsea 1937.
\item{[Si]} Silberger, A. J., 
``An elementary construction of the representations of 
${\rm SL}(2,\,{\rm GF}(q))$'', {\it Osaka J. Math.} {\bf 6}, (1969), 329--338. 
\item{[Ta]} Tanaka, S., ``Construction and Classification of Irreducible
Representations of Special Linear Group of the Second Order over a Finite
Field,'' {\it Osaka J. Math.} {\bf 4} (1967), 65--84.
\item{[Tei]} Teissier, B., 
``The hunting of invariants in the geometry of discriminants'', in
{\it Real and complex singularities}, 
(Proc. Ninth Nordic Summer School/NAVF Sympos. Math., Oslo,
1976), 565--678. 
Sijthoff and Noordhoff, Alphen aan den Rijn, 1977. 
\item{[Tjo]} Tjotta, E., ``Quantum cohomology of a Pfaffian 
Calabi-Yau variety: verifying mirror symmetry predictions'',
preprint {\tt math.AG/9906119}
\item{[Ve]} V\'elu, J., ``Courbes elliptiques munies d'un sous-groupe
$\boldz/n\boldz\times\mu_n$'', {\it Bull. Soc. Math. France},
M\'emoire {\bf 57}, (1978) 5--152.
\item{[Ver]} Verra, A., 
``A short proof of the unirationality of ${\cal A}_5$'', 
{\it Nederl. Akad. Wetensch. Indag. Math.} {\bf 46} (1984), 
no. 3, 339--355. 
\item{[We]} Werner, J., ``Kleine Aufl\"osungen spezieller dreidimensionaler
Variet\"aten'', {\it Bonner Mathematische Schriften} {\bf 186}, {\it
Universit\"at Bonn, Mathematisches Institut,} Bonn, 1987.

\end